\documentclass[11pt]{article}

\usepackage[english]{babel}
\usepackage[affil-it]{authblk}
\usepackage{amssymb,amsmath,amsthm,mathtools}
\usepackage{enumerate}
\usepackage[margin=2.6cm]{geometry}
\usepackage{graphicx}
\graphicspath{{./figures}}
\usepackage{subfig}
\usepackage{booktabs}
\usepackage{multirow}

\usepackage[hidelinks]{hyperref}
\usepackage{xcolor}

\newtheorem{lemma}{Lemma}

\theoremstyle{definition}

\newtheorem{remark}{Remark}

\newcommand{\N}{\mathbb{N}}
\newcommand{\Z}{\mathbb{Z}}
\newcommand{\R}{\mathbb{R}}
\newcommand{\C}{\mathbb{C}}
\newcommand{\nulmat}{\mathrm{O}}
\newcommand{\abs}[1]{\left| #1 \right|}
\newcommand\norm[1]{\|#1\|}

\newcommand{\mi}[1]{{\boldsymbol{#1}}}  % multi-index
\newcommand{\iu}{{\mathrm{i}}}
\DeclareMathOperator*{\diag}{diag}

\title{Block triangular preconditioning for inverse source problems in time-space fractional diffusion equations}

\author[1]{Monoswini Majumdar}
\author[1,2]{Stefano Serra-Capizzano\thanks{s.serracapizzano@uninsubria.it}}
\author[1,3]{Rosita L. Sormani\thanks{rositaluisa.sormani@uninsubria.it}}

\affil[1]{Department of Science and High Technology, University of Insubria, Como, Italy}
\affil[2]{Division of Scientific Computing, Department of Information Technology, Uppsala University, Uppsala, Sweden}
\affil[3]{Department of Mathematics and Computer Science, University of Cagliari, Cagliari, Italy}

\date{}

\begin{document}

\maketitle

\begin{abstract}
    The current work investigates the effectiveness of block triangular preconditioners in accelerating and stabilizing the numerical solution of inverse source problems governed by time-space fractional diffusion equations (TSFDEs). We focus on the recovery of an unknown spatial source function in a multi-dimensional TSFDE, incorporating Caputo time-fractional derivatives and the fractional Laplacian. The inherent ill-posedness is addressed via a quasi-boundary value regularization, followed by a finite difference discretization that leads to large, structured linear systems. We develop and analyze a block triangular preconditioning strategy that mimics the coefficient matrix, while simplifying its structure for computational efficiency. Numerical experiments using the GMRES solver demonstrate that the proposed preconditioner significantly improve convergence rates, robustness, and accuracy, making it well-suited for large-scale, real-world inverse problems involving fractional modeling.
\end{abstract}

\noindent Keywords: Fractional diffusion equations, inverse problems, preconditioning, GMRES, quasi-boundary value regularization

\noindent MSC[2020]: 65F10 (65F08, 15A18, 15B05, 65R32)

\section{Introduction}

In the present study, we consider the inverse problem of identifying the unknown source term $f(\mathbf{x})$ of the following multi-dimensional time-space fractional diffusion equation
\begin{equation} \label{eq:original-fde}
\begin{dcases}
    \gamma_1(\mathbf{x},t) D_t^\beta u(\mathbf{x},t) + \gamma_2(\mathbf{x},t) (-\Delta)^{\frac{\omega}{2}} u(\mathbf{x},t) = f(\mathbf{x})q(t),
    & \mathbf{x} \in \Omega, \; t \in (0,T), \\
    u(\mathbf{x},t) = 0,
    & \mathbf{x} \in \Omega^\mathsf{c}, \; t \in (0,T), \\
    u(\mathbf{x},0) = \rho(\mathbf{x}),
    & \mathbf{x} \in \Omega, \\
    u(\mathbf{x},T) = \mu(\mathbf{x}),
    & \mathbf{x} \in \Omega,
\end{dcases}
\end{equation}
in which $T>0$ is the upper bound of the time interval and $\Omega\subset\R^d$, with $d$ a positive integer, is an open and bounded domain with piecewise smooth boundary, whose complement in $\R^d$ is denoted by $\Omega^\mathsf{c}$. In the main equation, $\gamma_1(\mathbf{x},t)$ and $\gamma_2(\mathbf{x},t)$ are continuous functions representing variable coefficients, $q(t)>0$ is a time-dependent known term, and $\beta\in (0,1)$ is the order of the left Caputo fractional derivative, which is defined as
\begin{equation*}
    D_t^\beta u(\mathbf{x},t) := \frac{1}{\Gamma(1-\beta)} \int_0^t  \frac{u_\tau(\mathbf{x},\tau)}{(t-\tau)^\beta}
    \,\mathrm{d}\tau,
\end{equation*}
where $\Gamma$ is the Euler Gamma function and $u_\tau$ indicates the partial derivative of $u(\mathbf{x},\tau)$ with respect to $\tau$. Moreover, the fractional Laplacian of order $ \omega\in (1,2)$ is defined as the hypersingular integral
\begin{equation} \label{eq:fractional-laplacian}
    (-\Delta)^\frac{\omega}{2} u(\mathbf{x},t) := c_\omega \, \text{P.V.} \int_{\R^d} \frac{u(\mathbf{x},t)-u(\mathbf{y},t)}{\norm{\mathbf{x}-\mathbf{y}}_2^{d+\omega}} \,\mathrm{d}\mathbf{y},
\end{equation}
where
\begin{equation*}
    c_\omega := \frac{2^\omega \Gamma\left( \frac{\omega+d}{2} \right)}{\pi^{d/2}\abs{\Gamma(-\omega/2)}}
\end{equation*}
is a normalization constant, $\norm{\mathbf{x}-\mathbf{y}}_2$ is the Euclidean distance between the points $\mathbf{x}$ and $\mathbf{y}$, and P.V. denotes the Cauchy principal value, which is necessary due to the singularity of the integrand at $\mathbf{y}=\mathbf{x}$. Finally, $\rho(\mathbf{x})$ and $\mu(\mathbf{x})$ are given functions representing the initial condition and the final state observation, respectively. The unknowns in this system are the solution $u(\mathbf{x},t)$ and the source function $f(\mathbf{x})$.

Time-Space Fractional Diffusion Equations (TSFDEs) have received significant attention in recent years, due to their effectiveness in modeling complex and anomalous diffusion processes across various scientific domains. TSFDEs incorporate fractional derivatives in both time and space, allowing them to overcome some limitations of classical diffusion equations and to capture long-range dependencies and memory effects observed in heterogeneous and disordered media.

On the computational front, the numerical discretization of TSFDEs often leads to large-scale, dense, and ill-conditioned linear systems, with high computational and memory demands that necessitate the development of efficient numerical solvers. Numerous high-accuracy discretization schemes have been developed to handle the intrinsic nonlocality of both Caputo and Riesz derivatives. In \cite{chen2018} it was shown that banded block triangular preconditioners, designed via Kronecker product splittings, are particularly efficient when the time-fractional order is near one. These preconditioners exploit the Toeplitz and Hessenberg structure of the discretized matrices and preserve their sparsity to an extent, enabling fast matrix-vector products. Further contributions in \cite{luo2022} introduced a two-step splitting iteration method and the corresponding preconditioner. Their theoretical analysis indicated that the preconditioned matrix can be expressed as a combination of an identity matrix, a low-rank matrix, and a small-norm matrix, effectively clustering the spectrum and improving the convergence of Krylov subspace methods such as the GMRES. These advancements have enabled TSFDEs to be effectively applied in fields such as geophysics, biomedical imaging, quantitative finance, and image processing.

One key computational challenge in this domain lies in solving inverse source problems, where the objective is to recover an unknown spatial source term from noisy final time observations. These problems are typically ill-posed, requiring the use of regularization techniques like the quasi-boundary value method to ensure stable solutions, particularly where source terms are recovered from noisy final-time data. To address these challenges, block preconditioning techniques have emerged as effective tools to enhance the convergence of iterative solvers, such as GMRES, used in the numerical solution of these systems. In \cite{ruan2025}, a modified quasi-boundary regularization method was proposed. The authors constructed a matrix-based formulation with a block Toeplitz structure and employed efficient preconditioned solvers to obtain stable numerical solutions under data noise.

This paper investigates the role of block triangular preconditioners in accelerating and stabilizing the numerical solution of inverse source problems governed by TSFDEs, focusing on the multi-dimensional setting with variable coefficients. In the slightly simpler constant coefficient setting, this problem was studied in \cite{pang2024}, using block-diagonalizable preconditioners that mimick the coefficient matrix structure, and in \cite{qiao25} through a variational approach.

Here, we build on the foundation established in \cite{pang2024} to develop and analyze a block triangular preconditioning strategy for efficiently solving the inverse source problem \eqref{eq:original-fde}, addressing the recovery of the unknown spatial source function $f(\mathbf{x})$ from noisy terminal data. This mathematical model is ill-posed, therefore it is regularized using the quasi-boundary value method. Then, we discretize the regularized forward problem using a finite difference scheme, which leads to a large block linear system containing multi-level Toeplitz structures. A block triangular preconditioner is proposed, designed to mimic the structure of the system matrix but eliminating upper off-diagonal blocks for computational efficiency. The GMRES method is used to solve the linear system and the proposed preconditioners are shown to accelerate convergence significantly, confirming that block triangular preconditioning is a robust and efficient approach for solving this type of inverse source problems.

The paper is structured as follows. After briefly introducing some notation and preliminary definitions in Section \ref{sec:prelim}, in Section \ref{sec:reg} we regularize the continuous ill-posed problem through the quasi-boundary value method and in Section \ref{sec:discr} we discretize the resulting well-posed equation via a finite difference scheme. In Section \ref{sec:precond}, we present our proposed preconditioners and in Section \ref{sec:numer} we perform several numerical experiments and simulations under varying conditions. The results validate the effectiveness of the block triangular preconditioners in improving eigenvalue clustering and reducing GMRES iterations. Finally, in Section \ref{sec:concl} we draw conclusions.

\section{Notation and preliminaries} \label{sec:prelim}

In this work, we deal with multi-level matrix structures of size $n_1n_2\ldots n_d$ with $d,n_1,\ldots,n_d\in\N$, meaning that the entries can be partitioned into blocks of size $n_1$, which can be further subdivided into smaller blocks of size $n_2$, and so on, until the inner level of size $n_d$ is reached. To effectively address the elements, we adopt the multi-index notation, with the following definitions and rules.

\paragraph{Multi-index notation.} A multi-index is a vector $(m_1,\ldots,m_d)\in\Z^d$. It is denoted with the corresponding bold cursive letter $\mi{m}$, while $\mi{0}$, $\mi{1}$, etc., denote vectors of all zeros, ones, and so on. The size will be clear from context. Operations and relations are performed componentwise, and $\mi{h},\ldots,\mi{k}$ represents the multi-index range $\{\mi{j}\in\Z^d:\mi{h}\leq\mi{j}\leq\mi{k}\}$: when we write $\mi{j}=\mi{h},\ldots,\mi{k}$ we mean that $\mi{j}$ varies in this range according to the standard lexicographic ordering. Moreover, the product of all the components of $\mi{m}$ is denoted by $N(\mi{m}):=m_1m_2\ldots m_d$. Using these conventions, any $\mathbf{x}\in\C^{N(\mi{n})}$ is expressed as $\mathbf{x} = \big[\, x_{\mi{i}} \,\big]_{\mi{i}=\mi{1}}^\mi{n}$, while a $d$-level matrix $A\in\C^{N(\mi{n})\times N(\mi{n})}$ takes the compact form $ A = \begin{bmatrix} \; a_{\mi{i}\mi{j}} \;\, \end{bmatrix}_{\mi{i},\mi{j}=\mi{1}}^\mi{n}$. The multi-index notation aligns with the structure of a multi-level matrix in a natural way, since $(i_1,j_1)$ represents the outer block structure, and so on, until $(i_d,j_d)$ indexes the inner level. For instance, in a 2-level matrix the element $(\mi{i},\mi{j})$ is found as the $(i_2,j_2)$ entry of the block in position $(i_1,j_1)$.

\paragraph{Multi-level Toeplitz matrices.} A particular type of multi-level matrix structure is found in multi-level Toeplitz matrices, whose element in position $(\mi{i},\mi{j})$ depends only on the difference $\mi{i}-\mi{j}$. They can be expressed as $\begin{bmatrix} t_{\mi{i}-\mi{j}} \end{bmatrix}_{\mi{i},\mi{j} = \mi{1}}^{\mi{n}} \;\in \C^{N(\mi{n})\times N(\mi{n})}$ for any fixed $\mi{n}\in\N^d$. In the case where $d=1$, it is simply a Toeplitz matrix with constant entries along the diagonals, while for $d=2$ it is a block Toeplitz matrix in which each block has a Toeplitz structure. More in general, a $d$-level Toeplitz matrix is a block Toeplitz matrix whose blocks have $(d-1)$-level Toeplitz structure. When the elements $t_\mi{k}$, $\mi{k}=\mi{1}-\mi{n},\ldots,\mi{n}-\mi{1}$, are the Fourier coefficients of a given complex-valued function $f\in L^1\big([-\pi,\pi]^d\big)$, the Toeplitz matrix is said to be generated by $f$.

\section{Regularization} \label{sec:reg}

We start by regularizing the ill-posed problem given in \eqref{eq:original-fde}, in order to obtain an equivalent well-posed formulation. We follow the quasi-boundary value method applied in \cite{pang2024}, also known as the non-local boundary value method. It is a classical regularization technique used for stabilizing ill-posed inverse problems, particularly when the terminal data is corrupted by noise. Unlike standard Tikhonov regularization, this approach modifies the boundary condition rather than the operator, offering better preservation of problem structure. It has been successfully applied to time-fractional diffusion equations in works such as \cite{ruan2018}, and it is well-suited for problems where the final state measurement is uncertain or incomplete.

Suppose that the measurement of the final time data $\mu(\mathbf{x}) = u(\mathbf{x},T)$ in Equation \eqref{eq:original-fde} is contaminated by noise of level $\varepsilon>0$ and denote the noise-contaminated measurement by $\mu_\varepsilon$. We assume that
\begin{equation} \label{eq:noise}
    \norm{\mu - \mu_\varepsilon}_{L^2} \leq \varepsilon,
\end{equation}
where $\norm{\cdot}_{L^2}$ represents the $L^2$ norm on $\Omega$. Applying the quasi-boundary value method to \eqref{eq:original-fde} leads to the following well-posed regularized problem
\begin{equation} \label{eq:regularized-fde}
    \begin{dcases}
        \gamma_1(\mathbf{x},t) D_t^\beta v(\mathbf{x},t) + \gamma_2(\mathbf{x},t) (-\Delta)^{\frac{\omega}{2}} v(\mathbf{x},t) = \Tilde{f}_{\lambda,\varepsilon}(\mathbf{x})q(t),
        & \mathbf{x} \in \Omega, \; t \in (0,T), \\
        v(\mathbf{x},t) = 0,
        &\mathbf{x} \in \Omega^\mathsf{c}, \; t \in (0,T), \\
        v(\mathbf{x},0) = \rho(\mathbf{x}),
        & \mathbf{x} \in \Omega, \\
        v(\mathbf{x},T) = \mu_\varepsilon(\mathbf{x}) - \lambda \gamma_2(\mathbf{x},t) (-\Delta)^{\frac{\omega}{2}} \Tilde{f}_{\lambda,\varepsilon}(\mathbf{x}),
        & \mathbf{x} \in \Omega,
    \end{dcases}
\end{equation}
in which $\lambda>0$ is the regularization parameter, $v(\mathbf{x},t)$ represents an approximation of $u(\mathbf{x},t)$, and $\Tilde{f}_{\lambda,\varepsilon}$ is an approximation of $f(\mathbf{x})$ depending on $\lambda$ and $\varepsilon$.

The choice of $\lambda$ directly affects the balance between fitting the observed data and maintaining a stable, accurate solution. If $\lambda$ is too small, the method focuses too much on matching the data exactly, which can lead to overfitting, capturing noise and producing unstable or oscillatory results. On the other hand, if $\lambda$ is too large, the method emphasizes stability too much and may overly smooth the solution, causing important features of the true source to be lost or misrepresented. Optimal selection can be guided by techniques such as the L-curve method or the discrepancy principle, depending on the noise level in the data.

\section{Discretization} \label{sec:discr}

In this section, we follow the finite difference method used in \cite{pang2024} to discretize \eqref{eq:regularized-fde}, generalizing it to the variable coefficients setting and assuming that the spatial domain $\Omega$ is an open hyperrectangle in $\R^d$. The analysis can be extended to non-Cartesian domains using immersion techniques; see for instance \cite{mariarosa24} for an application in the setting of a general bounded convex domain.

Let $S$ be a positive integer representing the number of time steps and define the corresponding grid of step size $\Delta t$ as
\begin{equation*}
    \Delta t := \frac{T}{S}, \qquad t_s := s\Delta t,
    \qquad s = 0,1,\ldots,S.
\end{equation*}
Recalling that $\beta\in (0,1)$, the L1 formula \cite{fde-diff-book} is adopted to approximate the Caputo fractional derivative:
\begin{equation} \label{eq:L1-formula} %tag{4}
    D_{L1}^\beta v(\mathbf{x},t_s)
    = \frac{\Delta t^{-\beta}}{\Gamma(2-\beta)}
    \left[ b^{(\beta)}_0 v(\mathbf{x},t_s) - \sum_{m=1}^{s-1} \left(b^{(\beta)}_{s-m-1} - b^{(\beta)}_{s-m} \right) v(\mathbf{x}, t_m) - b^{(\beta)}_{s-1} v(\mathbf{x},t_0) \right],
\end{equation}
in which for $m = 0,1,2, \dots$ the coefficients $b_m^{(\beta)}$ are given by
\begin{equation*}
    b^{(\beta)}_m = (m+1)^{1-\beta} - m^{1-\beta}.
\end{equation*}

Under mild hypothesis, an error estimate is available for this formula, see \cite[Lemma 2.1]{pang2024} or \cite[Theorem 1.6.1]{fde-diff-book}.

To deal with the fractional Laplacian \eqref{eq:fractional-laplacian}, a fractional centered difference scheme is employed. Assuming that $\Omega := \prod_{i=1}^d (a_i,b_i)$ and setting, for simplicity, a unique step size $h$ in all the $d$ spatial directions, let $n_1,n_2,\ldots,n_d\in\N$ be the number of step size in each direction. We define the points in $\R^d$
\begin{equation*}
    x^{(1)}_{j_1} := j_1 h,
    \qquad
    x^{(2)}_{j_2} := j_2 h,
    \qquad \cdots \qquad
    x^{(d)}_{j_d} := j_d h,
    \qquad \mi{j}\in\Z^d,
\end{equation*}
and, making use of the multi-indices $\mi{n}=(n_1,\ldots,n_d)$ and $\mi{j}=(j_1,\ldots,j_d)$, we denote the corresponding discrete mesh in $\Omega$ as
\begin{equation*}
    \Omega_h := \big\{ \big(x^{(1)}_{j_1}, x^{(2)}_{j_2}, \dots, x^{(d)}_{j_d}\big) \mid \mi{j}=\mi{1},\ldots,\mi{n} \big\}.
\end{equation*}
% and the related index space as
% \begin{equation*}
%     \phi_\mi{n} := \big\{ \mi{j}\in\N^d \mid \big(x^{(1)}_{j_1}, x^{(2)}_{j_2}, \dots, x^{(d)}_{j_d}\big) \in \Omega_h \big\},
%     \qquad \mi{n}=(n_1,\ldots,n_d).
% \end{equation*}
Then, the discrete fractional Laplacian is constructed as
\begin{equation} \label{eq:discrete-laplacian} %\tag{5}
    (-\Delta_h)^{\frac{\omega}{2}} v\big(x^{(1)}_{j_1}, \dots, x^{(d)}_{j_d},t\big) =
    \frac{1}{h^\omega} \sum_{\mi{l}\in\Z^d} a^{(\omega)}_\mi{l} v\big(x^{(1)}_{j_1}+ l_1 h, \,\dots\,, x^{(d)}_{j_d}+ l_d h,\, t\big),
\end{equation}
where $\mi{l}$ represents a multi-index $\mi{l} = (l_1,l_2,\ldots,l_d)$ and $a^{(\omega)}_\mi{l}$, $\mi{l}\in\Z^d$, are the Fourier coefficients of the $d$-variate function
\begin{equation} \label{eq:gen-func}
    g(\mi{\theta}) := \left[ \sum_{i=1}^d 4\sin^2{\left(\frac{\theta_i}{2}\right)} \right]^\frac{\omega}{2},
    \qquad \mi{\theta} = (\theta_1,\theta_2,\ldots,\theta_d),
\end{equation}
defined by the formula
\begin{equation} \label{eq:fourier-coeff} %tag{6}
    a^{(\omega)}_\mi{l} =
    \frac{1}{{(2\pi)}^d} \int_{{[-\pi,\pi]}^d}
    \left[\sum_{i=1}^d {4\sin^2{\left(\frac{\theta_i}{2}\right)}}\right]^\frac{\omega}{2}
    e^{-\iu\mi{l}\cdot\mi{\theta}}
    \,\mathrm{d}\mi{\theta},
    \qquad \iu^2 = -1,
    \quad \mi{l}\cdot\mi{\theta} = \sum_{i=1}^d l_i\theta_i.
\end{equation}

It was proved in \cite{diff-lapl} that the approximation \eqref{eq:discrete-laplacian} converges to the fractional Laplacian with order $\gamma\leq 2$, under mild conditions concerning the smoothness of $v(\mathbf{x},t)$.

To compute the Fourier coefficients, in the one-dimensional case where $d=1$ we have the explicit formula
\begin{equation} \label{eq:fourier-coeff-1d} %\tag{7}
    a_l^{(\omega)} = \frac{1}{2\pi} \int_{-\pi}^{\pi} \left[ 4\sin^2 \left( \frac{\theta}{2} \right) \right]^\frac{\omega}{2} e^{-\iu l \theta} \,\mathrm{d}\theta
    = \frac{{(-1)}^l \Gamma(\omega+1)}{\Gamma\left(\frac{\omega}{2} - l +1\right) \Gamma\left(\frac{\omega}{2} +l+1 \right)}.
\end{equation}
Moreover, in this case the following properties are satisfied.

\begin{lemma}[{\cite[Lemma 2.2]{pang2024}}]
    Let $\omega\in (1,2)$ and let $a^{(\omega)}_l$ be defined as in \eqref{eq:fourier-coeff-1d} for $l\in\Z$. Then,
    \begin{enumerate}[(a)]
        \item $a^{(\omega)}_0$ is the only nonnegative coefficient, i.e.,
            \begin{equation*}
                a^{(\omega)}_0 = \frac{\Gamma(\omega+1)}{\Gamma(\frac{\omega}{2} +1)^2} \geq 0,
                \qquad a^{(\omega)}_l = a^{(\omega)}_{-l} < 0,
                \quad \forall l\neq 0;
            \end{equation*}
        \item for any $l\in\N$, the following recursive expression holds
            \begin{equation*}
                a^{(\omega)}_{l+1} = \left(1 - \frac{\omega+1}{\frac{\omega}{2} + l +1} \right) a^{(\omega)}_l;
            \end{equation*}
        \item the Fourier coefficients sum to zero, i.e.,
            \begin{equation*}
                \sum_{l=-\infty}^{\infty} a^{(\omega)}_l = 0.
            \end{equation*}
    \end{enumerate}
\end{lemma}

On the other hand, in multi-dimensional settings the coefficients $a_\mi{l}^{(\omega)}$ can be computed numerically through the Fast Fourier Transform (FFT).

To construct the finite difference scheme, we define
\begin{equation*}
    \begin{dcases}
        e_0^{(\beta)} = b_0^{(\beta)}, \\
        e_m^{(\beta)} = b_m^{(\beta)} - b_{m-1}^{(\beta)} = (m+1)^{1-\beta} - 2m^{1-\beta} + (m - 1)^{1-\beta},   &\qquad 1 \leq m \leq S - 1;
    \end{dcases}
\end{equation*}
then, combining \eqref{eq:regularized-fde}, \eqref{eq:L1-formula} and \eqref{eq:discrete-laplacian} and omitting the truncation errors, we obtain, for $s=1,\ldots,S$,
\begin{equation} \label{eq:numer-scheme} %\tag{8}
    \begin{dcases}
        \frac{\Delta t^{-\beta}\gamma_{1,\mi{j}}^{(s)}}{\Gamma(2-\beta)} \left[ e^{(\beta)}_0 v_\mi{j}^{(s)} + \sum_{m=1}^{s-1} e^{(\beta)}_{s-m} v_\mi{j}^{(m)} - b^{(\beta)}_{s-1} v_\mi{j}^{(0)}  \right] + \frac{\gamma_{2,\mi{j}}^{(s)}}{h^\omega} \sum_{\mi{l} \in \Z^d} a^{(\omega)}_\mi{l} v_{\mi{j}+\mi{l}}^{(s)} = \Tilde{f}_\mi{j} q^{(s)},
        & \mi{j}=\mi{1},\ldots,\mi{n},
        \\
        v_\mi{j}^{(s)} = 0,
        & \mi{j}\in\Z^d\setminus[\mi{1},\ldots,\mi{n}],
        \\
        v_\mi{j}^{(0)}= \rho_\mi{j},
        & \mi{j}=\mi{1},\ldots,\mi{n},
        \\
        v_\mi{j}^{(S)} = \mu_{\varepsilon,\mi{j}} - \frac{\lambda\gamma_{2,\mi{j}}^{(S)}}{h^\omega} \sum_{\mi{l}\in\Z^d} a^{(\omega)}_\mi{l} \Tilde{f}_{\mi{j}+\mi{l}},
        & \mi{j}=\mi{1},\ldots,\mi{n},
    \end{dcases}
\end{equation}
where $v_\mi{j}^{(s)}$ denotes the numerical approximation of $v\big(x^{(1)}_{j_1}, x^{(2)}_{j_2}, \dots, x^{(d)}_{j_d}, t_s\big)$ and similarly
\begin{alignat*}{2}
    \Tilde{f}_\mi{j} &:= \Tilde{f}_{\lambda,\varepsilon}\big(x^{(1)}_{j_1}, x^{(2)}_{j_2}, \dots, x^{(d)}_{j_d}\big),
    &\qquad
    q^{(s)} &:= q(t_s),
    \\
    \gamma_{1,\mi{j}}^{(s)} &:= \gamma_1\big(x^{(1)}_{j_1}, x^{(2)}_{j_2}, \dots, x^{(d)}_{j_d}, t_s\big),
    &\gamma_{2,\mi{j}}^{(s)} &:= \gamma_2\big(x^{(1)}_{j_1}, x^{(2)}_{j_2}, \dots, x^{(d)}_{j_d}, t_s\big),
    \\
    \rho_\mi{j} &:= \rho\big(x^{(1)}_{j_1}, x^{(2)}_{j_2}, \dots, x^{(d)}_{j_d}\big),
    &\mu_{\varepsilon,\mi{j}} &:= \mu_\varepsilon\big(x^{(1)}_{j_1}, x^{(2)}_{j_2}, \dots, x^{(d)}_{j_d}\big).
\end{alignat*}

\section{Matrix formulation and preconditioning} \label{sec:precond}

Now, let us write the scheme found above in matrix form. We define the coefficient
\begin{equation*}
    \eta = \eta(\beta,\Delta t)
    := \Delta t^\beta\Gamma(2-\beta),
\end{equation*}
the vectors
\begin{equation*}
    \mathbf{v}^{(s)} := \begin{bmatrix}
        v_\mi{j}^{(s)}
    \end{bmatrix}_{\mi{j} = \mi{1},\ldots,\mi{n}}^{T},
        \qquad
    \tilde{\mathbf{f}} := \begin{bmatrix}
        \Tilde{f}_\mi{j}
    \end{bmatrix}_{\mi{j} = \mi{1},\ldots,\mi{n}}^{T},
    \qquad
    \boldsymbol{\rho} := \begin{bmatrix}
        \rho_\mi{j}
    \end{bmatrix}_{\mi{j} = \mi{1},\ldots,\mi{n}}^{T},
        \qquad
    \boldsymbol{\mu}_\varepsilon := \begin{bmatrix}
        \mu_{\varepsilon,\mi{j}}
    \end{bmatrix}_{\mi{j} = \mi{1},\ldots,\mi{n}}^{T},
\end{equation*}
and the diagonal matrices
\begin{equation*}
    D_1^{(s)} := \diag_{\mi{j}=\mi{1},\ldots,\mi{n}} \gamma_{1,\mi{j}}^{(s)},
        \qquad
    D_2^{(s)} := \diag_{\mi{j}=\mi{1},\ldots,\mi{n}} \gamma_{2,\mi{j}}^{(s)},
        \qquad
    s=1,\ldots,S.
\end{equation*}
Therefore, the numerical scheme \eqref{eq:numer-scheme} can be reformulated as follows:
\begin{equation*} %tag{9}
    \begin{dcases}
        D_1^{(s)} \left( e_0^{(\beta)} \mathbf{v}^{(s)} + \sum_{m=1}^{s-1} e_{s-m}^{(\beta)} \mathbf{v}^{(m)} \right)
        + \frac{\eta}{h^\omega} D_2^{(s)} B_\mi{n} \mathbf{v}^{(s)} - \eta q^{(s)} \tilde{\mathbf{f}}
        = b_{s-1}^{(\beta)} D_1^{(s)}\boldsymbol{\rho},
        & s=1,\ldots,S,
        \\
        \mathbf{v}^{(S)} + \frac{\lambda}{ h^\omega} D_2^{(S)} B_\mi{n} \tilde{\mathbf{f}} = \boldsymbol{\mu}_\varepsilon,
    \end{dcases}
\end{equation*}
where $B_\mi{n}$ is the $d$-level Toeplitz matrix of size $N(\mi{n})\times N(\mi{n})$ representing the fractional Laplacian, with elements $a_\mi{l}^{(\omega)}$ given by \eqref{eq:fourier-coeff}. In other words, $B_\mi{n}$ is the $d$-level Toeplitz matrix generated by the function $g(\boldsymbol{\theta})$ in $\eqref{eq:gen-func}$. Note that, since $g(\boldsymbol{\theta})$ is real-valued and even, $B_\mi{n}$ is a real symmetric matrix.

The unknowns in the scheme above are $\tilde{\mathbf{f}}$ and $\mathbf{v}^{(s)}$ for $s=1,\ldots,S$. By gathering the equations for all the time levels and the regularization equation in one large linear system, we get
\begin{equation} \label{eq:lin-syst}
    A_{\mi{n},S} \mathbf{y} = \mathbf{z}
\end{equation}
where
\begin{equation*}
    \renewcommand*{\arraystretch}{1.4}
    A_{\mi{n},S} := \left[ \begin{array}{cccc|c}
        e_0 D_1^{(1)} + \frac{\eta}{h^\omega} D_2^{(1)} B_\mi{n} & \nulmat_\mi{n} & \dots & \nulmat_\mi{n} & -\eta q^{(1)} I_\mi{n} \\
        e_1 D_1^{(2)} & e_0 D_1^{(2)} + \frac{\eta}{h^\omega} D_2^{(2)} B_\mi{n} & \ddots & \vdots & -\eta q^{(2)} I_\mi{n} \\
        \vdots & \ddots & \ddots & \nulmat_\mi{n} & \vdots \\
        e_{S-1} D_1^{(S)} & \cdots & e_1 D_1^{(S)} & e_0 D_1^{(S)} + \frac{\eta}{h^\omega} D_2^{(S)} B_\mi{n} & -\eta q^{(S)} I_\mi{n} \\
            \hline
        \nulmat_\mi{n} & \cdots & \nulmat_\mi{n} & I_\mi{n} & \frac{\lambda}{h^\omega}  D_2^{(S)} B_\mi{n}
    \end{array} \right],
\end{equation*}
in which $I_\mi{n}$ is the identity matrix of size $N(\mi{n})$ and $\nulmat_\mi{n}$ is the zero matrix, while
\begin{equation*}
    \renewcommand*{\arraystretch}{1.4}
    \mathbf{y} := \left[ \begin{array}{c}
        \mathbf{v}^{(1)} \\ \mathbf{v}^{(2)} \\ \vdots \\ \mathbf{v}^{(S)} \\ \hline \tilde{\mathbf{f}}
    \end{array} \right],
    \qquad
    \mathbf{z}:= \left[ \begin{array}{c}
        b_0 D_1^{(1)} \boldsymbol{\rho} \\ b_1 D_1^{(2)} \boldsymbol{\rho} \\ \vdots \\ b_{S-1} D_1^{(S)} \boldsymbol{\rho} \\ \hline \boldsymbol{\mu}_\varepsilon
    \end{array} \right].
\end{equation*}

\paragraph{Implementation.} To solve the linear system, we employ the MATLAB built-in GMRES function, to which we supply a custom routine for the computation of the matrix-vector product $\mathbf{z} = A_{\mi{n},S} \mathbf{y}$, where
\begin{equation*}
    \mathbf{y} = \begin{bmatrix}
        \mathbf{y}^{(1)} \\ \vdots \\ \mathbf{y}^{(S+1)}
    \end{bmatrix},
    \qquad
    \mathbf{z} = \begin{bmatrix}
        \mathbf{z}^{(1)} \\ \vdots \\ \mathbf{z}^{(S+1)}
    \end{bmatrix},
    \qquad \mathbf{y}^{(s)},\mathbf{z}^{(s)}\in\C^{\mi{n}}.
\end{equation*}
The algorithm we implement follows this scheme:
\begin{enumerate}
    \item For $s=1,\ldots,S$
    \begin{enumerate}
        \item Compute $\mathbf{w}^{(s)} := D_1^{(s)} \big(e_{s-1} \mathbf{y}^{(1)} + \ldots + e_0 \mathbf{y}^{(s)}\big)$ in $O\big(s N(\mi{n})\big)$ operations;
        \item Compute $\frac{\eta}{h^\omega}D_2^{(s)}B_\mi{n}\mathbf{y}^{(s)}$ in $O\big(N(\mi{n})\log N(\mi{n})\big)$ operations without explicitly constructing the matrix, by embedding $B_\mi{n}$ in a larger $d$-level circulant matrix and exploiting the Fast Fourier Transform (FFT). Then, add the result to $\mathbf{w}^{(s)}$ in $O\big(N(\mi{n})\big)$ operations;
        \item Subtract $\eta q^{(s)} \mathbf{y}^{(s)}$ in $O\big(N(\mi{n})\big)$ operations to obtain $\mathbf{z}^{(s)}$.
    \end{enumerate}
    \item For $s=S+1$
    \begin{enumerate}
        \item Compute $\frac{\lambda}{h^\omega}D_2^{(S)}B_\mi{n}\mathbf{y}^{(S+1)}$ in $O\big(N(\mi{n})\log N(\mi{n})\big)$ operations using the FFT;
        \item Add $\mathbf{y}^{(S)}$ in $O\big(N(\mi{n})\big)$ operations to obtain $\mathbf{z}^{(S+1)}$.
    \end{enumerate}
\end{enumerate}
The overall computational complexity of this algorithm is $O\big( SN(\mi{n})\log SN(\mi{n})\big)$.

\subsection{Preconditioning proposal}

We now turn to our preconditioning strategy, aiming to reduce the condition number of the coefficient matrix or to cluster the spectrum at 1. Based on the structure of $A_{\mi{n},S}$, we propose the following lower block triangular preconditioner to ensure computational efficiency:
\begin{equation*}
    \renewcommand*{\arraystretch}{1.4}
    P_{\mi{n},S} := \left[ \begin{array}{cccc|c}
        e_0 D_1^{(1)} + \frac{\eta}{h^\omega} D_2^{(1)} B_\mi{n} & \nulmat_\mi{n} & \dots & \nulmat_\mi{n} & \nulmat_\mi{n} \\
        e_1 D_1^{(2)} & e_0 D_1^{(2)} + \frac{\eta}{h^\omega} D_2^{(2)} B_\mi{n} & \ddots & \vdots & \vdots \\
        \vdots & \ddots & \ddots & \nulmat_\mi{n} & \vdots \\
        e_{S-1} D_1^{(S)} & \cdots & e_1 D_1^{(S)} & e_0 D_1^{(S)} + \frac{\eta}{h^\omega} D_2^{(S)} B_\mi{n} & \nulmat_\mi{n} \\
            \hline
        \nulmat_\mi{n} & \cdots & \nulmat_\mi{n} & I_\mi{n} & \frac{\lambda}{h^\omega}  D_2^{(S)} B_\mi{n}
    \end{array} \right],
\end{equation*}
where we simply replaced the elements in the upper right section of $A_{\mi{n},S}$ with zeros.

$P_{\mi{n},S}$ clearly offers an accurate spectral approximation of the coefficient matrix, since it holds $A_{\mi{n},S} = P_{\mi{n},S} + R_{\mi{n},S}$, where $R_{\mi{n},S}$ is a matrix of rank at most $N(\mi{n})$. Therefore,
\begin{equation*}
    P_{\mi{n},S}^{-1}A_{\mi{n},S} = I_{\mi{n},S} + P_{\mi{n},S}^{-1}R_{\mi{n},S},
\end{equation*}
where $I_{\mi{n},S}$ is the identity of size $(S+1)N(\mi{n})$ and $P_{\mi{n},S}^{-1}R_{\mi{n},S}$ has rank at most $N(\mi{n})$. In other words, the preconditioned matrix is a low-rank correction of the identity, hence the eigenvalues will be clustered at 1.

\paragraph{Implementation.} The block triangular structure allows us to implement a forward substitution algorithm to solve linear systems of the form $P_{\mi{n},S} \mathbf{y} = \mathbf{z}$, with
\begin{equation*}
    \mathbf{y} = \begin{bmatrix}
        \mathbf{y}^{(1)} \\ \vdots \\ \mathbf{y}^{(S+1)}
    \end{bmatrix},
    \qquad
    \mathbf{z} = \begin{bmatrix}
        \mathbf{z}^{(1)} \\ \vdots \\ \mathbf{z}^{(S+1)}
    \end{bmatrix},
    \qquad \mathbf{y}^{(s)},\mathbf{z}^{(s)}\in\C^{\mi{n}}.
\end{equation*}
The algorithm is structured as follows:
\begin{enumerate}
    \item For $s=1$
    \begin{enumerate}
        \item Compute $e_0 D_1^{(1)} + \frac{\eta}{h^\omega} D_2^{(1)} B_\mi{n}$ in $O\big(N(\mi{n})^2\big)$ flops;
        \item Solve the associated linear system with right hand term $\mathbf{z}^{(1)}$, using the MATLAB backslash operator, to obtain $\mathbf{y}^{(1)}$ in $O\big(N(\mi{n})^3\big)$ flops.
    \end{enumerate}
    \item For $s=2,\ldots,S$
    \begin{enumerate}
        \item Compute $e_0 D_1^{(s)} + \frac{\eta}{h^\omega} D_2^{(s)} B_\mi{n}$ in $O\big(N(\mi{n})^2\big)$ flops;
        \item Compute the right hand term $\mathbf{z}^{(s)} - D_1^{(s)} \big( e_{s-1} \mathbf{y}^{(1)} + \ldots + e_1 \mathbf{y}^{(s-1)} \big)$ in $O\big(s N(\mi{n})\big)$ flops;
        \item Solve the associated linear system with the backslash operator to obtain $\mathbf{y}^{(s)}$ in $O\big(N(\mi{n})^3\big)$ flops.
    \end{enumerate}
    \item For $s=S+1$
    \begin{enumerate}
        \item Compute $\frac{\lambda}{h^\omega} D_2^{(S)} B_\mi{n}$ in $O\big(N(\mi{n})^2\big)$ flops;
        \item Compute the right hand term $\mathbf{z}^{(S+1)} - \mathbf{y}^{(S)}$ in $N(\mi{n})$ flops;
        \item Solve the associated linear system with the backslash operator to obtain $\mathbf{y}^{(S+1)}$ in $O\big(N(\mi{n})^3\big)$ flops.
    \end{enumerate}
\end{enumerate}
Note that both the coefficient matrix and the preconditioner need not be explicitly assembled, allowing us to bypass the construction phase, except for the Toeplitz matrix $B_\mi{n}$ which can be constructed just once during the setup.

\begin{remark}
    The most computationally demanding steps in the above procedure are those involving the solution of linear systems that contain the dense Toeplitz matrix $B_\mi{n}$. In the constant coefficient case, i.e., when $D_1^{(s)} = \gamma_1 I_\mi{n}$ and $D_2^{(s)} = \gamma_2 I_\mi{n}$ for all the time instants $s=1,\ldots,S$, we can overcome this difficulty by substituting $B_\mi{n}$ with structured matrices that spectrally approximate $B_\mi{n}$ and allow for more efficient inversion through fast transforms. A detailed analysis of this constant coefficient subcase is provided in \cite{pang2024}, where the authors propose preconditioners that preserve the structure of the coefficient matrix by replacing every occurrence of $B_\mi{n}$ in $A_{\mi{n},S}$ with its $\tau$ and Strang preconditioners \cite{Bini-tau}.

    This approach can be directly extended to our block triangular preconditioning strategy by substituting each instance of $B_\mi{n}$ in the definition of $P_{\mi{n},S}$ with its $\tau$ and Strang preconditioners. Let us denote the preconditioners constructed in this way as $T_{\mi{n},S}$ and $S_{\mi{n},S}$, respectively. This procedure allows one to solve the aforementioned linear systems by leveraging the diagonalizability of circulant and $\tau$ matrices, respectively via the Fourier and Sine matrices, and the associated fast transforms \cite{vanloan-book}. Clearly, $P_{\mi{n},S}$ provides a more accurate spectral approximation of the coefficient matrix than $T_{\mi{n},S}$ and $S_{\mi{n},S}$, since the preconditioned matrix is a low-rank correction of the identity. As a result, it is expected that preconditioning with $P_{\mi{n},S}$ would lead to a lower number of iterations when using the GMRES method. However, thanks to the computational efficiency offered by fast transforms, $T_{\mi{n},S}$ and $S_{\mi{n},S}$ would reasonably yield lower CPU times than $P_{\mi{n},S}$, making them the most competitive. In particular, $T_{\mi{n},S}$ is expected to deliver a superior performance, as the $\tau$ preconditioner more accurately captures the behavior of the small eigenvalues of $B_\mi{n}$ compared to the circulant Strang preconditioner. Given that this subcase has already been extensively studied in \cite{pang2024}, with minor differences with respect to our setting, we will not delve further into it.

    Unfortunately, the computational benefits of the Strang and $\tau$ preconditioners cannot be exploited in the variable coefficients case. In fact, the circulant or $\tau$ structure is lost when they are multiplied by $D_2^{(s)}$. Hence, we restrict our attention to the preconditioner $P_{\mi{n},S}$, since it provides the most accurate spectral approximation of $A_{\mi{n},S}$.
\end{remark}

\section{Numerical experiments} \label{sec:numer}

We conclude our analysis with a set of numerical experiments designed to evaluate the efficiency of the proposed strategy. Our primary objective is to assess how effectively the preconditioner reduces the number of iterations required by the GMRES method. In support of this goal, we first conduct a spectral analysis of both the nonpreconditioned and preconditioned coefficient matrix, to visualize the preconditioner’s ability to cluster eigenvalues around 1 and improve the condition number. In addition, in the one dimensional case, we test the method’s accuracy in reconstructing the source term. All tests were run on MATLAB R2018b.

\subsection{One dimensional setting} \label{sssec:numer-1d}

Let us consider the one-dimensional version of the time-space fractional diffusion equation \eqref{eq:regularized-fde}, with
\begin{gather*}
    \Omega=(0,\pi), \qquad T=1, \\
    \gamma_1(x,t) = t e^x, \qquad \gamma_2(x,t) = t^2 x, \\
    q(t) = t^2, \qquad \rho(x)=0.
\end{gather*}
To construct the final time data $\mu(x)=v(x,T)$, we fix the exact space source term
\begin{equation*}
    f(x) = x\sin(x)
\end{equation*}
and numerically solve the direct problem
\begin{equation*}
\begin{dcases}
    \gamma_1(x,t) D_t^\beta u(x,t) + \gamma_2(x,t) (-\Delta)^{\frac{\omega}{2}} u(x,t) = f(x)q(t),
    & x \in \Omega, \; t \in (0,T), \\
    u(x,t) = 0,
    & x \in \Omega^\mathsf{c}, \; t \in (0,T), \\
    u(x,0) = \rho(x),
    & x \in \Omega.
\end{dcases}
\end{equation*}
In other words, we solve the linear system obtained by applying the same discretization technique used in Section \ref{sec:discr}, which is
\begin{equation*}
    \widehat{A}_{\mi{n},S} \mathbf{y} = \widehat{\mathbf{z}}
\end{equation*}
where
\begin{equation*}
    \renewcommand*{\arraystretch}{1.4}
    \widehat{A}_{\mi{n},S} := \left[ \begin{array}{cccc}
        e_0 D_1^{(1)} + \frac{\eta}{h^\omega} D_2^{(1)} B_\mi{n} & \nulmat_\mi{n} & \dots & \nulmat_\mi{n} \\
        e_1 D_1^{(2)} & e_0 D_1^{(2)} + \frac{\eta}{h^\omega} D_2^{(2)} B_\mi{n} & \ddots & \vdots \\
        \vdots & \ddots & \ddots & \nulmat_\mi{n} \\
        e_{S-1} D_1^{(S)} & \cdots & e_1 D_1^{(S)} & e_0 D_1^{(S)} + \frac{\eta}{h^\omega} D_2^{(S)} B_\mi{n}
    \end{array} \right],
\end{equation*}
\begin{equation*}
    \renewcommand*{\arraystretch}{1.4}
    \mathbf{y} := \left[ \begin{array}{c}
        \mathbf{v}^{(1)} \\ \mathbf{v}^{(2)} \\ \vdots \\ \mathbf{v}^{(S)}
    \end{array} \right],
    \qquad
    \widehat{\mathbf{z}}:= \left[ \begin{array}{c}
        b_0 D_1^{(1)} \boldsymbol{\rho} + \eta q^{(1)}\mathbf{f} \\ b_1 D_1^{(2)} \boldsymbol{\rho} + \eta q^{(2)} \mathbf{f} \\ \vdots \\ b_{S-1} D_1^{(S)}\boldsymbol{\rho} + \eta q^{(S)} \mathbf{f}
    \end{array} \right],
\end{equation*}
where $\mathbf{f}$ denotes the sampling of the source function over the spatial mesh.

Then, to simulate the contamination of $\mu(x)$ by noise as in \eqref{eq:noise}, we perturb it by adding a random quantity $\varepsilon\delta (x)$, where $\delta (x)$ is randomly selected in the interval $(-1,1)$ and the noise level is set to $\varepsilon = 0.01\cdot\norm{\mu}_2$. Therefore, we have $\mu_\varepsilon(x) = \mu(x) + \varepsilon\delta (x)$ and our regularized model problem takes the form
\begin{equation*}
    \begin{dcases}
        \gamma_1(x,t) D_t^\beta v(x,t) + \gamma_2(x,t) (-\Delta)^{\frac{\omega}{2}} v(x,t) = t^2 \tilde{f}_{\lambda,\varepsilon}(x),
        & x\in (0,\pi), \; t \in (0,1), \\
        v(x,t) = 0,
        & x\in\R\setminus (0,\pi), \; t \in (0,1), \\
        v(x,0) = 0,
        & x\in (0,\pi), \\
        v(x,1) = \mu_\varepsilon(x) - \gamma_2(x,t) \lambda(-\Delta)^{\frac{\omega}{2}} \tilde{f}_{\lambda,\varepsilon}(x),
        & x\in (0,\pi).
    \end{dcases}
\end{equation*}

\subsubsection{Spectral analysis}

First, let us examine the effectiveness of the preconditioner $P_{n,S}$ in reducing the condition number and clustering the eigenvalues of the coefficient matrix at 1. In what follows, the regularization parameter is set to $\lambda= 5\cdot10^{-3}$.

\begin{figure}[ht]
    \centering
    \subfloat[][$\beta=0.1$, $\omega=1.9$]
    {\includegraphics[width=0.33\textwidth]{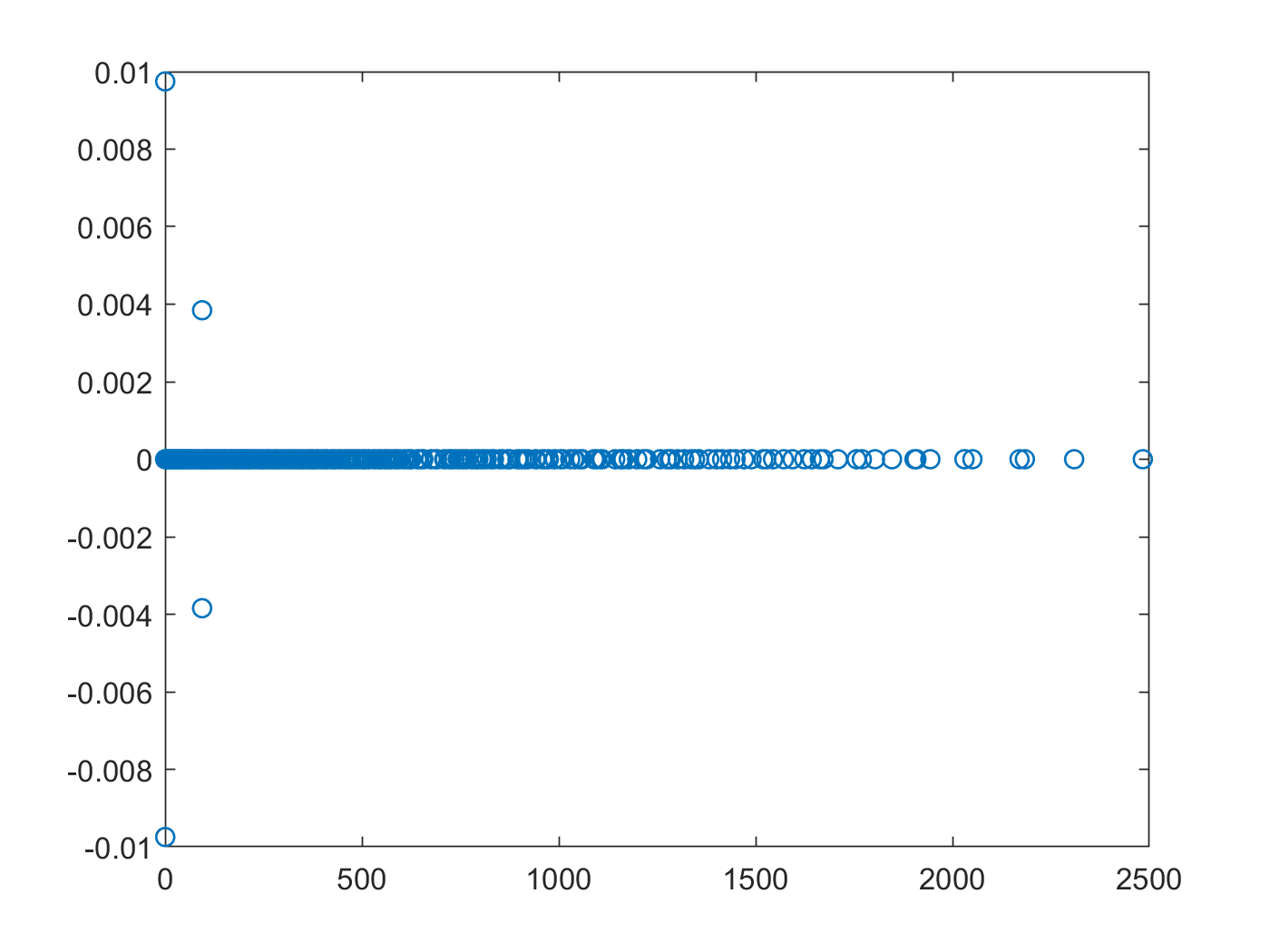}}
    \hfill
    \subfloat[][$\beta=0.1$, $\omega=1.9$]
    {\includegraphics[width=0.33\textwidth]{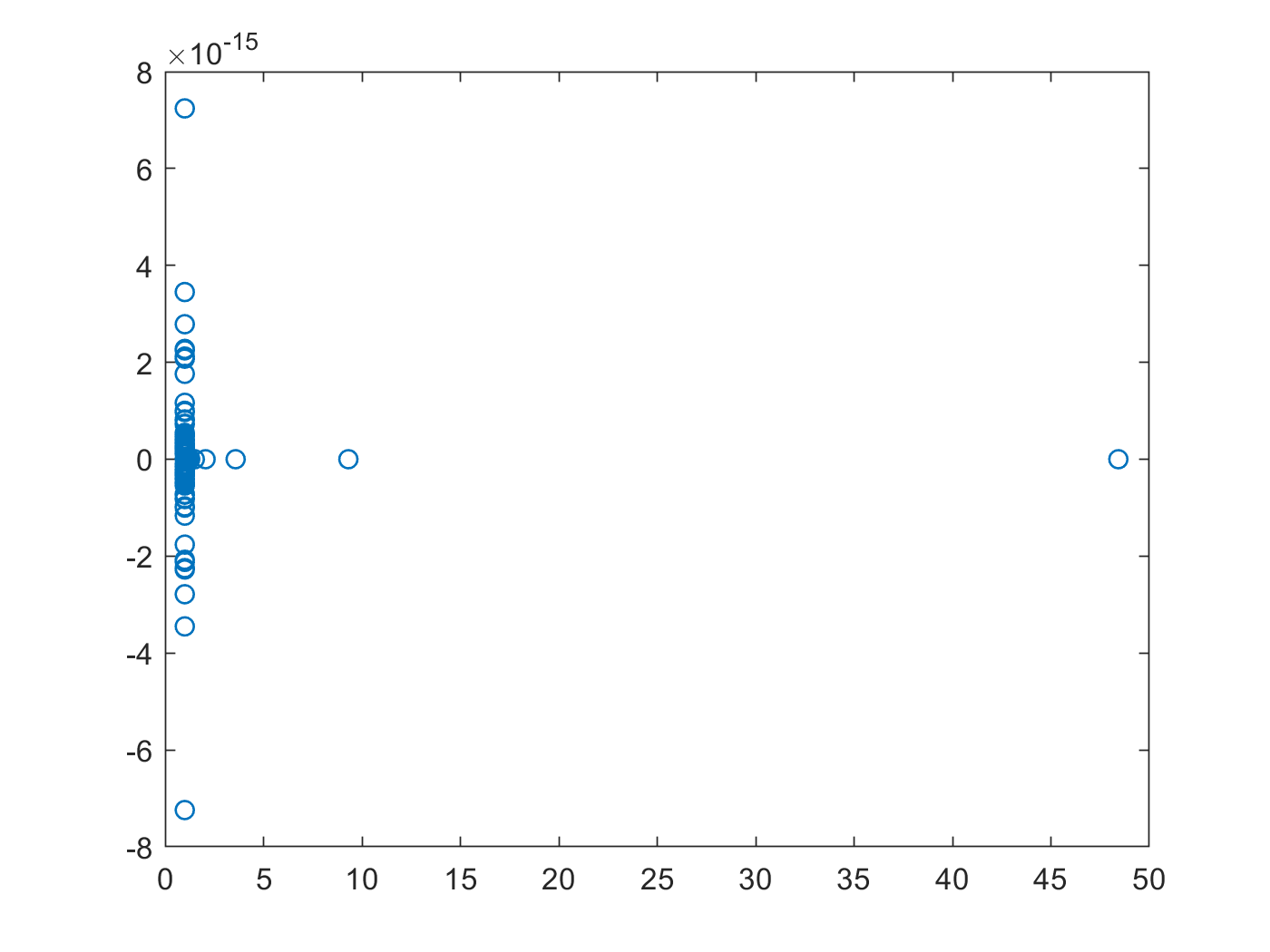}}
    \hfill
    \subfloat[][$\beta=0.1$, $\omega=1.9$]
    {\includegraphics[width=0.33\textwidth]{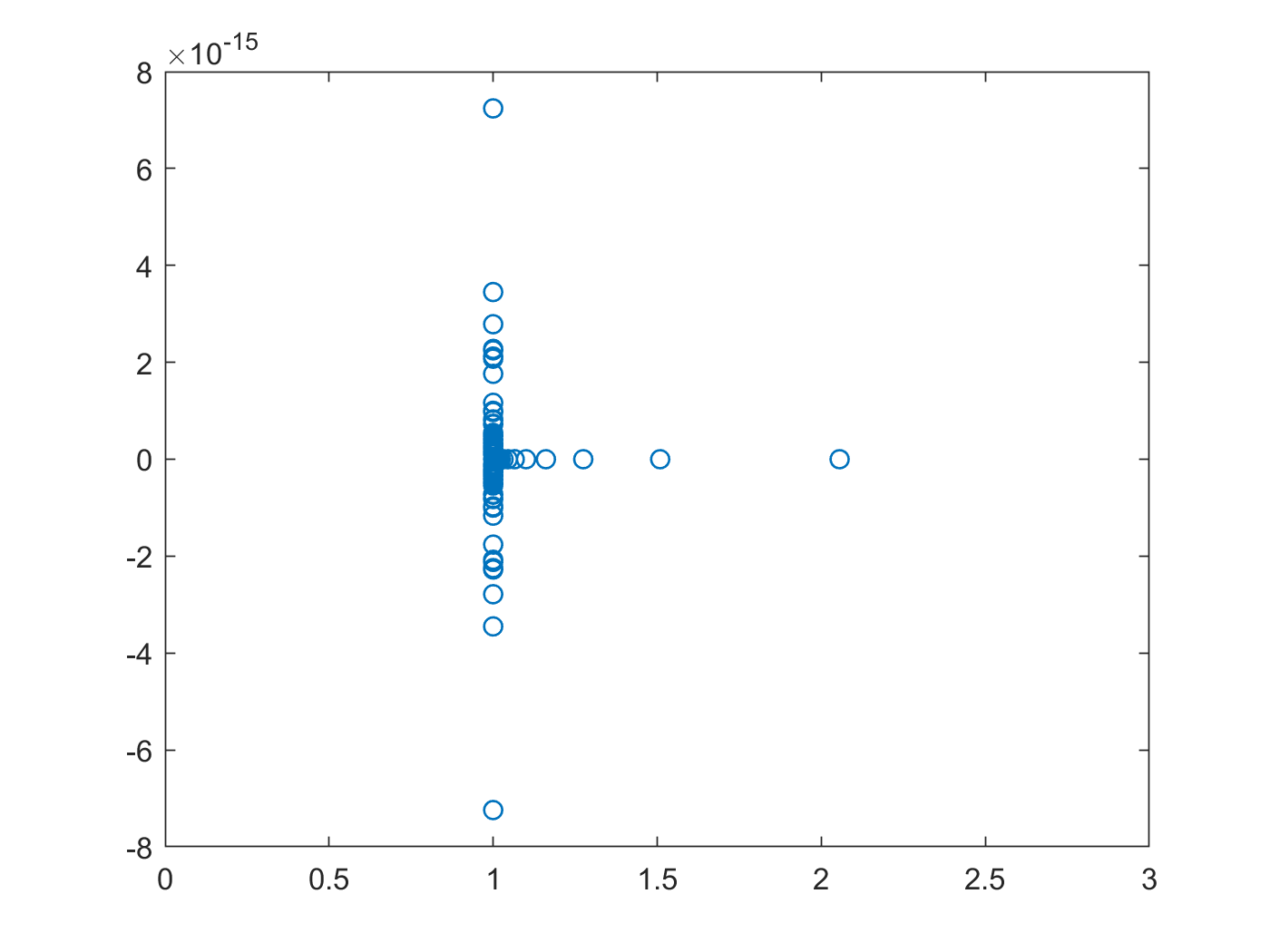}}
    \\
    \subfloat[][$\beta=0.5$, $\omega=1.5$]
    {\includegraphics[width=0.33\textwidth]{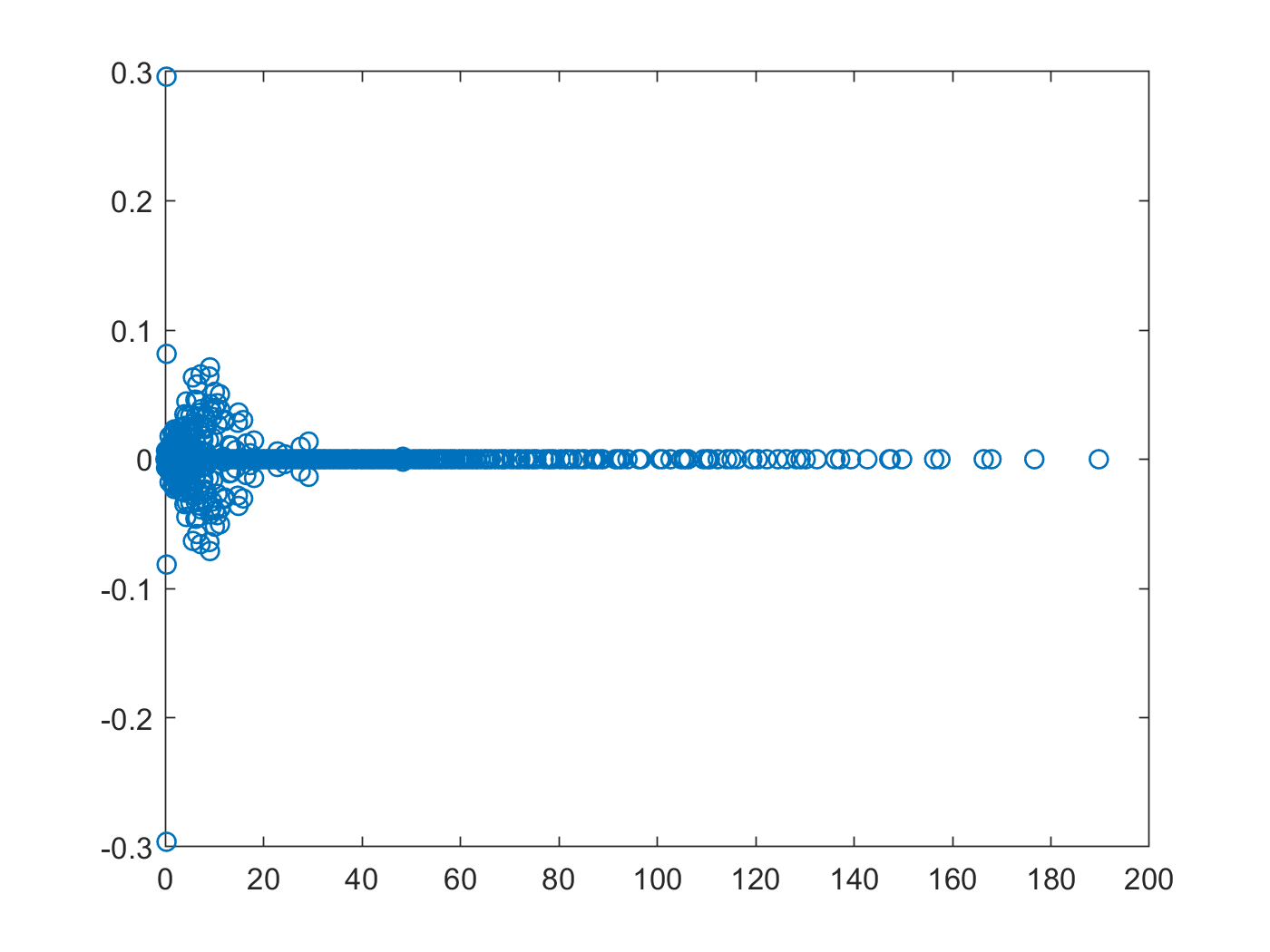}}
    \hfill
    \subfloat[][$\beta=0.5$, $\omega=1.5$]
    {\includegraphics[width=0.33\textwidth]{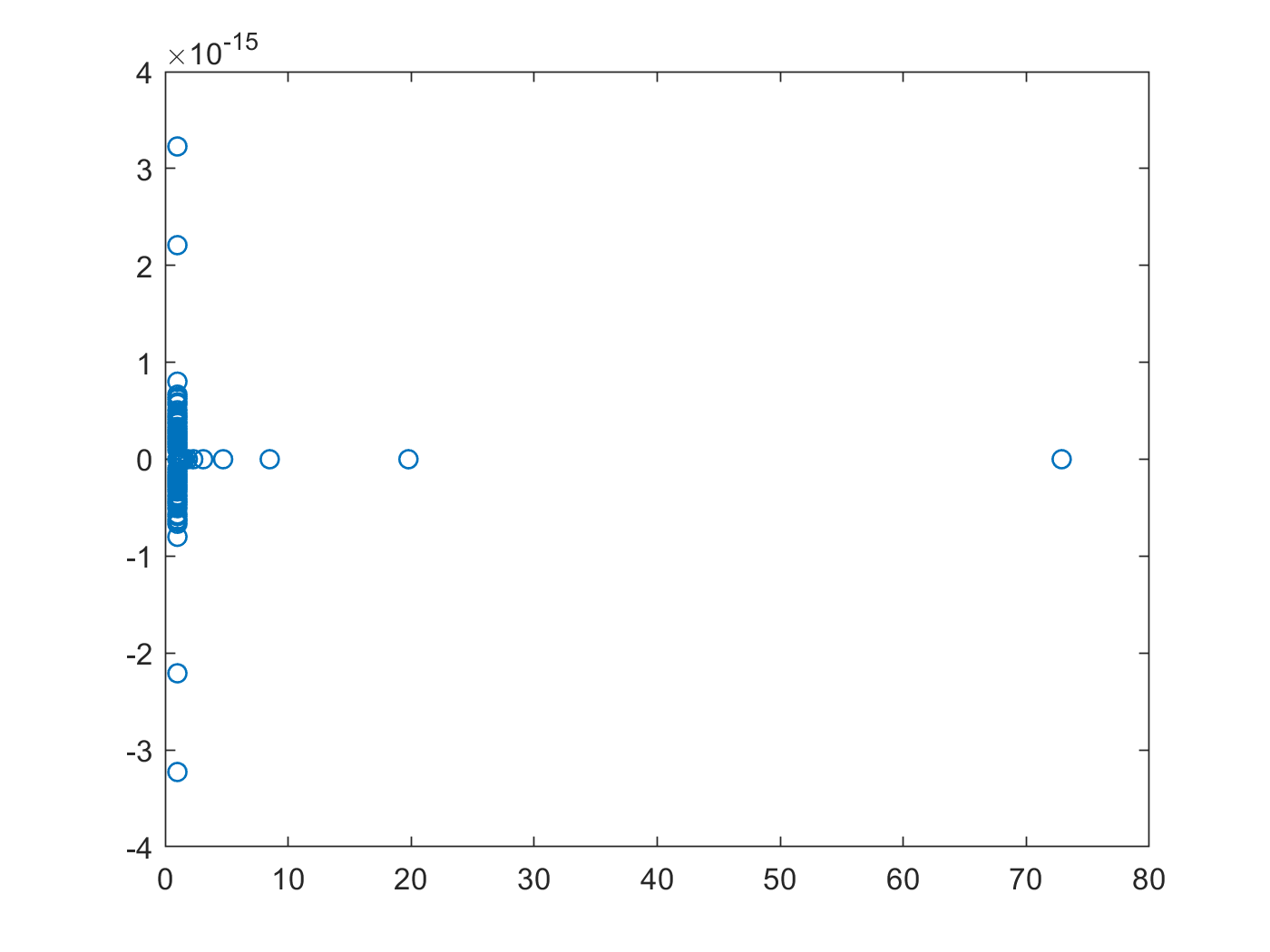}}
    \hfill
    \subfloat[][$\beta=0.5$, $\omega=1.5$]
    {\includegraphics[width=0.33\textwidth]{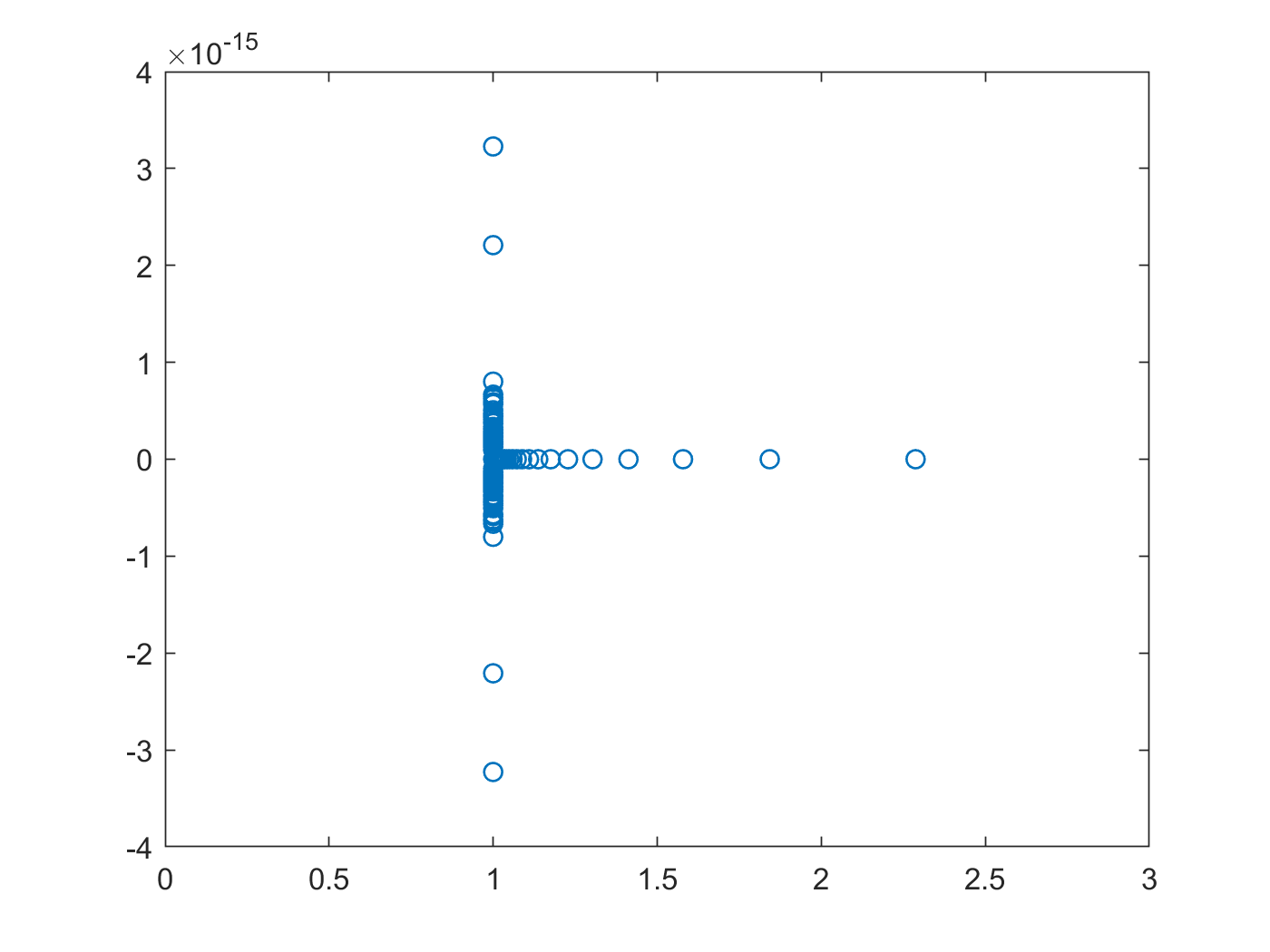}}
    \\
    \subfloat[][$\beta=0.9$, $\omega=1.1$]
    {\includegraphics[width=0.33\textwidth]{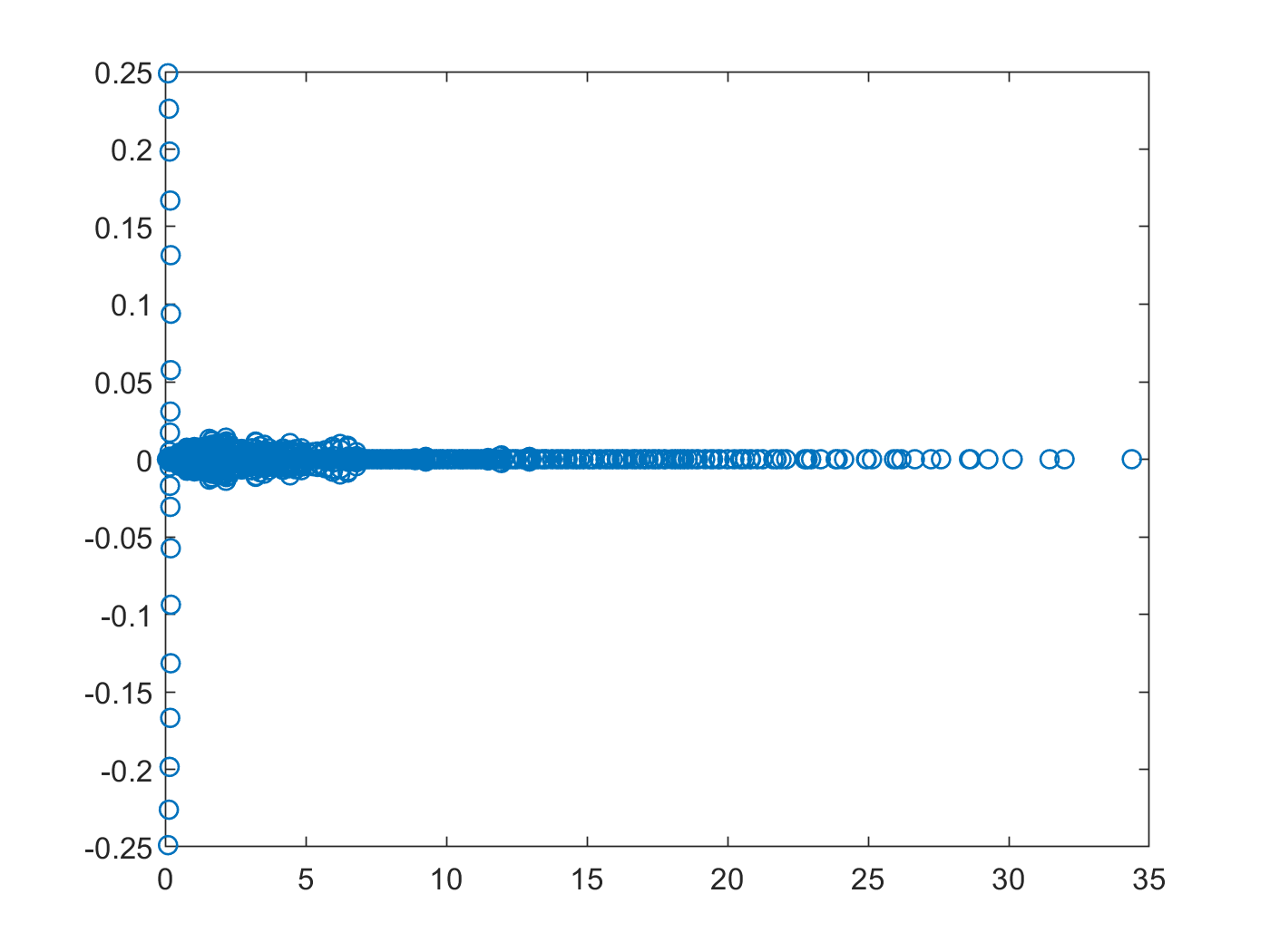}}
    \hfill
    \subfloat[][$\beta=0.9$, $\omega=1.1$]
    {\includegraphics[width=0.33\textwidth]{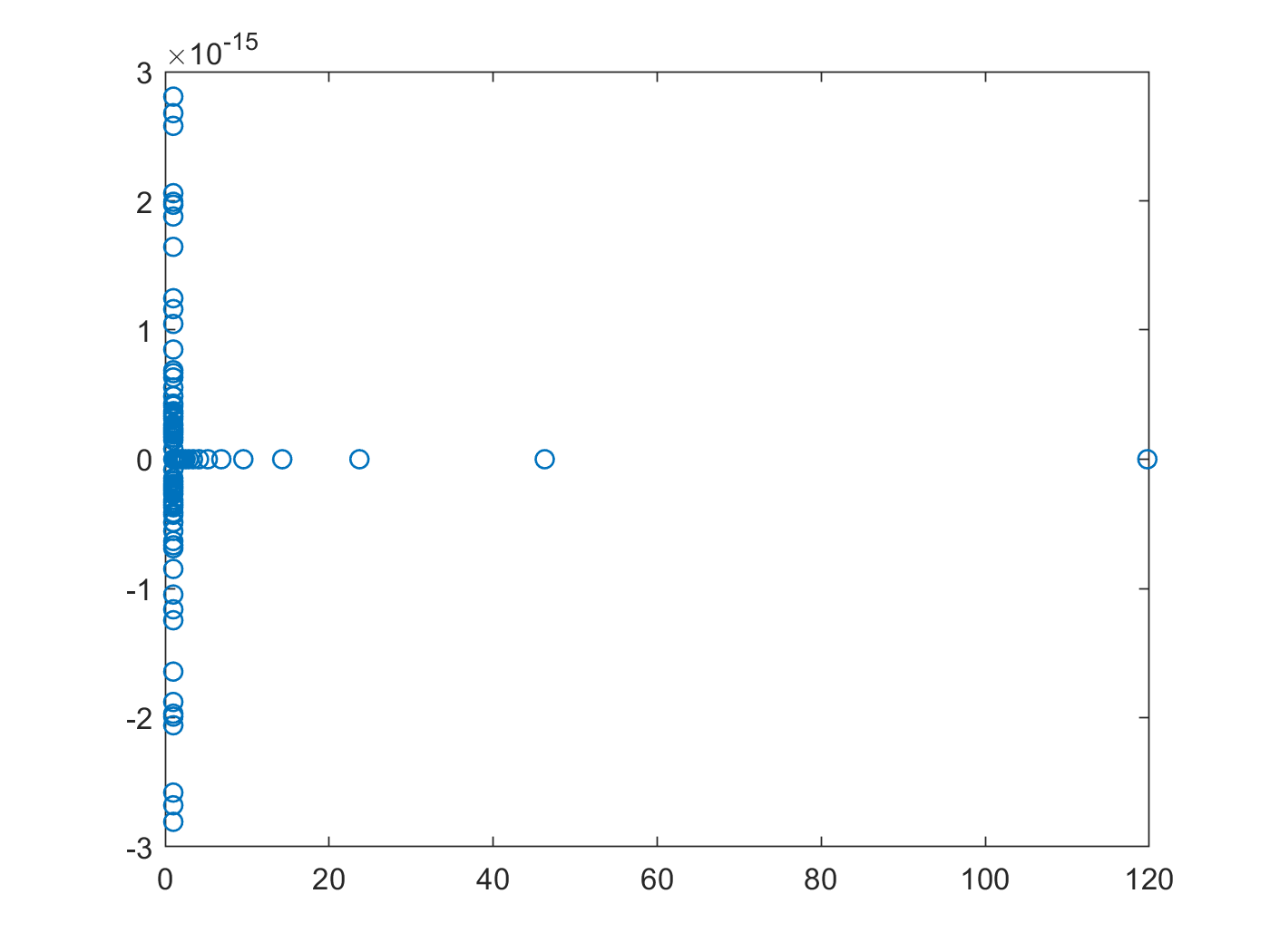}}
    \subfloat[][$\beta=0.9$, $\omega=1.1$]
    {\includegraphics[width=0.33\textwidth]{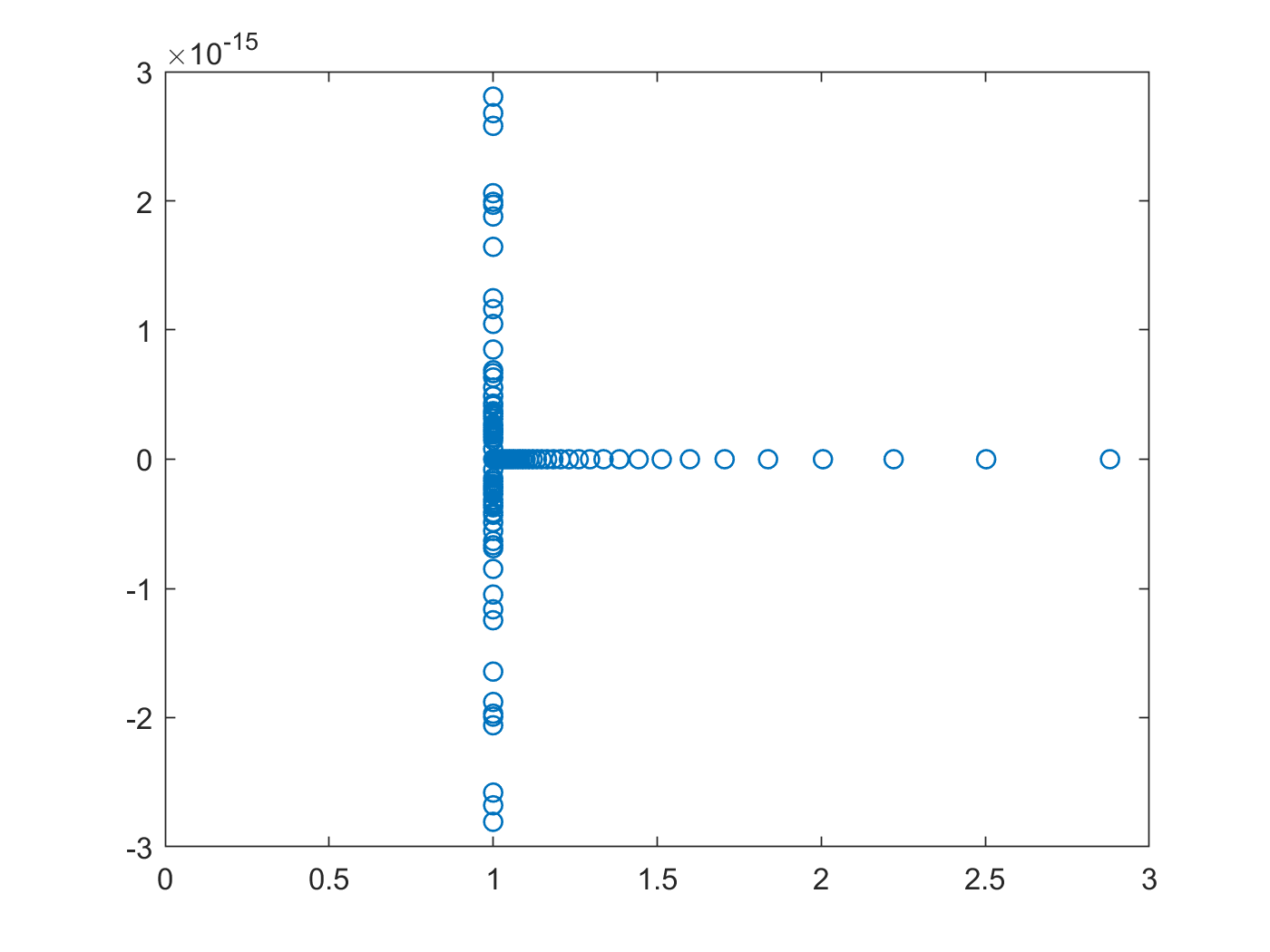}}
    \\
    \caption{1D case - Eigenvalues of $A_{n,S}$ (left column) and $P_{n,S}^{-1}A_{n,S}$ (central column) in the complex plane, for different choices of the fractional orders, with $n=2^6$ and $S=2^4$. For clarity, we also display magnified versions of the central graphs, focused on the cluster at 1 (right column).}
    \label{fig:1D_eigenvalues_n50s10}
\end{figure}

Figure \ref{fig:1D_eigenvalues_n50s10} shows the distribution of the eigenvalues in the complex plane for various combinations of the fractional orders $\beta$ and $\omega$, with fixed problem dimensions $n=2^6$ and $S=2^4$. These relatively small sizes were chosen to ensure that the spectral patterns are clearly visible. For larger dimensions, the spectral behavior remains qualitatively unchanged and is therefore not included. The plots align with the theoretical predictions, showing eigenvalues accumulating near 0 for $A_{n,S}$ and a tight cluster at 1 for the preconditioned matrix, with no eigenvalues close to 0 and few outliers. Since the cluster at 1 is difficult to discern at the original scale, in the third column we included a magnified graph focused on this accumulation point.

Table \ref{table:1D_condnum} presents the 2-norm condition numbers of the nonpreconditioned and preconditioned coefficient matrix, for the same values of $\beta$ and $\omega$ and varying $n$ and $S$. Clearly, $A_{n,S}$ is severely ill-conditioned, even for small problem sizes. However, the proposed $P_{n,S}$ proves to be very effective in reducing the condition number, maintaining a stable performance as $n$ and $S$ grow. It performs especially well when $\beta$ is close to 1 and $\omega$ is close to 2, even though in this case $A_{n,S}$ is particularly ill-conditioned.

\begin{table}[ht]
    \centering
    \caption{1D case - 2-norm condition numbers of the nonpreconditioned and preconditioned coefficient matrix.}
    \label{table:1D_condnum}
    \begin{tabular}{cccccccc}
    \toprule
    & & \multicolumn{2}{c}{$\beta=0.1$, $\omega=1.9$} & \multicolumn{2}{c}{$\beta=0.5$, $\omega=1.5$} & \multicolumn{2}{c}{$\beta=0.9$, $\omega=1.1$} \\
    \cmidrule(lr){3-4}\cmidrule(lr){5-6}\cmidrule(lr){7-8}
    $n$ & $S$ & - & $P_{n,S}$ & - & $P_{n,S}$ & - & $P_{n,S}$ \\
    \midrule
    $2^{4}$ & $2^{4}$ & 2021 & 62 & 731 & 86 & 1371 & 108 \\
    $2^{4}$ & $2^{5}$ & 4125 & 63 & 1500 & 86 & 3655 & 107 \\
    $2^{4}$ & $2^{6}$ & 8176 & 64 & 3084 & 88 & 9879 & 108 \\
    \midrule
    $2^{5}$ & $2^{4}$ & 7656 & 66 & 1511 & 98 & 1688 & 136 \\
    $2^{5}$ & $2^{5}$ & 15644 & 67 & 2787 & 99 & 4430 & 136 \\
    $2^{5}$ & $2^{6}$ & 31110 & 68 & 5206 & 101 & 11878 & 137 \\
    \midrule
    $2^{6}$ & $2^{4}$ & 29513 & 68 & 3763 & 107 & 2098 & 163 \\
    $2^{6}$ & $2^{5}$ & 60215 & 69 & 6383 & 108 & 5199 & 163 \\
    $2^{6}$ & $2^{6}$ & 119676 & 71 & 10877 & 110 & 13497 & 164 \\
    \midrule
    $2^{7}$ & $2^{4}$ & 113517 & 70 & 10326 & 114 & 2857 & 189 \\
    $2^{7}$ & $2^{5}$ & 231376 & 71 & 16832 & 116 & 6392 & 188 \\
    $2^{7}$ & $2^{6}$ & 459245 & 73 & 27133 & 118 & 15500 & 191 \\
    \midrule
    $2^{8}$ & $2^{4}$ & 433429 & 71 & 29317 & 119 & 4461 & 212 \\
    $2^{8}$ & $2^{5}$ & 883159 & 72 & 47086 & 120 & 8733 & 211 \\
    $2^{8}$ & $2^{6}$ & 1751554 & 74 & 74164 & 123 & 19073 & 213 \\
    \bottomrule
    \end{tabular}
\end{table}

\subsubsection{GMRES method performance}

Now, we solve the linear system of the form \eqref{eq:lin-syst} arising from the discretization of the problem, evaluating the number of iterations and CPU time required by the nonpreconditioned and preconditioned GMRES method. The results are gathered in Table \ref{table:1D_gmres}, which reports the number of iterations and CPU times (in seconds) for the nonpreconditioned and preconditioned GMRES, using the same values of $\beta$ and $\omega$ as in the previous subsection. We set the regularization parameter to $\lambda= 5\cdot10^{-3}$ and the noise level to $\varepsilon=0.01$. The initial guess is the zero vector, the tolerance is $\texttt{tol}=10^{-8}$, and the maximum number of iterations is equal to the minimum between the size of the matrix and 1000. When the maximum number of iterations is reached before convergence, in the table we write $>1000$, without reporting the corresponding CPU time.

In the columns corresponding to the nonpreconditioned linear system, we observe high iteration counts, reflecting the severe ill-conditioning of the problem, and consequently we have high CPU times. Both the iteration counts and the CPU times are more sensitive to increases in $n$ than in $S$, in the sense that high values of $n$ result in more iterations and longer CPU times, even when the total matrix size is the same or even smaller due to a smaller $S$. This behavior is reasonable, considering that larger $n$ values lead to a denser coefficient matrix and, as discussed in the previous section, the computational cost is primarily driven by the spatial variable. These trends also align with Table \ref{table:1D_condnum}, which shows that the condition number is generally more influenced by the increase in $n$ than in $S$.

Turning to the columns related to the preconditioner, the iteration counts across the various choices of fractional orders align well with the spectral analysis. $P_{n,S}$ demonstrates an excellent performance, even showing nearly optimal convergence, in the sense that the number of iterations is almost independent from the matrix size and remains practically constant as $n$ and $S$ increase. The corresponding CPU times remain very low, reflecting this good behavior. In accordance with the spectral analysis presented in the previous subsection, the number of iterations required by the preconditioned system is especially low when $\beta$ is close to 1 and $\omega$ is close to 2, while in the nonpreconditioned case we have particularly high iteration counts.

\begin{table}[ht]
    \centering
    \caption{1D case - Iterations and CPU times required for solving the nonpreconditioned and preconditioned linear system with the GMRES method.}
    \label{table:1D_gmres}
    \begin{tabular}{cccccccccccccc}
    \toprule
    & & \multicolumn{4}{c}{$\beta=0.1$, $\omega=1.9$} & \multicolumn{4}{c}{$\beta=0.5$, $\omega=1.5$} & \multicolumn{4}{c}{$\beta=0.9$, $\omega=1.1$} \\
    \cmidrule(lr){3-6}\cmidrule(lr){7-10}\cmidrule(lr){11-14}
    & & \multicolumn{2}{c}{Iter} & \multicolumn{2}{c}{CPU} & \multicolumn{2}{c}{Iter} & \multicolumn{2}{c}{CPU} & \multicolumn{2}{c}{Iter} & \multicolumn{2}{c}{CPU} \\
    \cmidrule(lr){3-4}\cmidrule(lr){5-6}\cmidrule(lr){7-8}\cmidrule(lr){9-10}\cmidrule(lr){11-12}\cmidrule(lr){13-14}
    $n$ & $S$ & - & $P_{n,S}$ & - & $P_{n,S}$ & - & $P_{n,S}$ & - & $P_{n,S}$ & - & $P_{n,S}$ & - & $P_{n,S}$ \\
    \midrule
    $2^{4}$ & $2^{4}$ & 175 & 9 & 0.19 & 0.04 & 117 & 10 & 0.06 & 0.01 & 112 & 13 & 0.07 & 0.02 \\
    $2^{4}$ & $2^{5}$ & 258 & 9 & 0.22 & 0.02 & 151 & 10 & 0.11 & 0.02 & 157 & 13 & 0.12 & 0.02 \\
    $2^{4}$ & $2^{6}$ & 372 & 8 & 0.69 & 0.03 & 185 & 11 & 0.25 & 0.04 & 277 & 13 & 0.43 & 0.05 \\
    $2^{4}$ & $2^{7}$ & 505 & 8 & 1.82 & 0.06 & 283 & 11 & 0.78 & 0.08 & 536 & 13 & 1.69 & 0.09 \\
    $2^{4}$ & $2^{8}$ & 657 & 9 & 11.61 & 0.17 & 413 & 11 & 4.71 & 0.20 & $>$1000 & 13 & - & 0.24 \\
    \midrule
    $2^{5}$ & $2^{4}$ & 318 & 8 & 0.21 & 0.01 & 172 & 11 & 0.08 & 0.02 & 151 & 15 & 0.07 & 0.02 \\
    $2^{5}$ & $2^{5}$ & 465 & 9 & 0.67 & 0.03 & 223 & 11 & 0.22 & 0.04 & 192 & 15 & 0.18 & 0.05 \\
    $2^{5}$ & $2^{6}$ & 690 & 9 & 2.38 & 0.06 & 264 & 11 & 0.58 & 0.07 & 339 & 15 & 0.87 & 0.10 \\
    $2^{5}$ & $2^{7}$ & 923 & 9 & 22.93 & 0.15 & 334 & 11 & 3.23 & 0.18 & 663 & 15 & 11.39 & 0.23 \\
    $2^{5}$ & $2^{8}$ & $>$1000 & 9 & - & 0.37 & 463 & 12 & 11.95 & 0.46 & $>$1000 & 16 & - & 0.60 \\
    \midrule
    $2^{6}$ & $2^{4}$ & 600 & 8 & 0.95 & 0.04 & 288 & 11 & 0.27 & 0.04 & 192 & 17 & 0.14 & 0.06 \\
    $2^{6}$ & $2^{5}$ & 911 & 9 & 3.26 & 0.06 & 349 & 11 & 0.66 & 0.07 & 246 & 17 & 0.42 & 0.11 \\
    $2^{6}$ & $2^{6}$ & $>$1000 & 9 & - & 0.14 & 425 & 11 & 4.10 & 0.16 & 393 & 17 & 3.44 & 0.23 \\
    $2^{6}$ & $2^{7}$ & $>$1000 & 9 & - & 0.27 & 489 & 12 & 10.58 & 0.35 & 779 & 17 & 24.19 & 0.50 \\
    $2^{6}$ & $2^{8}$ & $>$1000 & 9 & - & 0.69 & 575 & 12 & 24.16 & 0.82 & $>$1000 & 17 & - & 1.16 \\
    \midrule
    $2^{7}$ & $2^{4}$ & $>$1000 & 9 & - & 0.12 & 479 & 12 & 1.07 & 0.16 & 250 & 18 & 0.38 & 0.22 \\
    $2^{7}$ & $2^{5}$ & $>$1000 & 9 & - & 0.33 & 590 & 11 & 7.88 & 0.37 & 306 & 18 & 2.32 & 0.47 \\
    $2^{7}$ & $2^{6}$ & $>$1000 & 9 & - & 0.49 & 710 & 12 & 19.49 & 0.61 & 454 & 18 & 8.57 & 0.86 \\
    $2^{7}$ & $2^{7}$ & $>$1000 & 9 & - & 1.03 & 808 & 12 & 37.41 & 1.34 & 864 & 18 & 42.82 & 1.98 \\
    $2^{7}$ & $2^{8}$ & $>$1000 & 9 & - & 2.28 & 840 & 12 & 72.49 & 2.88 & $>$1000 & 18 & - & 4.29 \\
    \midrule
    $2^{8}$ & $2^{4}$ & $>$1000 & 9 & - & 0.40 & 839 & 12 & 17.34 & 0.57 & 351 & 19 & 3.18 & 0.75 \\
    $2^{8}$ & $2^{5}$ & $>$1000 & 9 & - & 0.78 & 985 & 12 & 37.30 & 1.00 & 400 & 19 & 6.83 & 1.34 \\
    $2^{8}$ & $2^{6}$ & $>$1000 & 9 & - & 1.40 & $>$1000 & 12 & - & 1.99 & 521 & 19 & 15.84 & 2.87 \\
    $2^{8}$ & $2^{7}$ & $>$1000 & 9 & - & 2.92 & $>$1000 & 12 & - & 3.69 & 980 & 19 & 82.34 & 5.33 \\
    $2^{8}$ & $2^{8}$ & $>$1000 & 9 & - & 5.99 & $>$1000 & 13 & - & 8.10 & $>$1000 & 19 & - & 11.10 \\
    \bottomrule
    \end{tabular}
\end{table}

\subsubsection{Accuracy of the reconstruction}

Finally, we assess the accuracy of the method
in reconstructing $f(x)$.

\begin{figure}[ht]
    \centering
    \subfloat[][$\lambda=10^{-5}$, $\varepsilon=0.001$]
    {\includegraphics[width=0.33\textwidth]{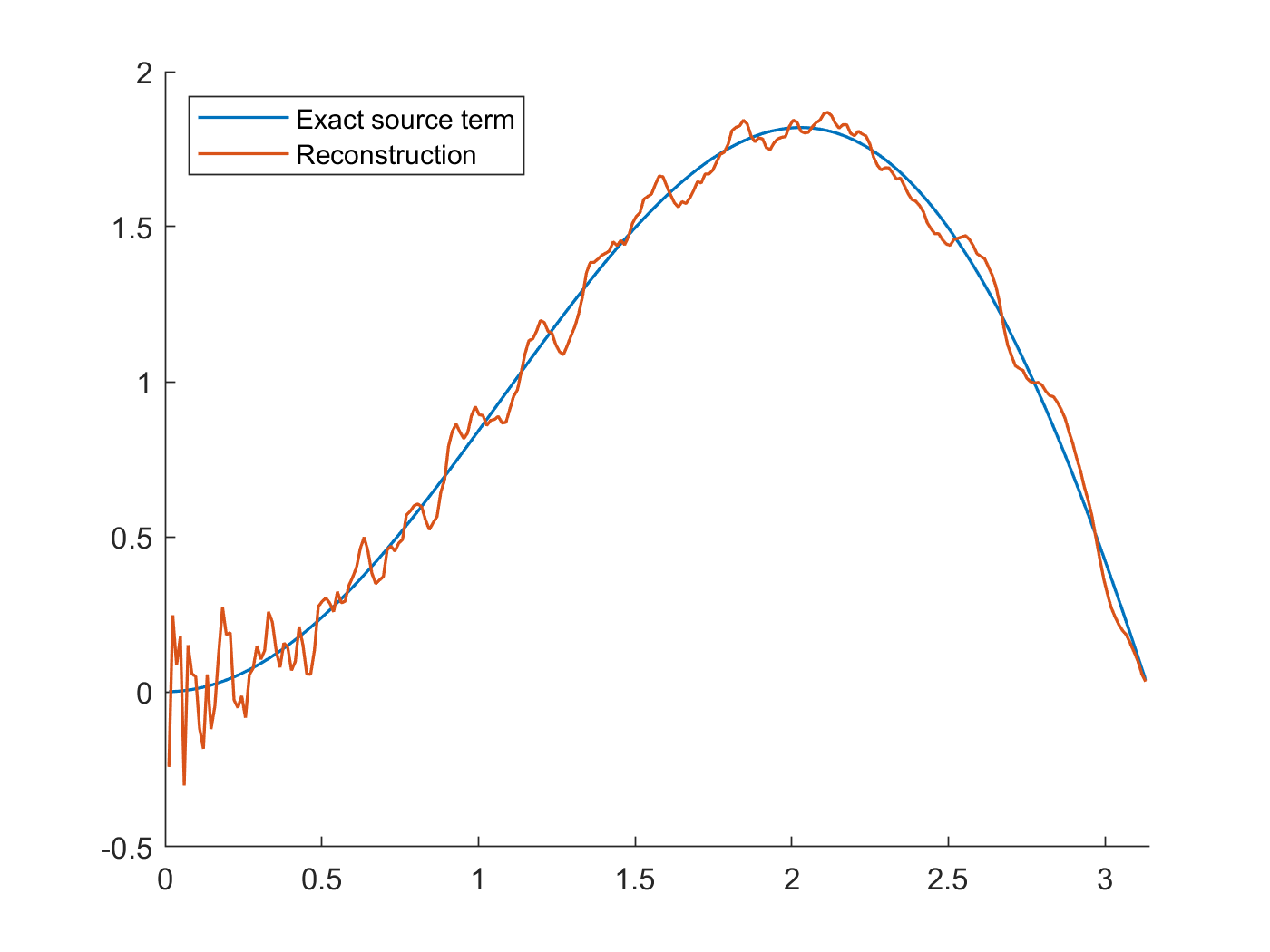}}
    \hfill
    \subfloat[][$\lambda=10^{-4}$, $\varepsilon=0.001$]
    {\includegraphics[width=0.33\textwidth]{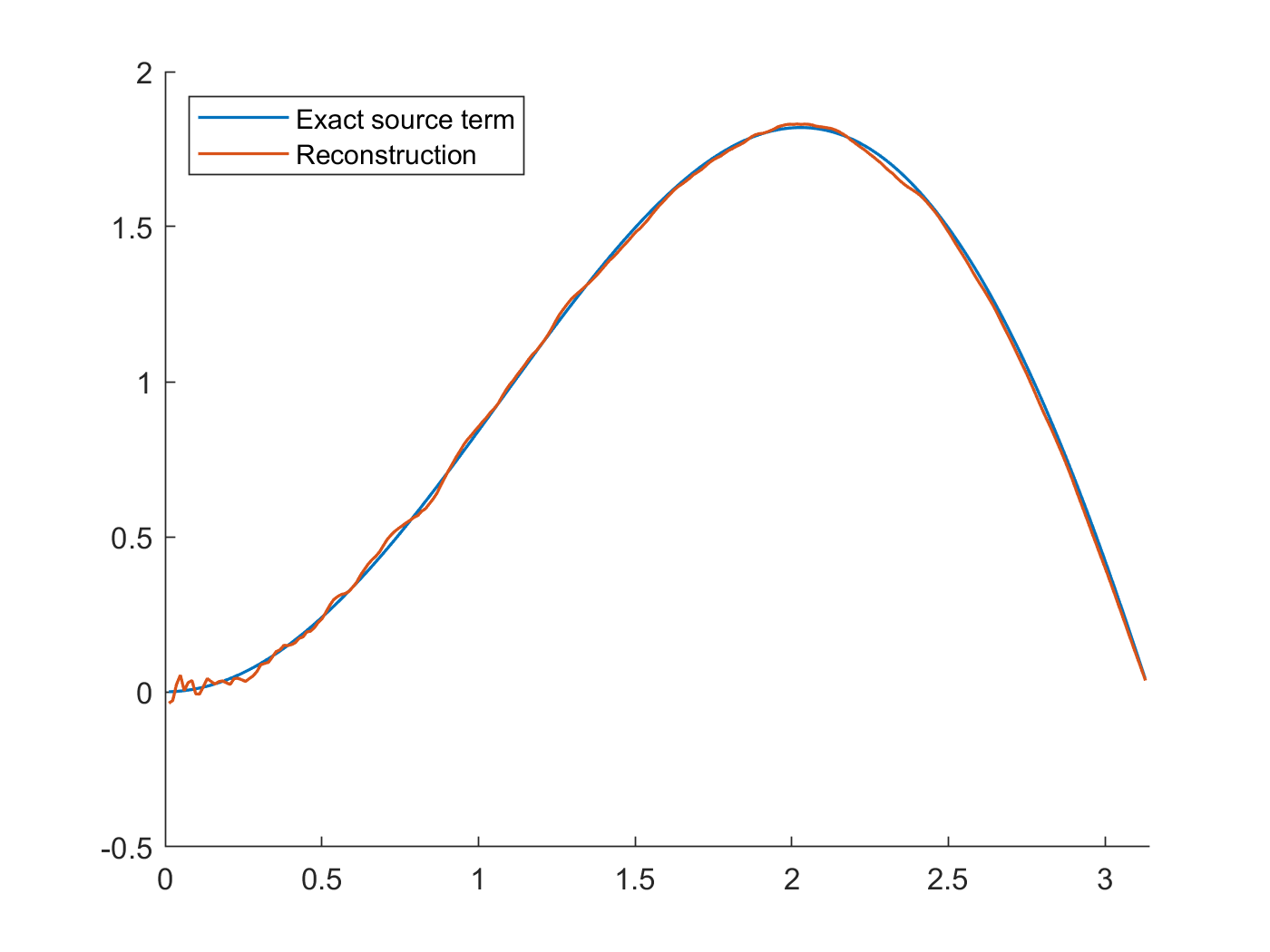}}
    \hfill
    \subfloat[][$\lambda=10^{-3}$, $\varepsilon=0.001$]
    {\includegraphics[width=0.33\textwidth]{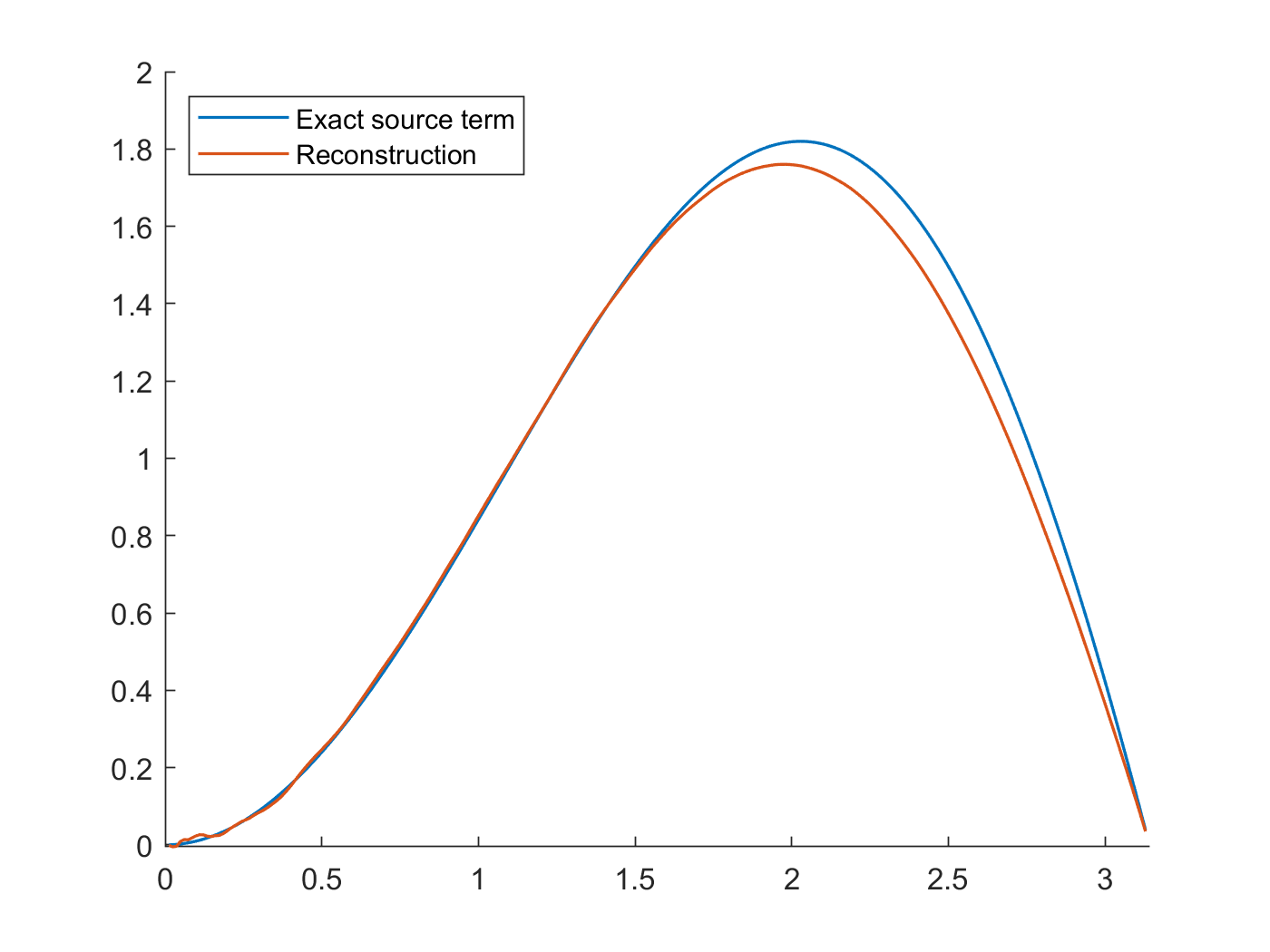}}
    \\
    \subfloat[][$\lambda=10^{-4}$, $\varepsilon=0.01$]
    {\includegraphics[width=0.33\textwidth]{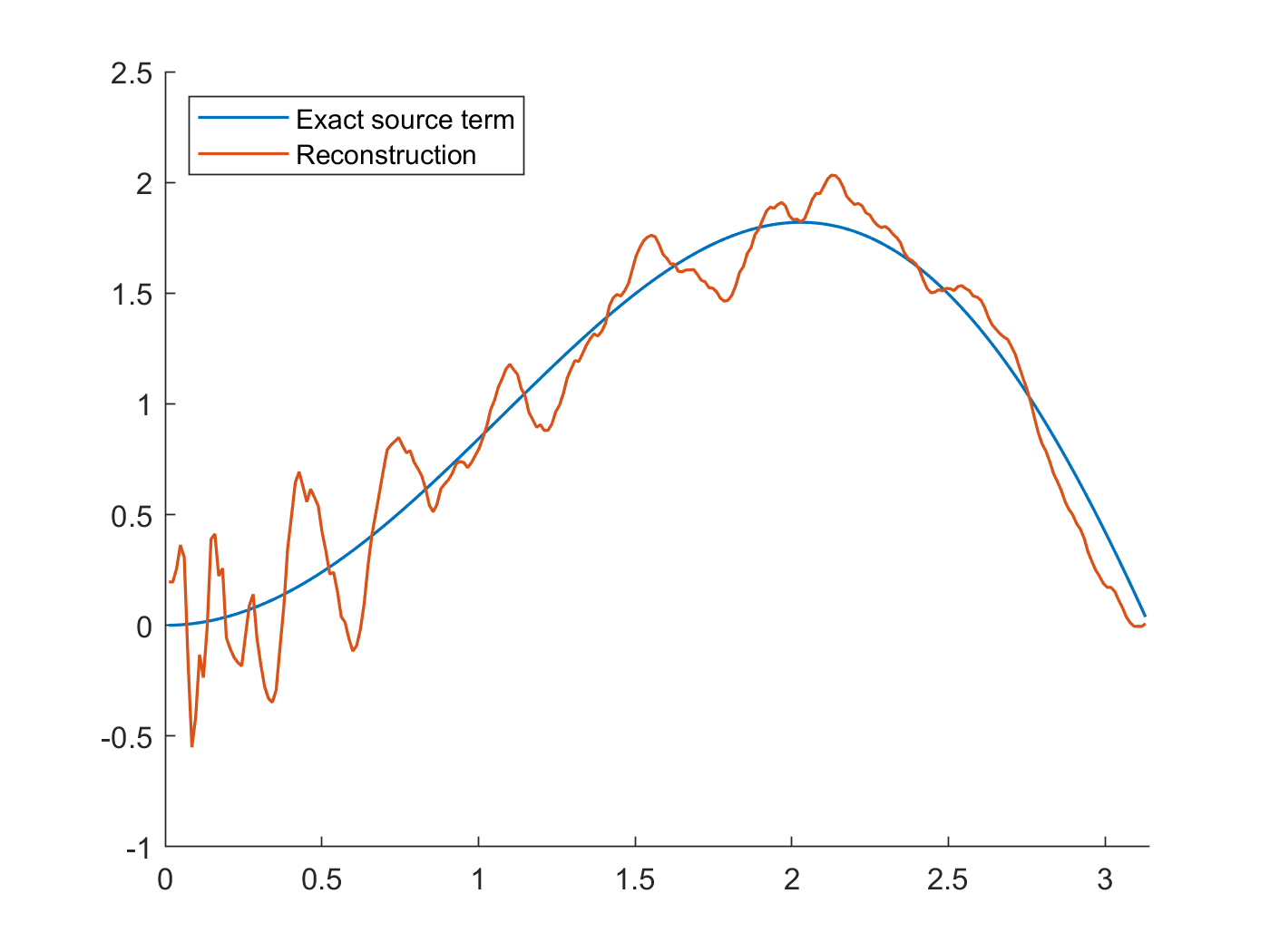}}
    \hfill
    \subfloat[][$\lambda=10^{-3}$, $\varepsilon=0.01$]
    {\includegraphics[width=0.33\textwidth]{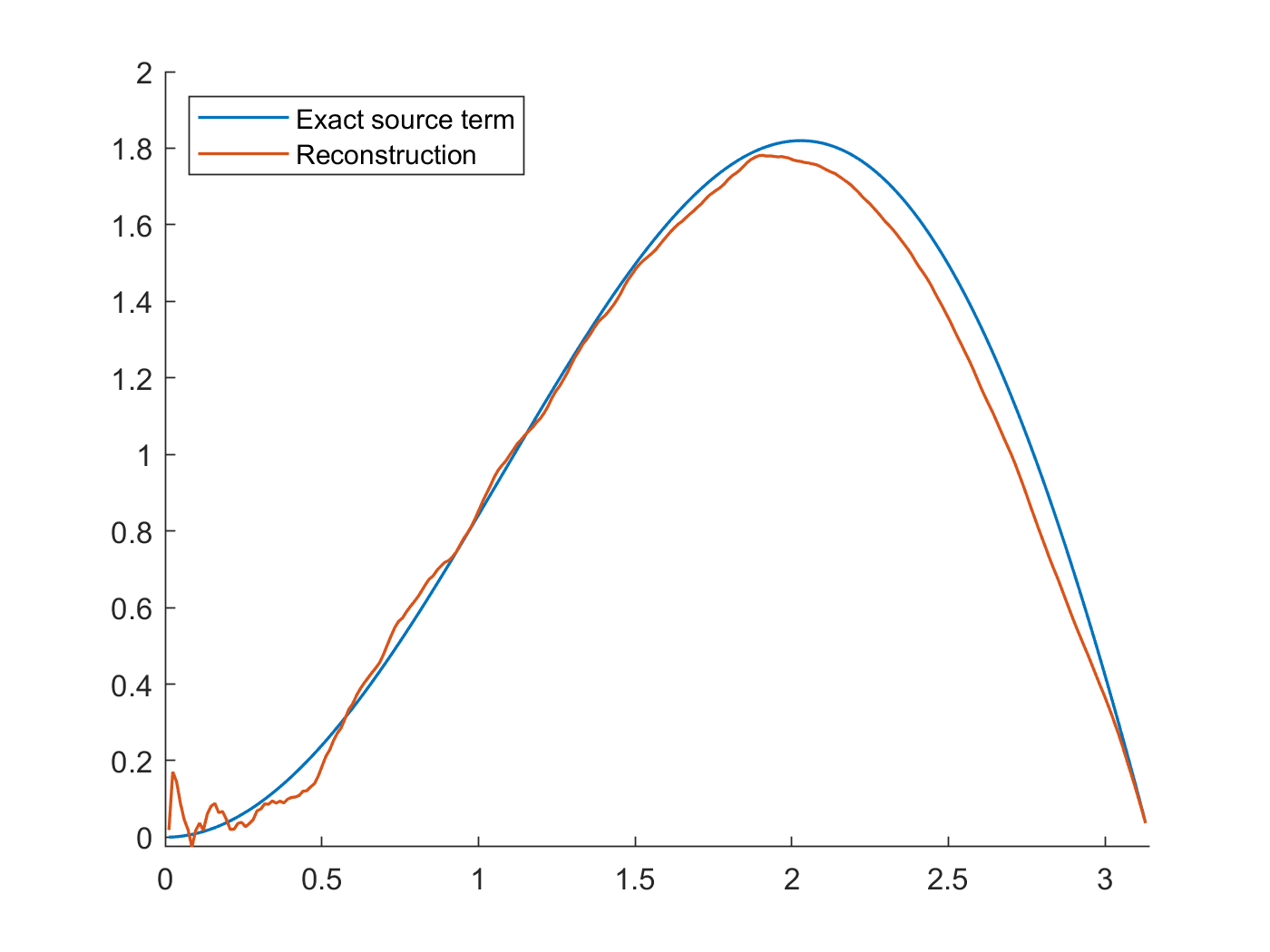}}
    \hfill
    \subfloat[][$\lambda=10^{-2}$, $\varepsilon=0.01$]
    {\includegraphics[width=0.33\textwidth]{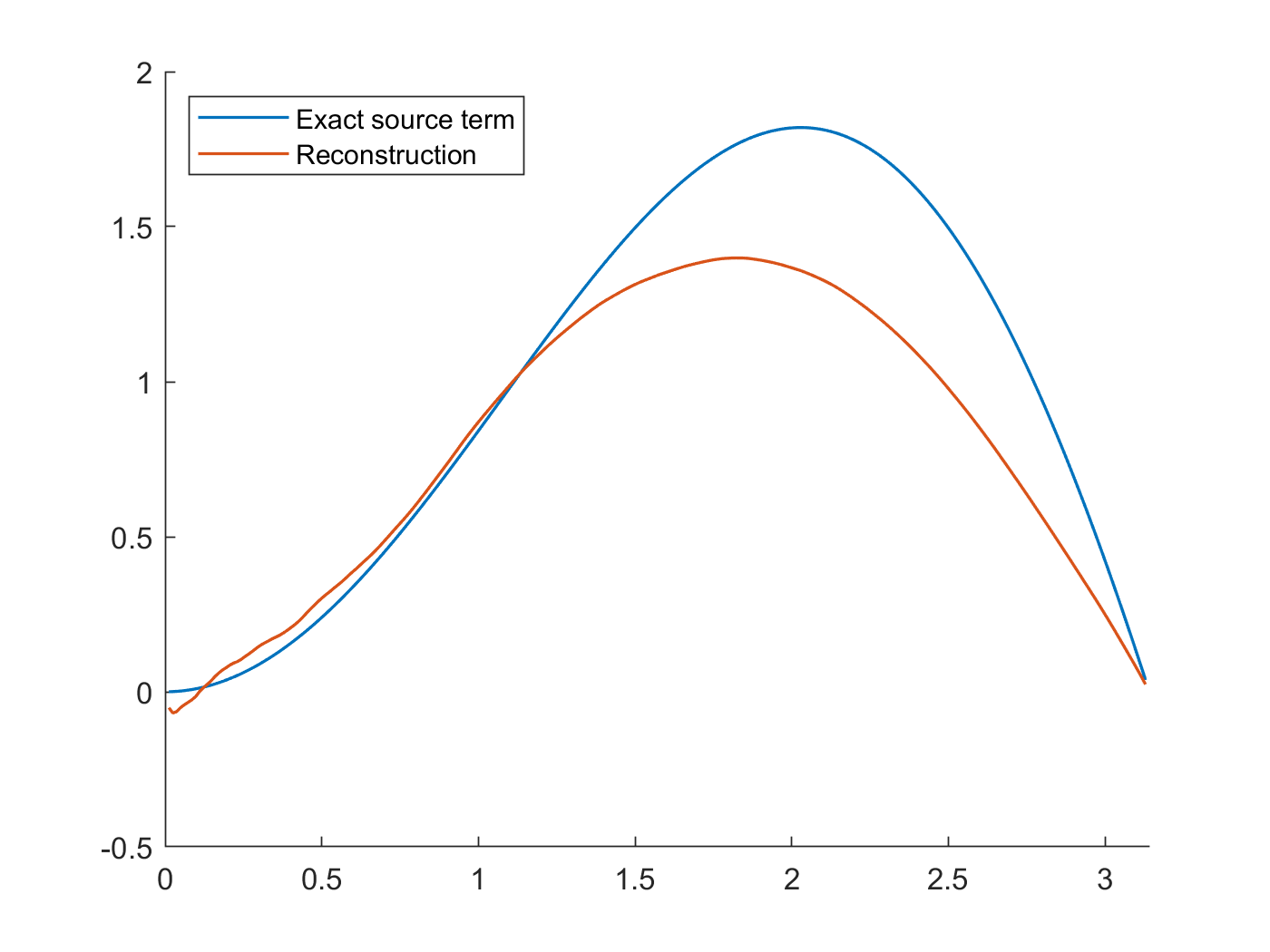}}
    \caption{1D case - Comparison between the exact source function $f(x)$, in blue, and its reconstruction, in red, for different choices of the regularization parameter $\lambda$ and noise level $\varepsilon$. The problem sizes are $n=2^8$, $S=2^6$ and the fractional orders are $\beta=0.1$, $\omega=1.9$.}
    \label{fig:1D_accuracy_beta0.1-omega1.9}
\end{figure}

Figure \ref{fig:1D_accuracy_beta0.1-omega1.9} shows the function $f(x)$ and its reconstruction under varying noise levels $\varepsilon$ and choices of regularization parameter $\lambda$, with the fractional orders fixed at $\beta=0.1$, $\omega=1.9$ and sizes $n=2^8$, $S=2^6$. In the first row, corresponding to a low noise level $\varepsilon=0.001$, the reconstruction achieves high accuracy when $\lambda=10^{-4}$. A smaller value, $\lambda=10^{-5}$, results in overfitting to noise, introducing oscillations. On the other hand, with a larger value $\lambda=10^{-3}$ we begin to see an oversmoothing effect. As the noise level increases to $\varepsilon=0.01$ (second row),  the reconstruction becomes more challenging. Nonetheless, $\lambda=10^{-3}$ provides a reasonable balance between fidelity and regularity. With the other two values we obtain similar deficiencies as before: $\lambda=10^{-4}$ introduces oscillations, while $\lambda=10^{-2}$ causes oversmoothing.

\begin{figure}[ht]
    \centering
    \subfloat[][$\lambda=10^{-5}$, $\varepsilon=0.001$]
    {\includegraphics[width=0.33\textwidth]{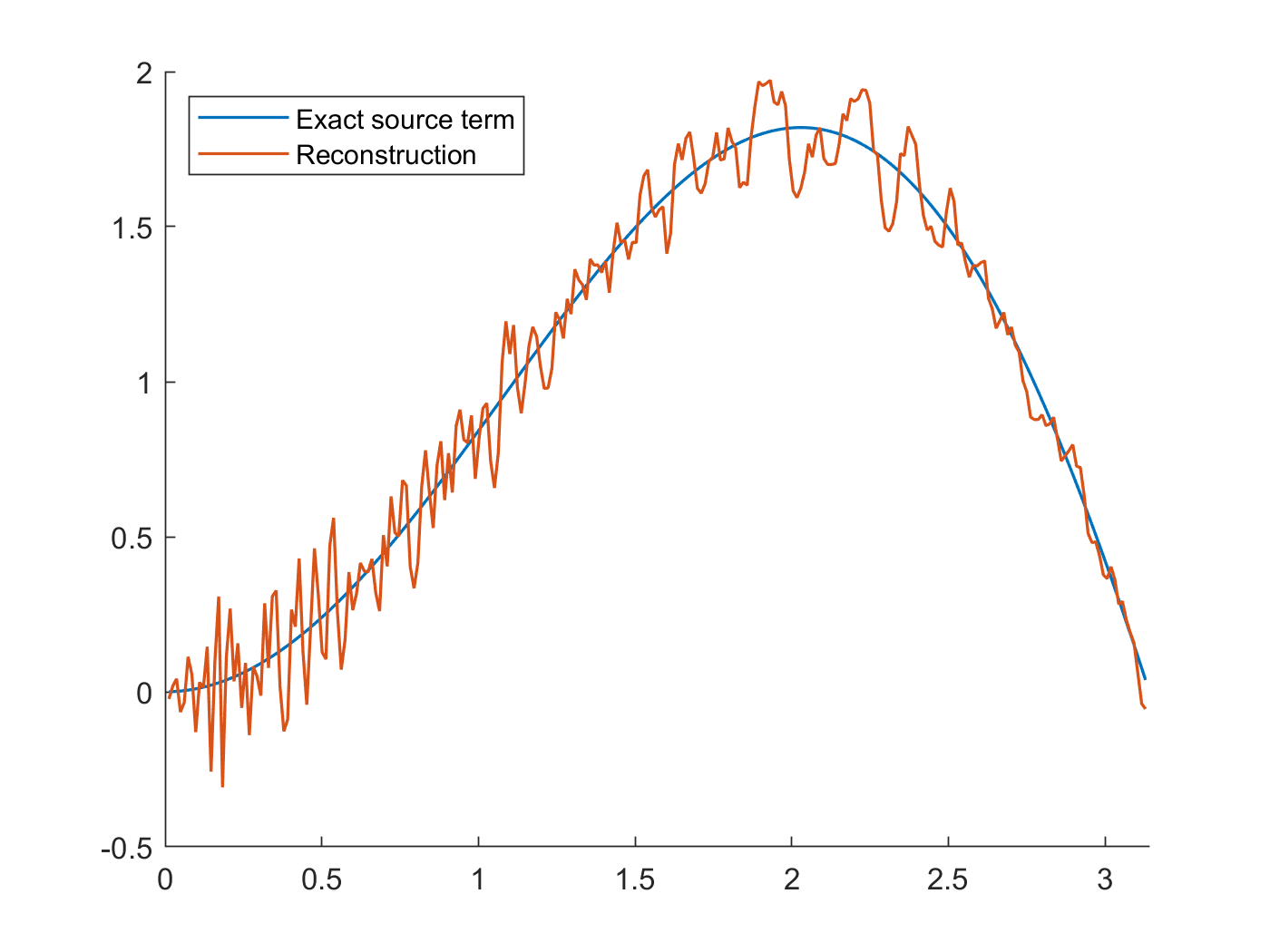}}
    \hfill
    \subfloat[][$\lambda=10^{-4}$, $\varepsilon=0.001$]
    {\includegraphics[width=0.33\textwidth]{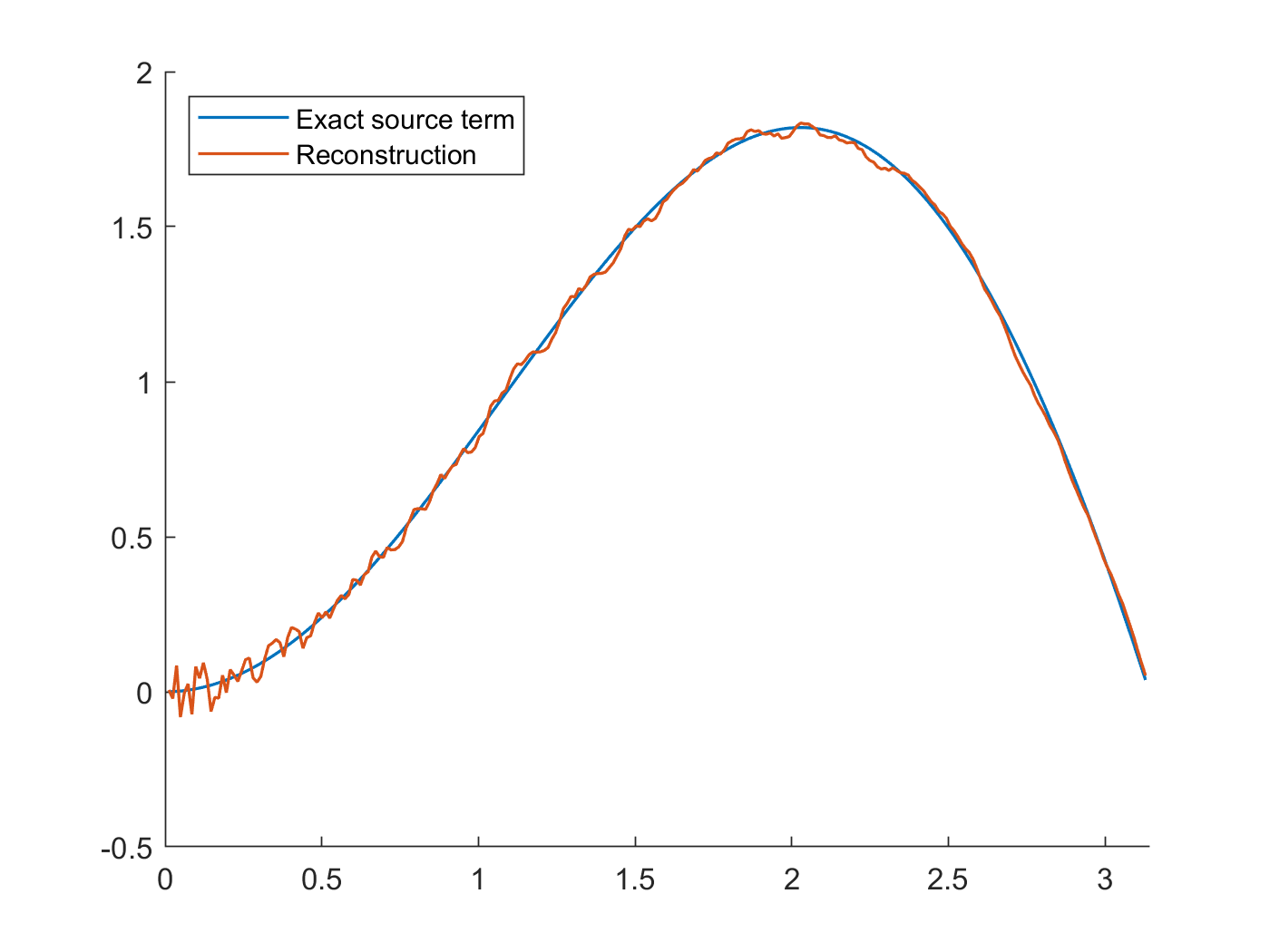}}
    \hfill
    \subfloat[][$\lambda=10^{-3}$, $\varepsilon=0.001$]
    {\includegraphics[width=0.33\textwidth]{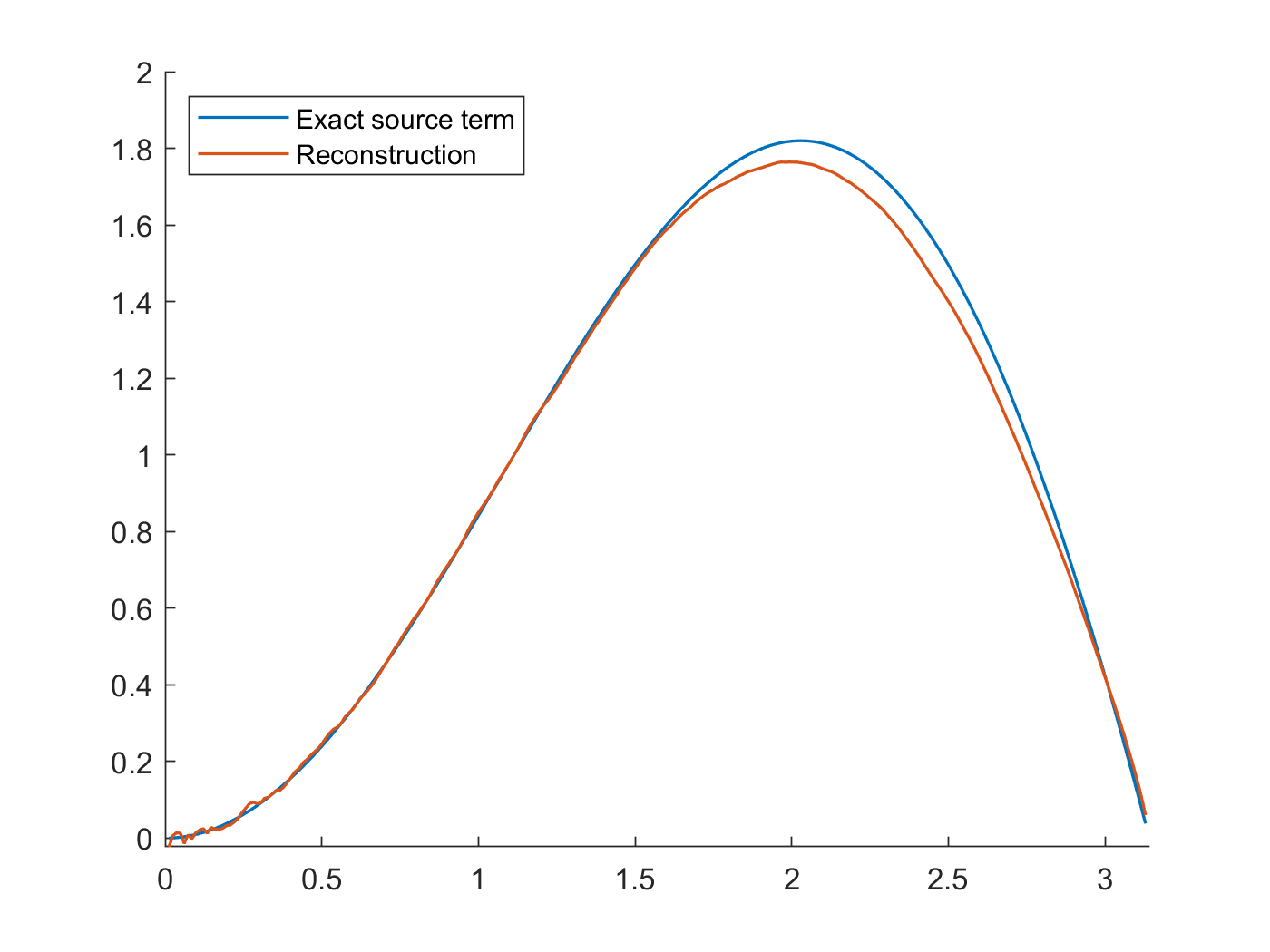}}
    \\
    \subfloat[][$\lambda=10^{-4}$, $\varepsilon=0.01$]
    {\includegraphics[width=0.33\textwidth]{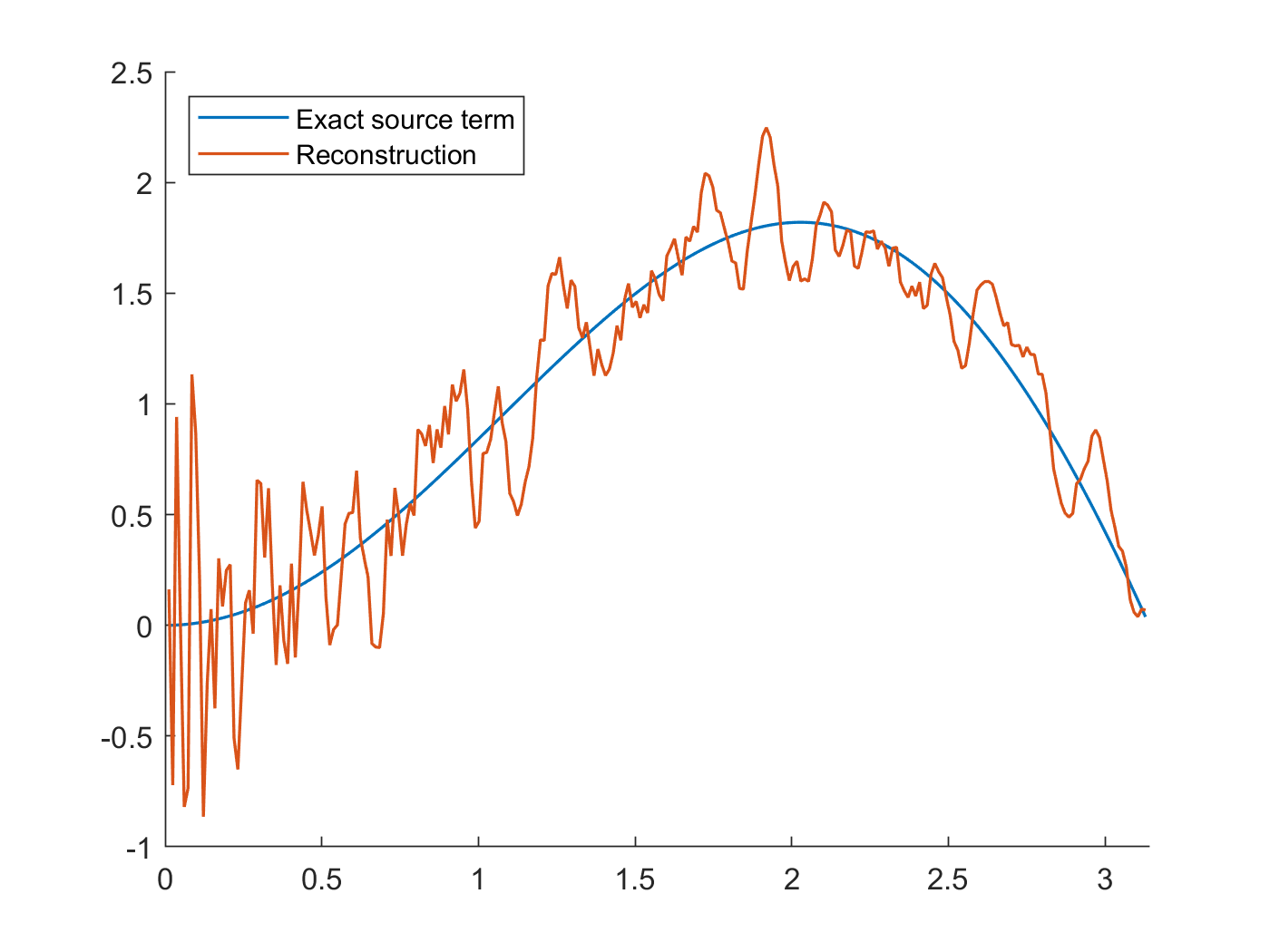}}
    \hfill
    \subfloat[][$\lambda=10^{-3}$, $\varepsilon=0.01$]
    {\includegraphics[width=0.33\textwidth]{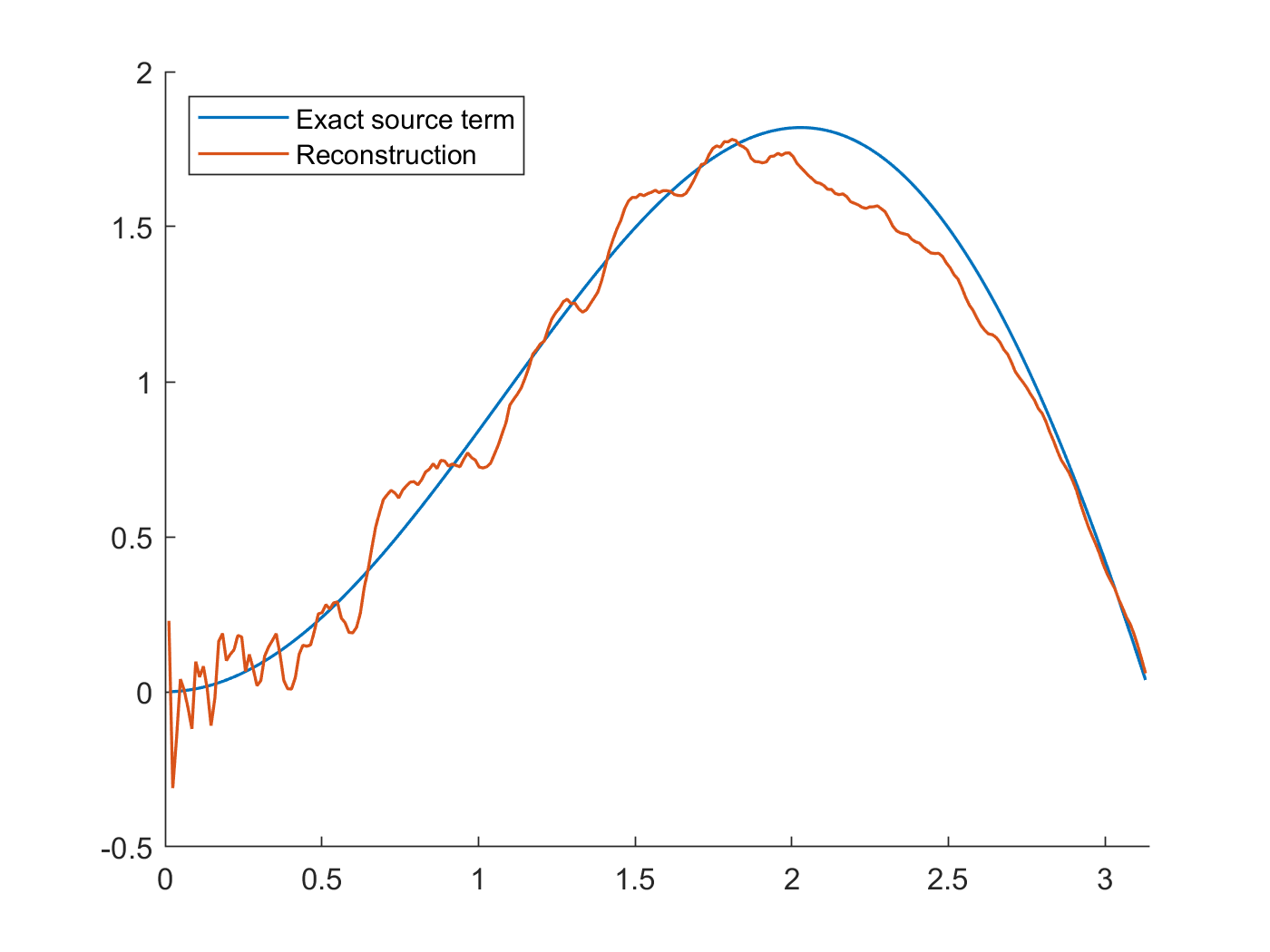}}
    \hfill
    \subfloat[][$\lambda=10^{-2}$, $\varepsilon=0.01$]
    {\includegraphics[width=0.33\textwidth]{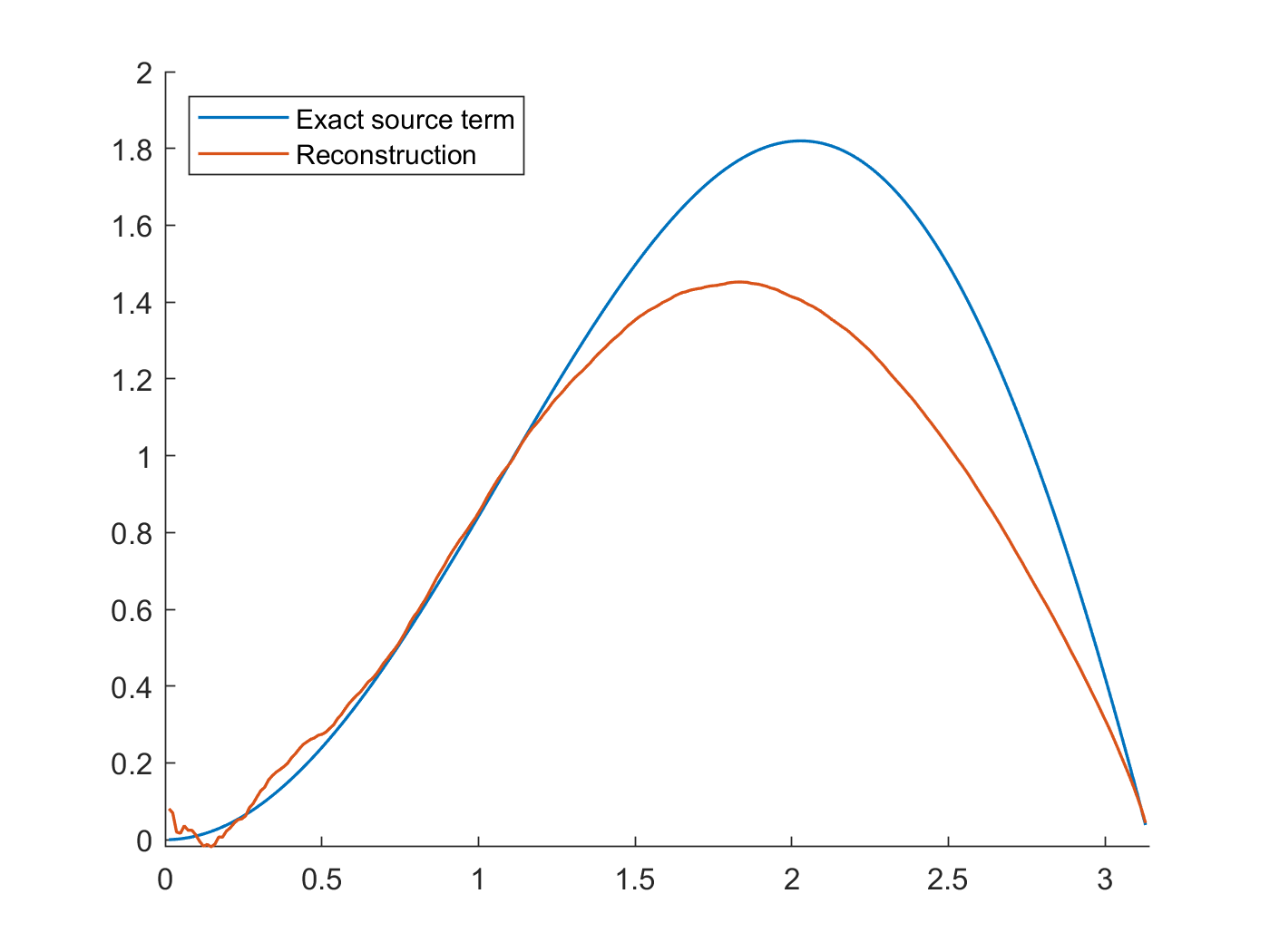}}
    \\
    \caption{1D case - Comparison between the exact source function $f(x)$, in blue, and its reconstruction, in red, for different choices of the regularization parameter $\lambda$ and noise level $\varepsilon$. The problem sizes are $n=2^8$, $S=2^6$ and the fractional orders are $\beta=0.5$, $\omega=1.5$.}
    \label{fig:1D_accuracy_beta0.5-omega1.5}
\end{figure}

\begin{figure}[ht]
    \centering
    \subfloat[][$\lambda=10^{-5}$, $\varepsilon=0.001$]
    {\includegraphics[width=0.33\textwidth]{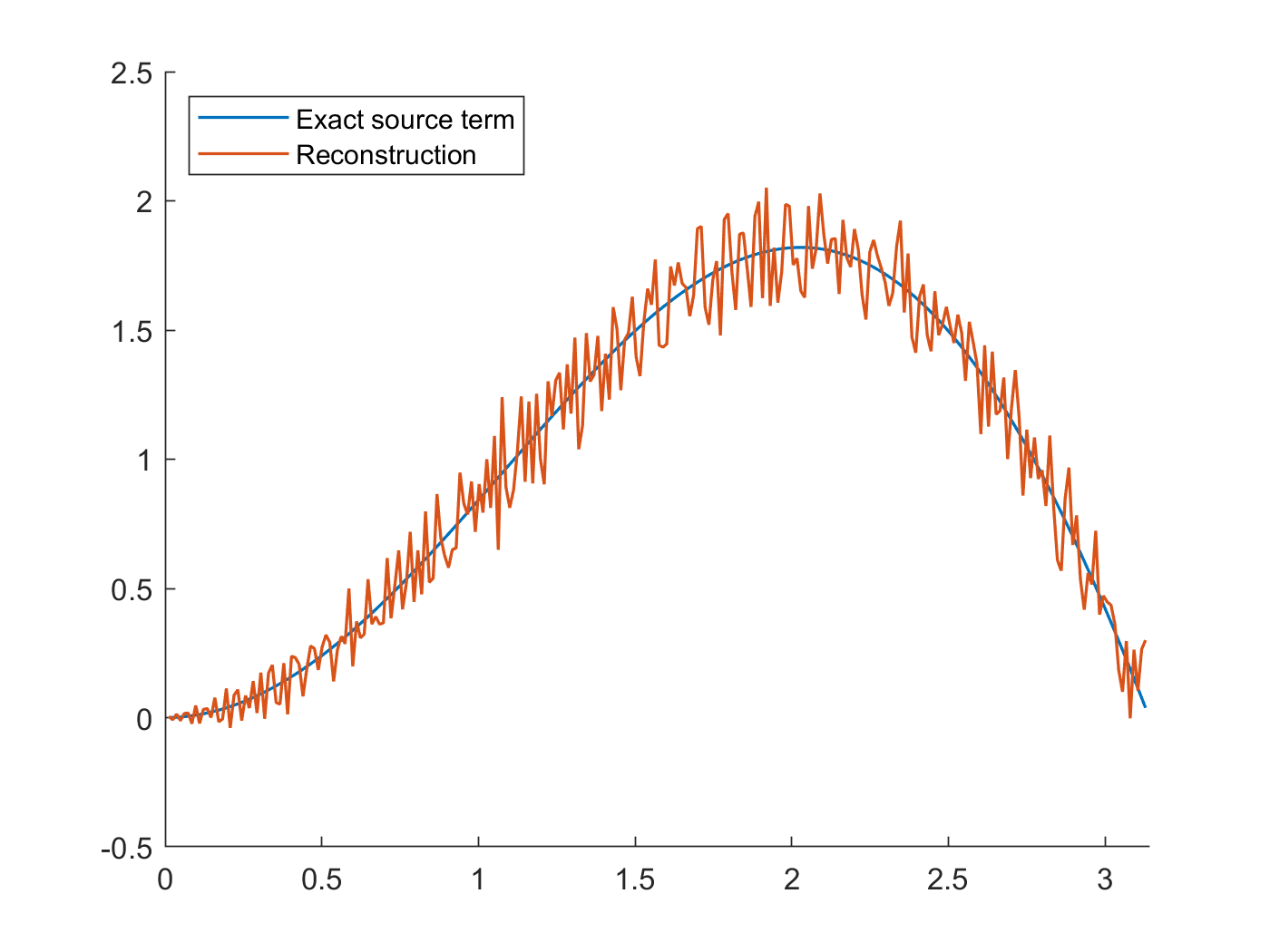}}
    \hfill
    \subfloat[][$\lambda=10^{-4}$, $\varepsilon=0.001$]
    {\includegraphics[width=0.33\textwidth]{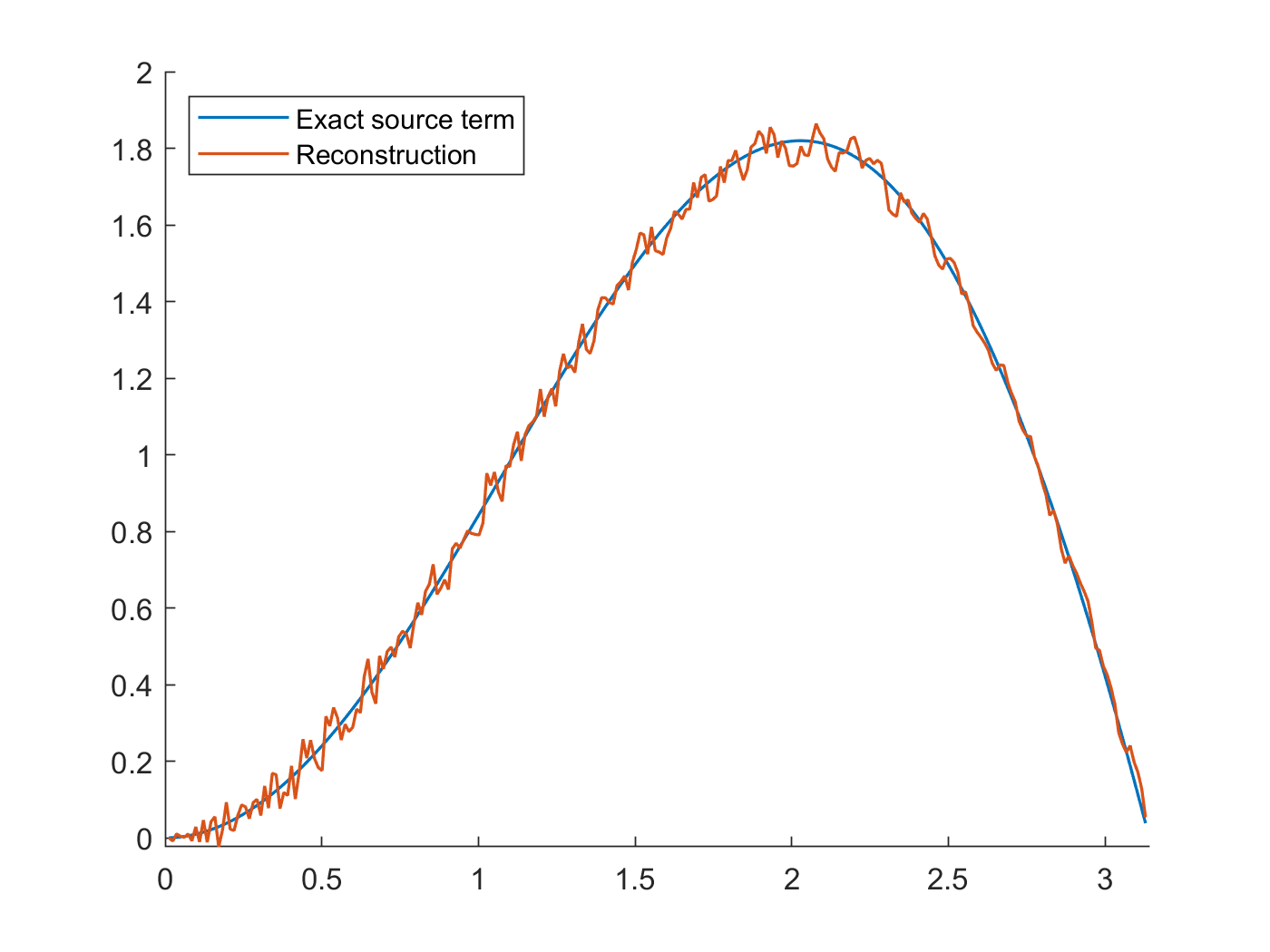}}
    \hfill
    \subfloat[][$\lambda=10^{-3}$, $\varepsilon=0.001$]
    {\includegraphics[width=0.33\textwidth]{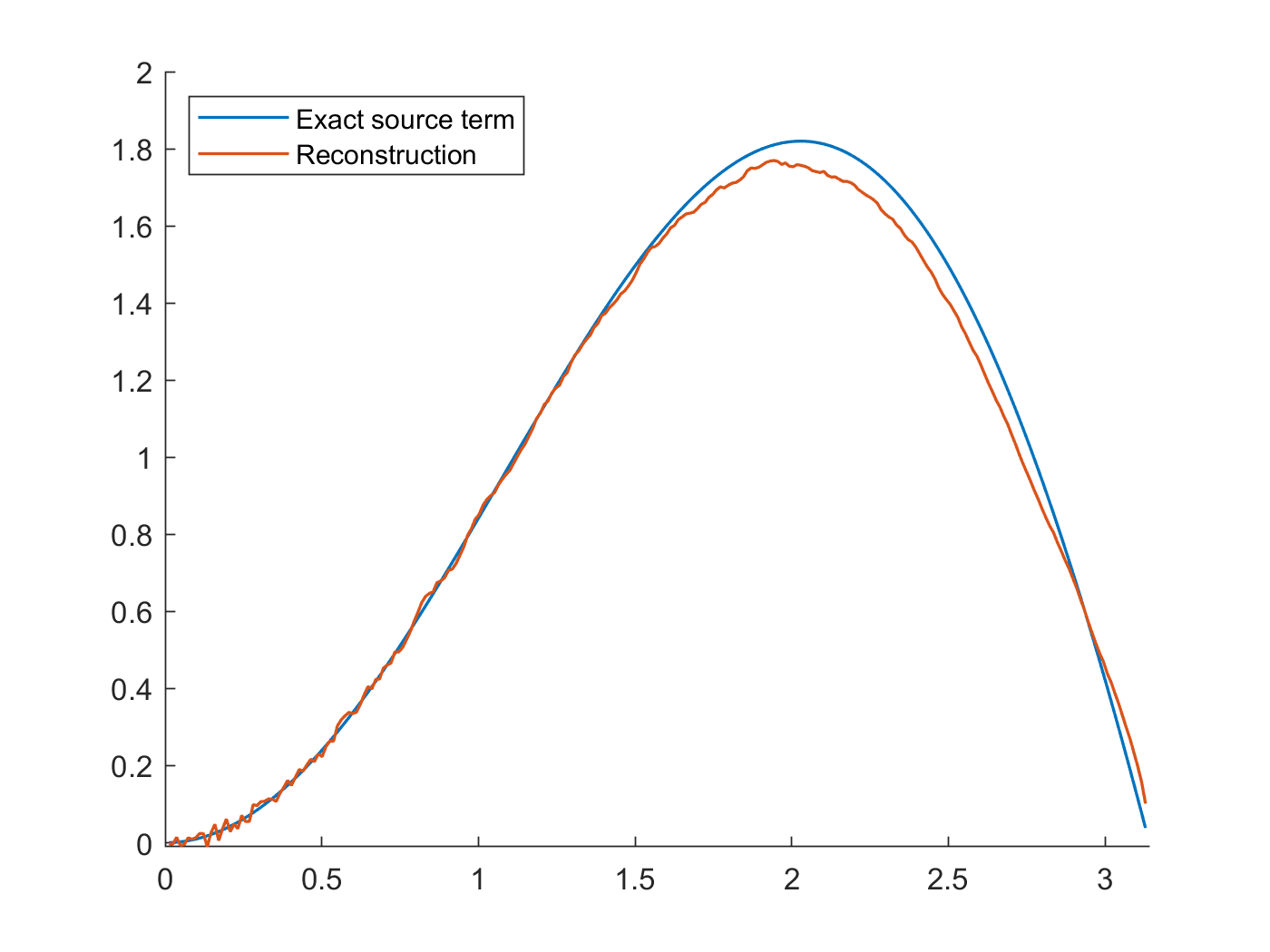}}
    \\
    \subfloat[][$\lambda=10^{-4}$, $\varepsilon=0.01$]
    {\includegraphics[width=0.33\textwidth]{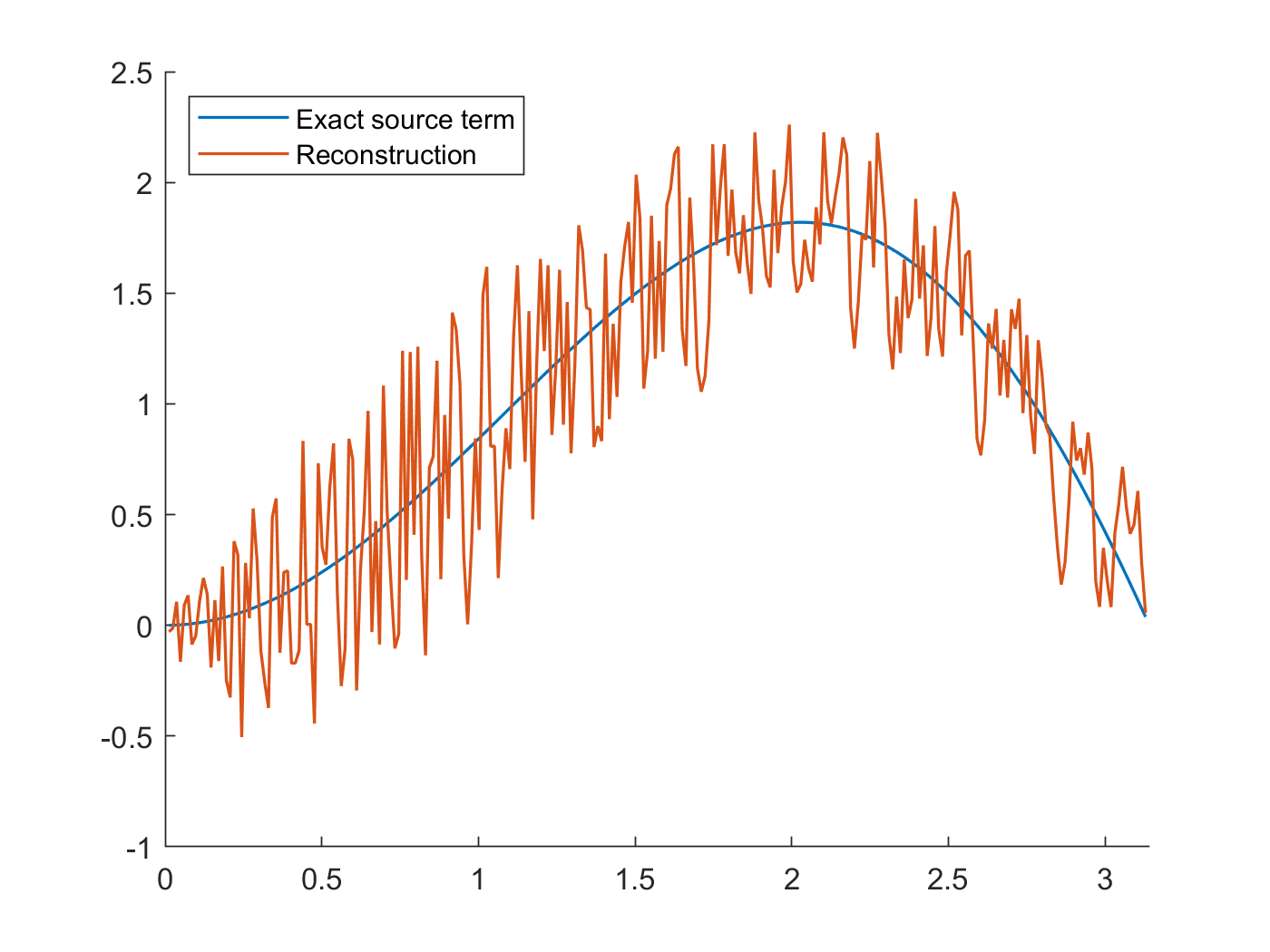}}
    \hfill
    \subfloat[][$\lambda=10^{-3}$, $\varepsilon=0.01$]
    {\includegraphics[width=0.33\textwidth]{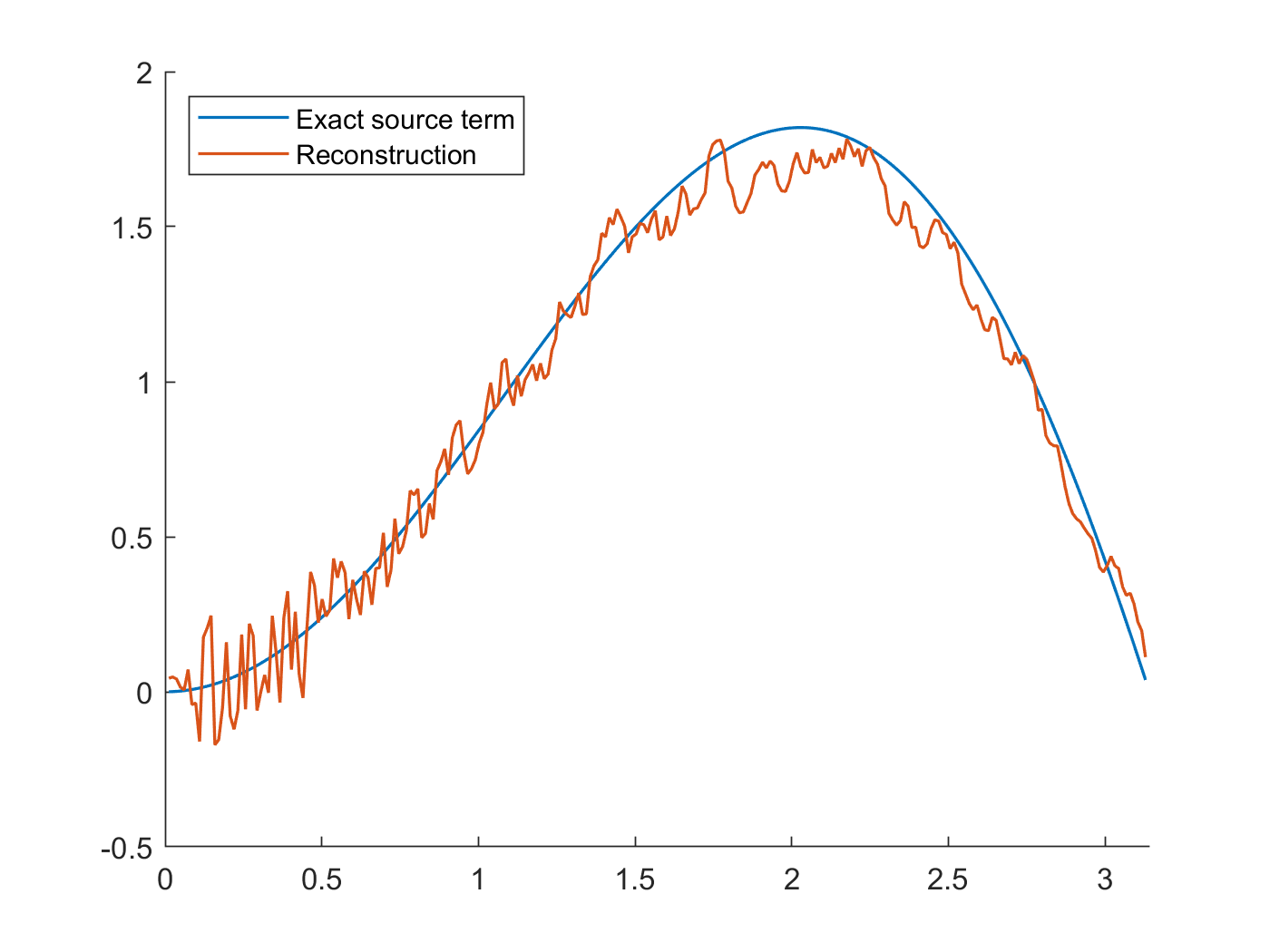}}
    \hfill
    \subfloat[][$\lambda=10^{-2}$, $\varepsilon=0.01$]
    {\includegraphics[width=0.33\textwidth]{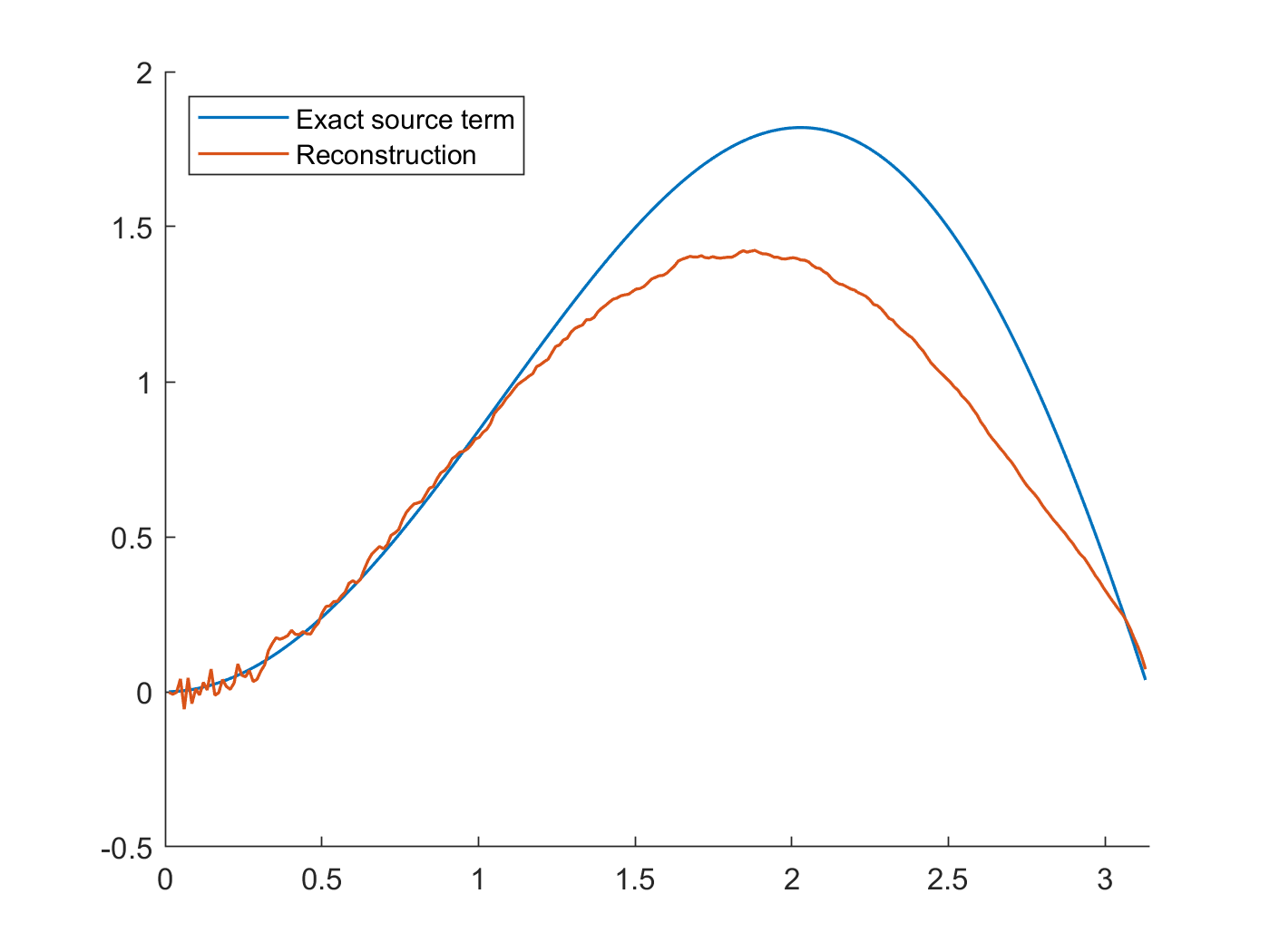}}
    \\
    \caption{1D case - Comparison between the exact source function $f(x)$, in blue, and its reconstruction, in red, for different choices of the regularization parameter $\lambda$ and noise level $\varepsilon$. The problem sizes are $n=2^8$, $S=2^6$ and the fractional orders are $\beta=0.9$, $\omega=1.1$.}
    \label{fig:1D_accuracy_beta0.9-omega1.1}
\end{figure}

Similar patterns are observed in Figures \ref{fig:1D_accuracy_beta0.5-omega1.5} and \ref{fig:1D_accuracy_beta0.9-omega1.1}, which depict the results for the fractional orders $\beta=0.5$, $\omega=1.5$ and $\beta=0.9$, $\omega=1.1$, respectively.

Furthermore, Tables \ref{table:1D_accuracy_beta0.1-omega1.9}--\ref{table:1D_accuracy_beta0.9-omega1.1} report the 2-norm relative error between the exact source function $f(x)$ and the reconstruction, computed as
\begin{equation*}
    \texttt{error} = \frac{\norm{\mathbf{f}-\tilde{\mathbf{f}}}_2}{\norm{\mathbf{f}}_2}.
\end{equation*}
Indeed, these quantitative results confirm the visual observations from the plots. For each noise level, the minimum error occurs at the intermediate value of $\lambda$ and, unsurprisingly, the error generally increases with the noise level.

\begin{table}[ht]
    \centering
    \caption{1D case - 2-norm relative error between the exact source function $f(x)$ and its reconstruction, for different choices of the regularization parameter $\lambda$ and noise level $\varepsilon$. The fractional orders are $\beta=0.1$, $\omega=1.9$.}
    \label{table:1D_accuracy_beta0.1-omega1.9}
    \begin{tabular}{cccccccc}
    \toprule
    & & \multicolumn{3}{c}{$\varepsilon=0.001$} & \multicolumn{3}{c}{$\varepsilon=0.01$} \\
    \cmidrule(lr){3-5}\cmidrule(lr){6-8}
    & & \multicolumn{3}{c}{$\lambda$} & \multicolumn{3}{c}{$\lambda$} \\
    \cmidrule(lr){3-5}\cmidrule(lr){6-8}
    $n$ & $S$ & $10^{-5}$ & $10^{-4}$ & $10^{-3}$ & $10^{-4}$ & $10^{-3}$ & $10^{-2}$ \\
    \midrule
    $2^{4}$ & $2^{4}$ & 0.02 & 0.01 & 0.05 & 0.11 & 0.07 & 0.26 \\
    $2^{4}$ & $2^{5}$ & 0.03 & 0.01 & 0.05 & 0.10 & 0.07 & 0.25 \\
    $2^{4}$ & $2^{6}$ & 0.03 & 0.01 & 0.05 & 0.12 & 0.05 & 0.25 \\
    \midrule
    $2^{5}$ & $2^{4}$ & 0.05 & 0.02 & 0.05 & 0.14 & 0.06 & 0.25 \\
    $2^{5}$ & $2^{5}$ & 0.05 & 0.02 & 0.05 & 0.12 & 0.05 & 0.25 \\
    $2^{5}$ & $2^{6}$ & 0.05 & 0.02 & 0.05 & 0.13 & 0.04 & 0.25 \\
    \midrule
    $2^{6}$ & $2^{4}$ & 0.07 & 0.02 & 0.05 & 0.13 & 0.06 & 0.25 \\
    $2^{6}$ & $2^{5}$ & 0.06 & 0.02 & 0.05 & 0.14 & 0.08 & 0.26 \\
    $2^{6}$ & $2^{6}$ & 0.06 & 0.01 & 0.05 & 0.15 & 0.07 & 0.25 \\
    \midrule
    $2^{7}$ & $2^{4}$ & 0.05 & 0.02 & 0.05 & 0.17 & 0.07 & 0.25 \\
    $2^{7}$ & $2^{5}$ & 0.06 & 0.02 & 0.05 & 0.17 & 0.07 & 0.26 \\
    $2^{7}$ & $2^{6}$ & 0.07 & 0.02 & 0.05 & 0.17 & 0.08 & 0.24 \\
    \midrule
    $2^{8}$ & $2^{4}$ & 0.06 & 0.02 & 0.05 & 0.14 & 0.07 & 0.25 \\
    $2^{8}$ & $2^{5}$ & 0.06 & 0.02 & 0.05 & 0.15 & 0.05 & 0.24 \\
    $2^{8}$ & $2^{6}$ & 0.06 & 0.02 & 0.05 & 0.14 & 0.05 & 0.25 \\
    \bottomrule
    \end{tabular}
\end{table}

\begin{table}[ht]
    \centering
    \caption{1D case - 2-norm relative error between the exact source function $f(x)$ and its reconstruction, for different choices of the regularization parameter $\lambda$ and noise level $\varepsilon$. The fractional orders are $\beta=0.5$, $\omega=1.5$.}
    \label{table:1D_accuracy_beta0.5-omega1.5}
    \begin{tabular}{cccccccc}
    \toprule
    & & \multicolumn{3}{c}{$\varepsilon=0.001$} & \multicolumn{3}{c}{$\varepsilon=0.01$} \\
    \cmidrule(lr){3-5}\cmidrule(lr){6-8}
    & & \multicolumn{3}{c}{$\lambda$} & \multicolumn{3}{c}{$\lambda$} \\
    \cmidrule(lr){3-5}\cmidrule(lr){6-8}
    $n$ & $S$ & $10^{-5}$ & $10^{-4}$ & $10^{-3}$ & $10^{-4}$ & $10^{-3}$ & $10^{-2}$ \\
    $2^{4}$ & $2^{4}$ & 0.01 & 0.01 & 0.04 & 0.07 & 0.05 & 0.22 \\
    $2^{4}$ & $2^{5}$ & 0.01 & 0.01 & 0.04 & 0.08 & 0.06 & 0.23 \\
    $2^{4}$ & $2^{6}$ & 0.01 & 0.01 & 0.04 & 0.11 & 0.04 & 0.22 \\
    \midrule
    $2^{5}$ & $2^{4}$ & 0.04 & 0.01 & 0.04 & 0.15 & 0.07 & 0.23 \\
    $2^{5}$ & $2^{5}$ & 0.02 & 0.02 & 0.04 & 0.15 & 0.07 & 0.21 \\
    $2^{5}$ & $2^{6}$ & 0.04 & 0.02 & 0.04 & 0.14 & 0.05 & 0.23 \\
    \midrule
    $2^{6}$ & $2^{4}$ & 0.07 & 0.02 & 0.04 & 0.19 & 0.06 & 0.22 \\
    $2^{6}$ & $2^{5}$ & 0.07 & 0.02 & 0.04 & 0.19 & 0.04 & 0.22  \\
    $2^{6}$ & $2^{6}$ & 0.07 & 0.02 & 0.04 & 0.16 & 0.05 & 0.23 \\
    \midrule
    $2^{7}$ & $2^{4}$ & 0.08 & 0.02 & 0.04 & 0.23 & 0.07 & 0.22 \\
    $2^{7}$ & $2^{5}$ & 0.09 & 0.02 & 0.04 & 0.21 & 0.06 & 0.21 \\
    $2^{7}$ & $2^{6}$ & 0.08 & 0.02 & 0.04 & 0.19 & 0.06 & 0.22 \\
    \midrule
    $2^{8}$ & $2^{4}$ & 0.09 & 0.02 & 0.04 & 0.22 & 0.05 & 0.22 \\
    $2^{8}$ & $2^{5}$ & 0.09 & 0.03 & 0.04 & 0.19 & 0.05 & 0.23 \\
    $2^{8}$ & $2^{6}$ & 0.10 & 0.02 & 0.04 & 0.20 & 0.07 & 0.22 \\
    \bottomrule
    \end{tabular}
\end{table}

\begin{table}[ht]
    \centering
    \caption{1D case - 2-norm relative error between the exact source function $f(x)$ and its reconstruction, for different choices of the regularization parameter $\lambda$ and noise level $\varepsilon$. The fractional orders are $\beta=0.9$, $\omega=1.1$.}
    \label{table:1D_accuracy_beta0.9-omega1.1}
    \begin{tabular}{cccccccc}
    \toprule
    & & \multicolumn{3}{c}{$\varepsilon=0.001$} & \multicolumn{3}{c}{$\varepsilon=0.01$} \\
    \cmidrule(lr){3-5}\cmidrule(lr){6-8}
    & & \multicolumn{3}{c}{$\lambda$} & \multicolumn{3}{c}{$\lambda$} \\
    \cmidrule(lr){3-5}\cmidrule(lr){6-8}
    $n$ & $S$ & $10^{-5}$ & $10^{-4}$ & $10^{-3}$ & $10^{-4}$ & $10^{-3}$ & $10^{-2}$ \\
    \midrule
    $2^{4}$ & $2^{4}$ & 0.01 & 0.01 & 0.04 & 0.06 & 0.05 & 0.23 \\
    $2^{4}$ & $2^{5}$ & 0.00 & 0.01 & 0.04 & 0.04 & 0.06 & 0.24 \\
    $2^{4}$ & $2^{6}$ & 0.01 & 0.00 & 0.04 & 0.04 & 0.05 & 0.24 \\
    \midrule
    $2^{5}$ & $2^{4}$ & 0.01 & 0.01 & 0.04 & 0.10 & 0.06 & 0.23 \\
    $2^{5}$ & $2^{5}$ & 0.01 & 0.01 & 0.04 & 0.10 & 0.05 & 0.23 \\
    $2^{5}$ & $2^{6}$ & 0.01 & 0.01 & 0.04 & 0.08 & 0.06 & 0.23 \\
    \midrule
    $2^{6}$ & $2^{4}$ & 0.03 & 0.02 & 0.04 & 0.16 & 0.06 & 0.24 \\
    $2^{6}$ & $2^{5}$ & 0.03 & 0.02 & 0.04 & 0.15 & 0.08 & 0.23 \\
    $2^{6}$ & $2^{6}$ & 0.03 & 0.02 & 0.04 & 0.17 & 0.06 & 0.23 \\
    \midrule
    $2^{7}$ & $2^{4}$ & 0.07 & 0.02 & 0.04 & 0.25 & 0.08 & 0.22 \\
    $2^{7}$ & $2^{5}$ & 0.06 & 0.03 & 0.04 & 0.23 & 0.09 & 0.23 \\
    $2^{7}$ & $2^{6}$ & 0.07 & 0.02 & 0.04 & 0.26 & 0.08 & 0.24 \\
    \midrule
    $2^{8}$ & $2^{4}$ & 0.10 & 0.03 & 0.04 & 0.31 & 0.09 & 0.22 \\
    $2^{8}$ & $2^{5}$ & 0.11 & 0.03 & 0.04 & 0.30 & 0.07 & 0.22 \\
    $2^{8}$ & $2^{6}$ & 0.11 & 0.03 & 0.04 & 0.28 & 0.07 & 0.23 \\
    \bottomrule
    \end{tabular}
\end{table}

\subsection{Two dimensional setting}

Now we proceed to the two-dimensional version of \eqref{eq:regularized-fde}, setting
\begin{gather*}
    \Omega=(0,\pi)^2, \qquad T=1, \\
    \gamma_1(x_1,x_2,t) = t e^{x_1+x_2}, \qquad \gamma_2(x_1,x_2,t) = t^2 x_1x_2, \\
    q(t) = t^2, \qquad \rho(x_1,x_2)=0.
\end{gather*}
We construct the final time data $\mu(x_1,x_2)=v(x_1,x_2,T)$ by fixing
\begin{equation*}
   f(x_1,x_2) = x_1 x_2 \sin (x_1x_2)
\end{equation*}
and numerically solving the corresponding direct problem, as we have done in Subsection \ref{sssec:numer-1d}. Likewise, we simulate the noisy data by adding a random quantity $\varepsilon\delta (x_1,x_2)$ to $\mu(x_1,x_2)$, with $\delta (x_1,x_2)$ randomly selected in $(-1,1)^2$ and $\varepsilon = 0.01\cdot\norm{\mu}_2$. Moreover, we always set $n_1=n_2=n$ for simplicity, so that $\mi{n} =(n,n)$.

\subsubsection{Spectral analysis}

We examine the effectiveness of $P_{\mi{n},S}$ in reducing the condition number and clustering the eigenvalues of the coefficient matrix at 1. The regularization parameter is set to $\lambda= 5\cdot10^{-3}$.

\begin{figure}[ht]
    \centering
    \subfloat[][$\beta=0.1$, $\omega=1.9$]
    {\includegraphics[width=0.33\textwidth]{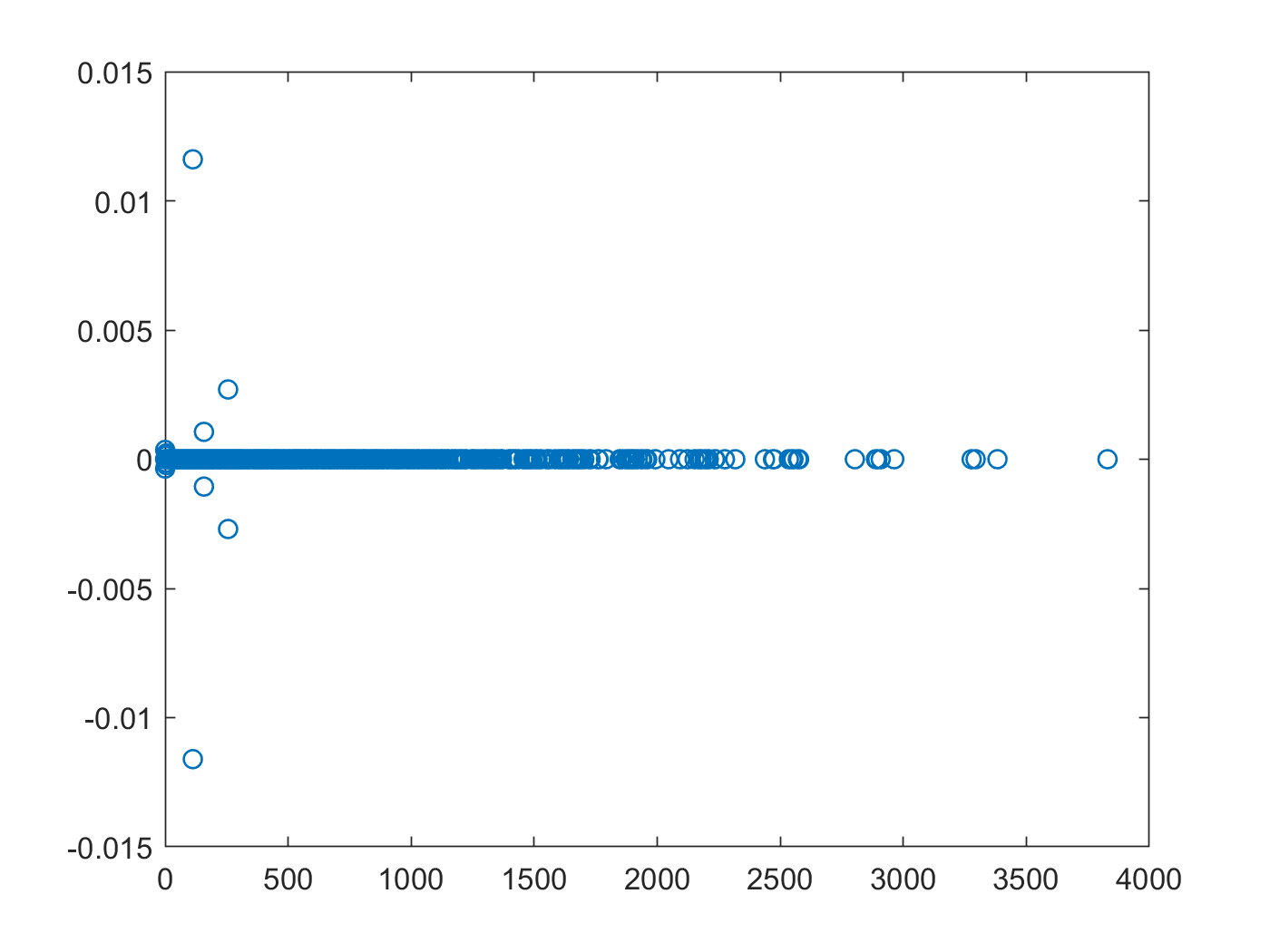}}
    \hfill
    \subfloat[][$\beta=0.1$, $\omega=1.9$]
    {\includegraphics[width=0.33\textwidth]{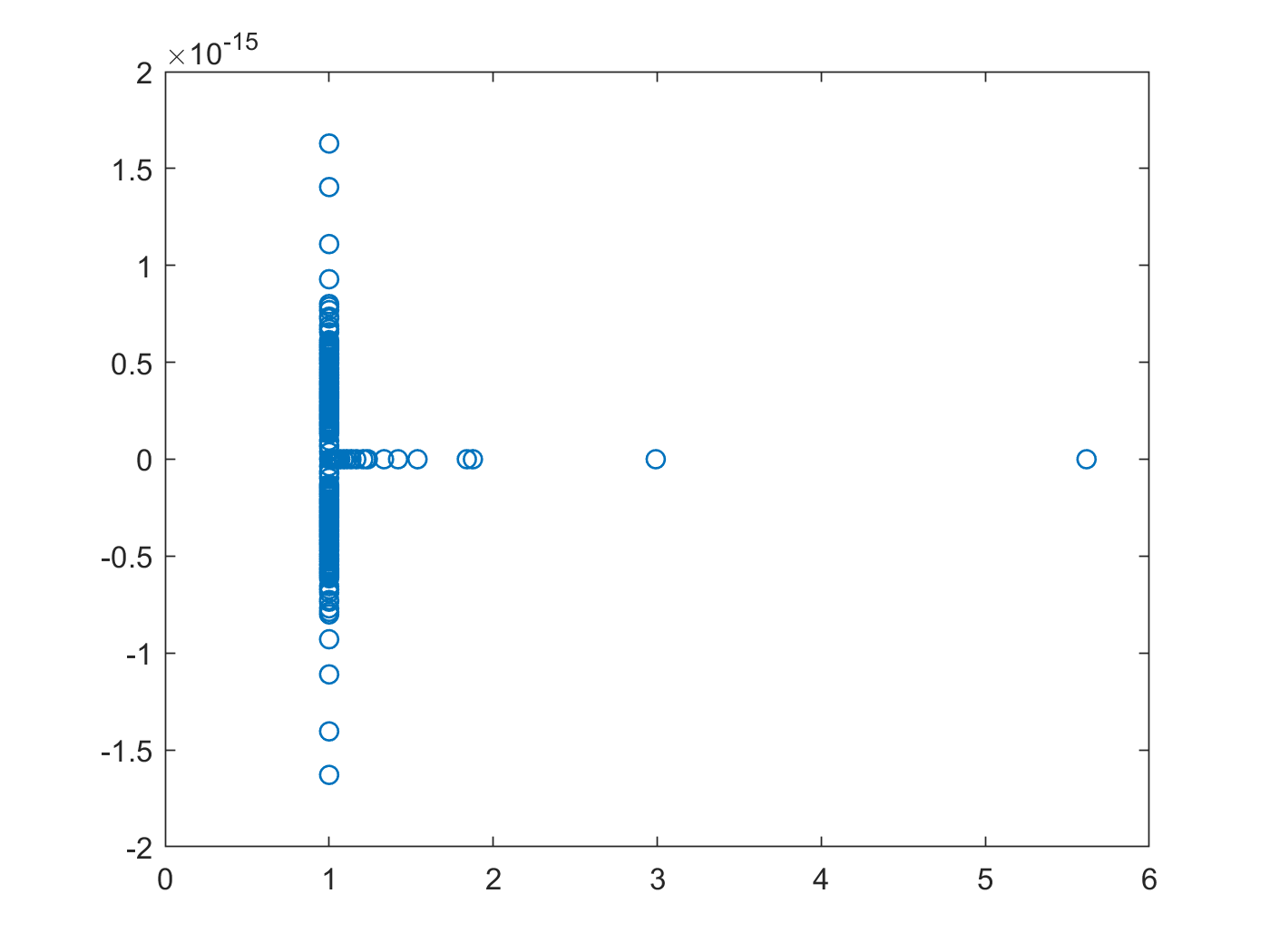}}
    \hfill
    \subfloat[][$\beta=0.1$, $\omega=1.9$]
    {\includegraphics[width=0.33\textwidth]{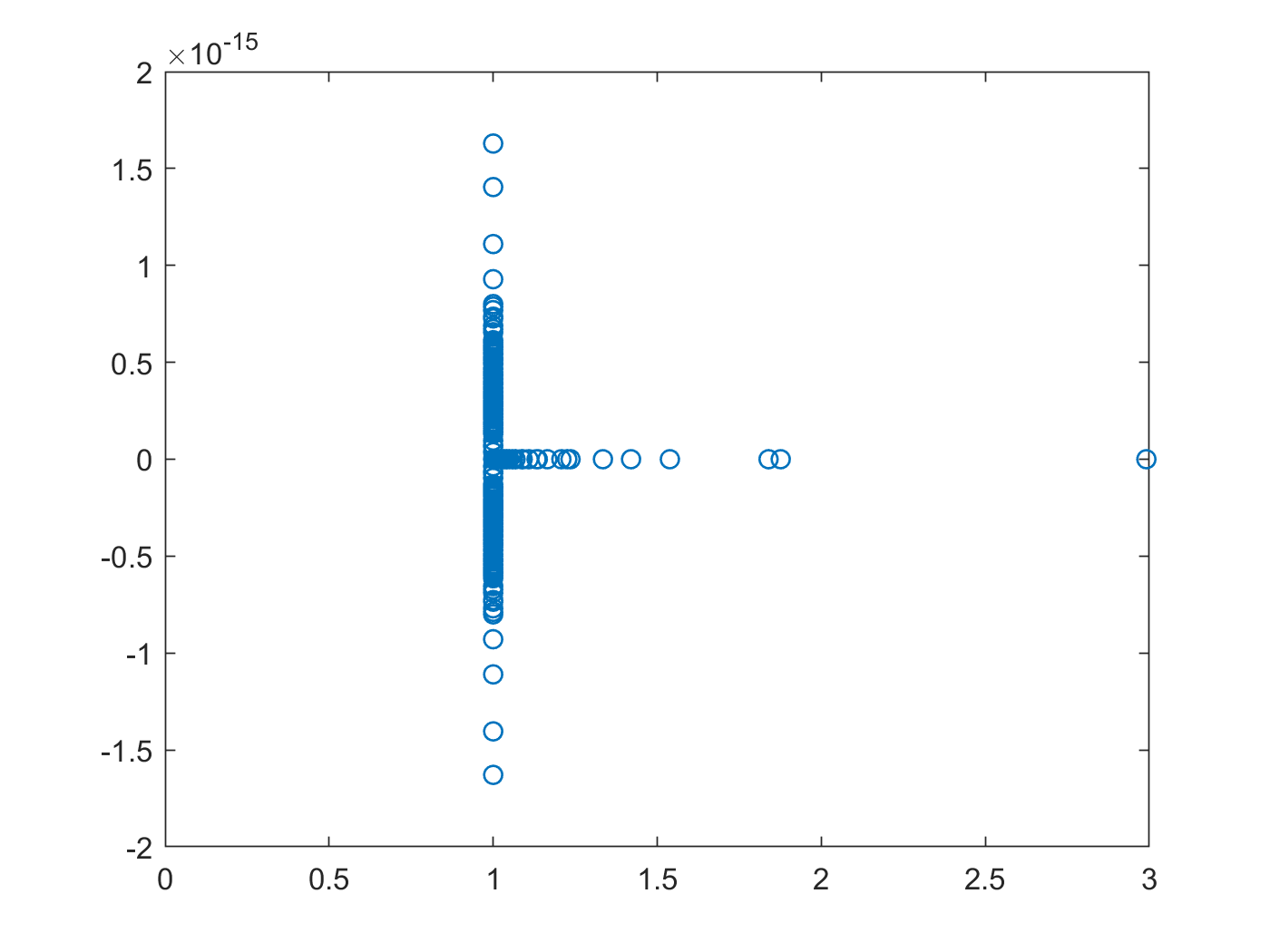}}
    \\
    \subfloat[][$\beta=0.5$, $\omega=1.5$]
    {\includegraphics[width=0.33\textwidth]{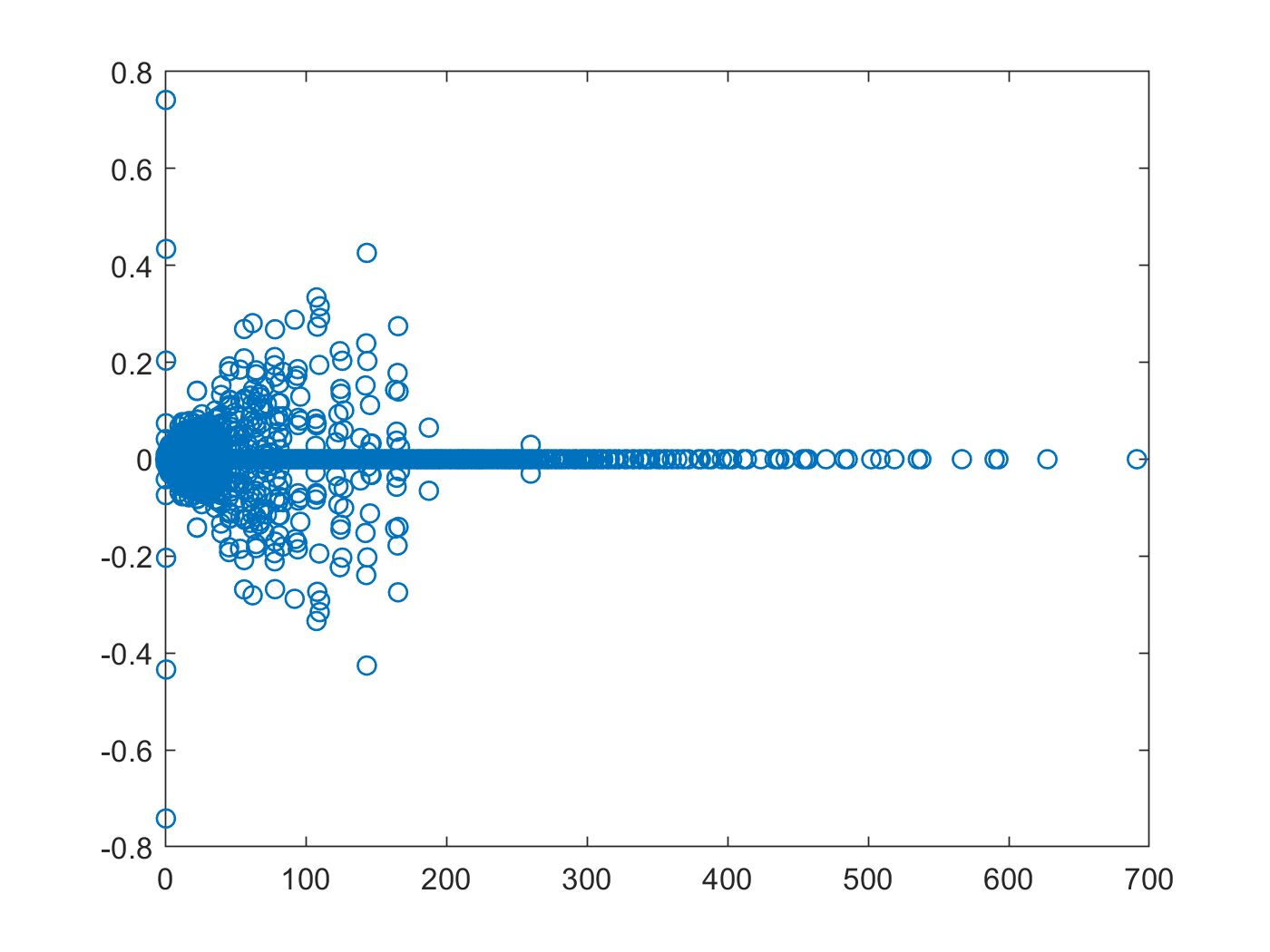}}
    \hfill
    \subfloat[][$\beta=0.5$, $\omega=1.5$]
    {\includegraphics[width=0.33\textwidth]{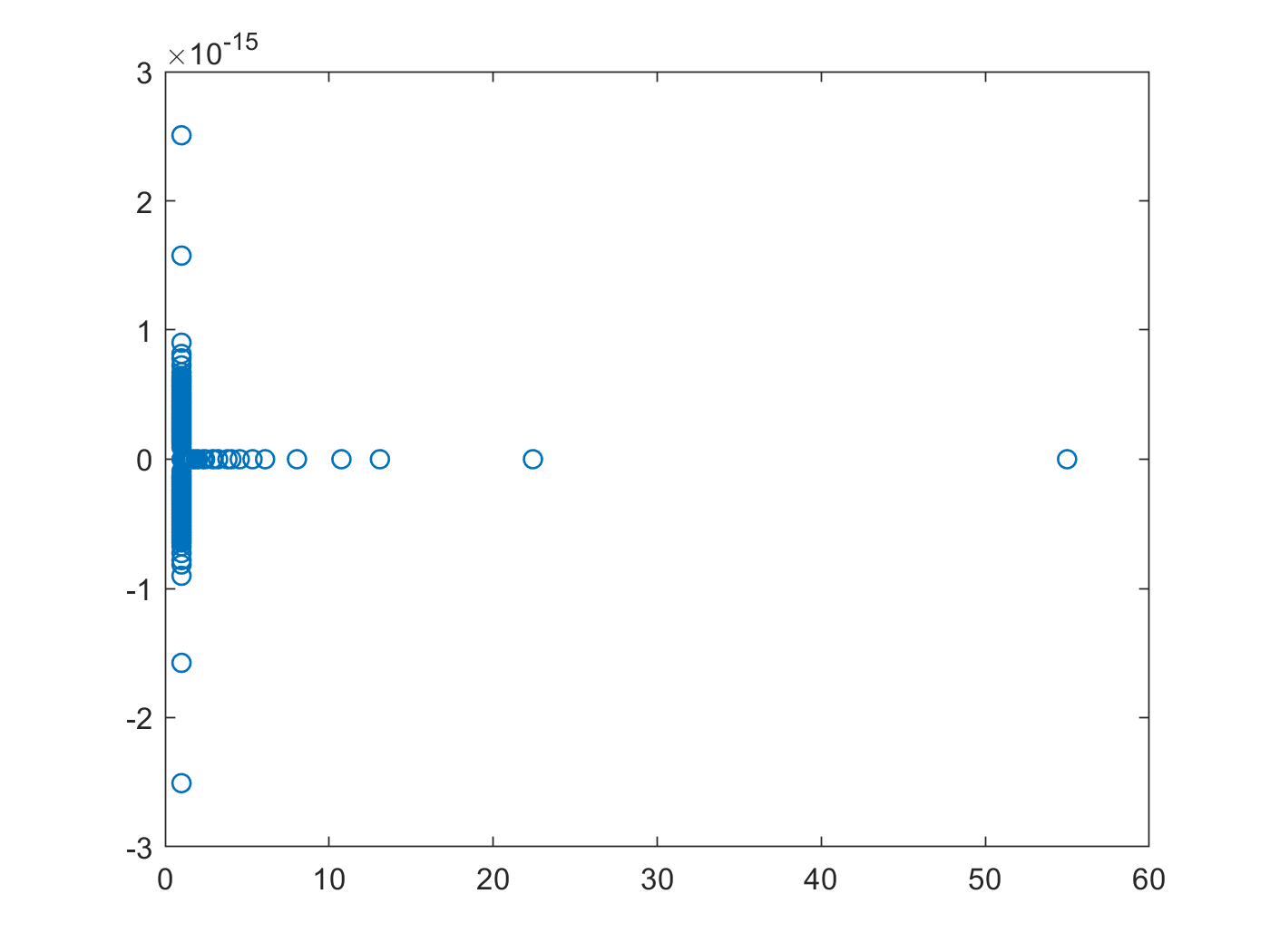}}
    \hfill
    \subfloat[][$\beta=0.5$, $\omega=1.5$]
    {\includegraphics[width=0.33\textwidth]{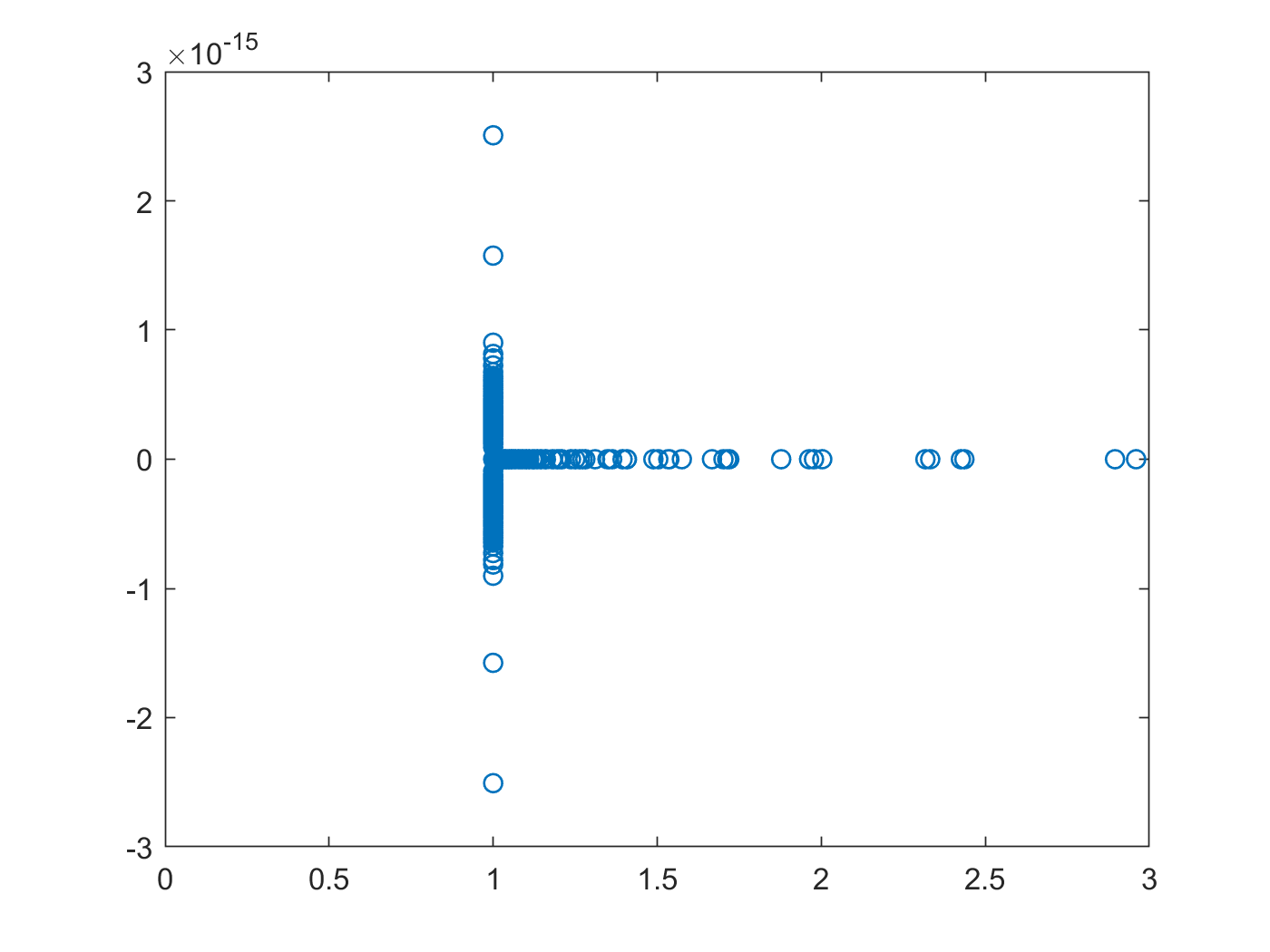}}
    \\
    \subfloat[][$\beta=0.9$, $\omega=1.1$]
    {\includegraphics[width=0.33\textwidth]{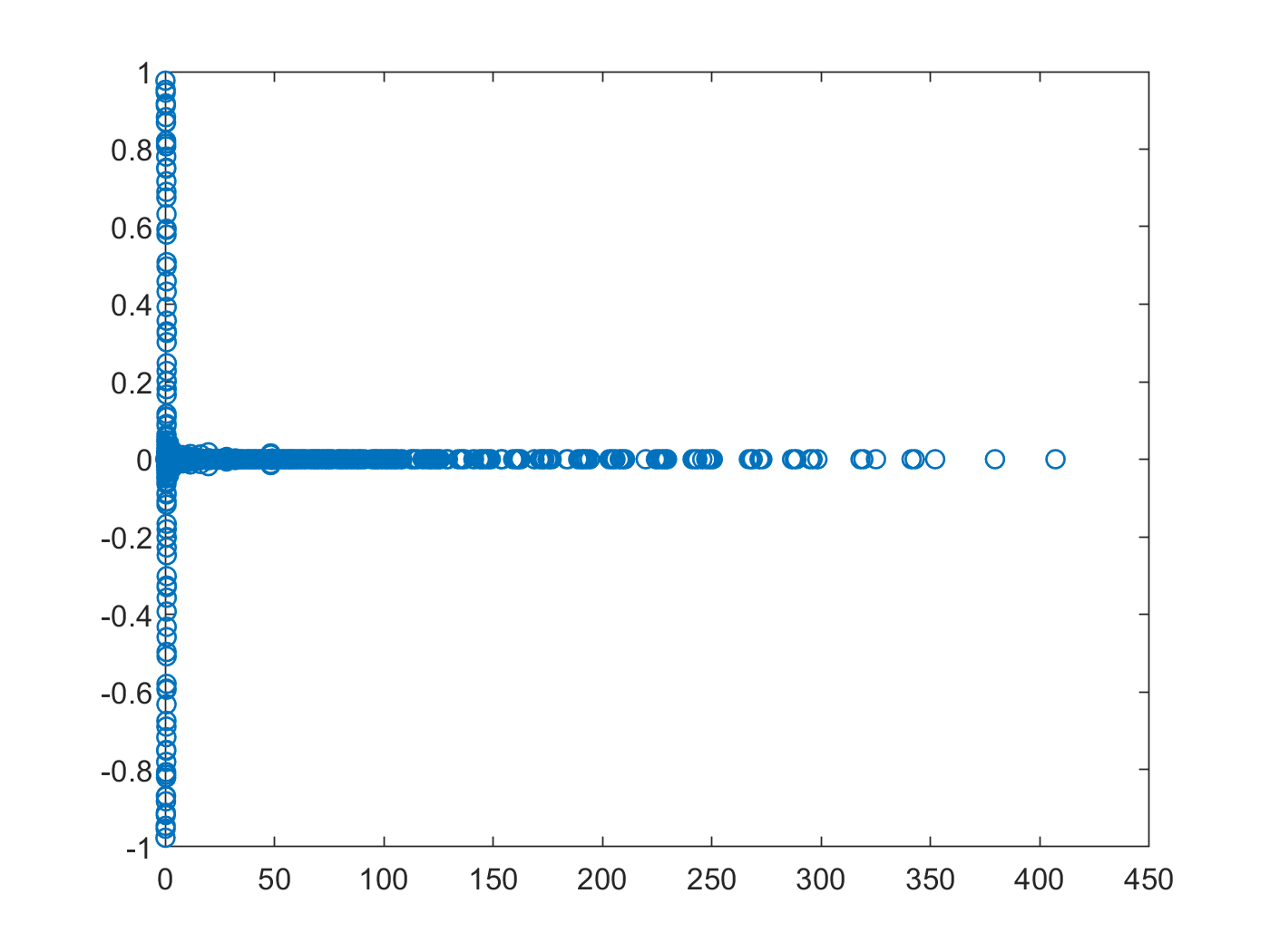}}
    \hfill
    \subfloat[][$\beta=0.9$, $\omega=1.1$]
    {\includegraphics[width=0.33\textwidth]{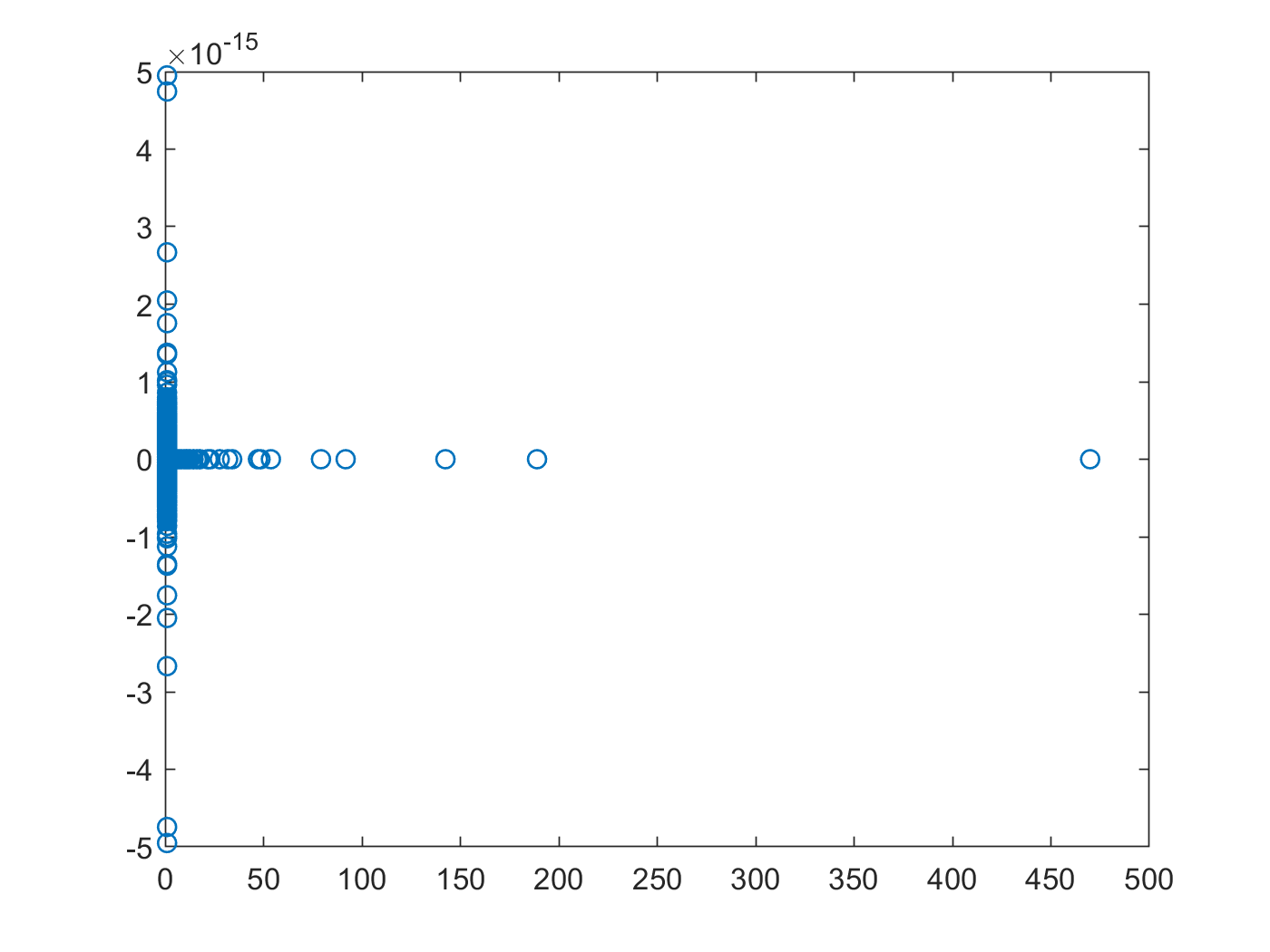}}
    \subfloat[][$\beta=0.9$, $\omega=1.1$]
    {\includegraphics[width=0.33\textwidth]{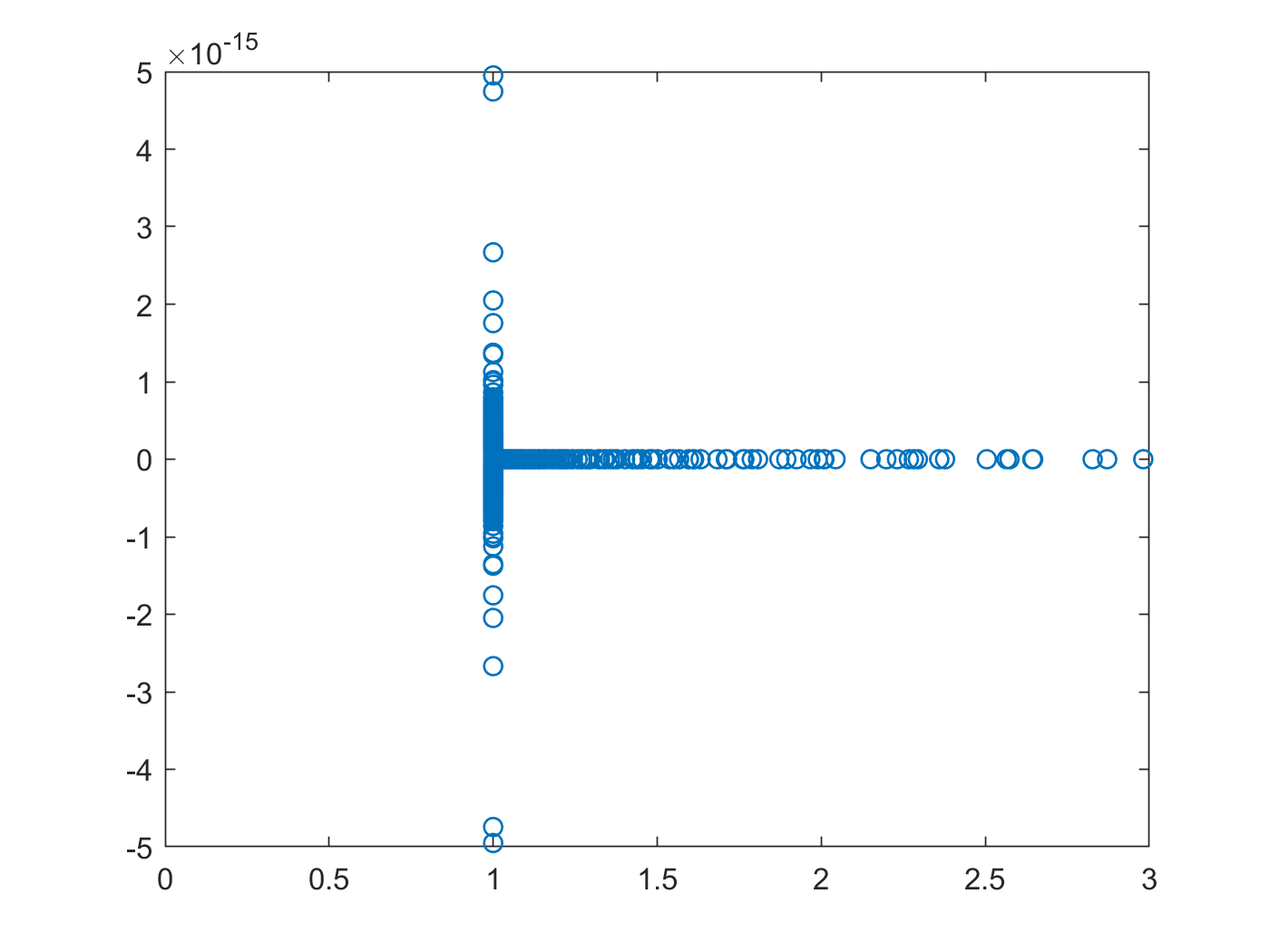}}
    \\
    \caption{2D case - Eigenvalues of $A_{\mi{n},S}$ (left column) and $P_{\mi{n},S}^{-1}A_{\mi{n},S}$ (central column) in the complex plane, for different choices of the fractional orders, with $\mi{n}=(2^4,2^4)$ and $S=2^4$. For clarity, we also display magnified versions of the central graphs, focused on the cluster at 1 (right column).}
    \label{fig:2D_eigenvalues_n4s4}
\end{figure}

Figure \ref{fig:2D_eigenvalues_n4s4} shows the eigenvalues in the complex plane for various choices of $\beta$ and $\omega$, with fixed problem dimensions $n=2^4$ (so that $\mi{n}=(2^4,2^4)$) and $S=2^4$. The plots are consistent with the theoretical discussion, showing that the eigenvalues accumulate near 0 for $A_{n,S}$. In the second and especially the third column, containing magnified graphs like in the one dimensional case, a tight cluster around 1 is clearly visible for the preconditioned matrix, with no eigenvalues near zero across all parameter settings.

Table \ref{table:2D_condnum} contains the 2-norm condition numbers of the nonpreconditioned and preconditioned coefficient matrix, for different choices of $\mi{n}$ and $S$. Once again, $A_{\mi{n},S}$ is extremely ill-conditioned, even for such small problem sizes. Similarly to the one dimensional case, $P_{\mi{n},S}$ is effective in reducing the condition number and performs very well when $A_{n,S}$ is particularly ill-conditioned, for $\beta\to 1$ and $\omega\to 2$.

\begin{table}[ht]
    \centering
    \caption{2D case - 2-norm condition numbers of the nonpreconditioned and preconditioned coefficient matrix.}
    \label{table:2D_condnum}
    \begin{tabular}{cccccccc}
    \toprule
    & & \multicolumn{2}{c}{$\beta=0.1$, $\omega=1.9$} & \multicolumn{2}{c}{$\beta=0.5$, $\omega=1.5$} & \multicolumn{2}{c}{$\beta=0.9$, $\omega=1.1$} \\
    \cmidrule(lr){3-4}\cmidrule(lr){5-6}\cmidrule(lr){7-8}
    $\mi{n}$ & $S$ & - & $P_{\mi{n},S}$ & - & $P_{\mi{n},S}$ & - & $P_{\mi{n},S}$ \\
    \midrule
    $2^{4}$ & $2^{4}$ & 34428 & 7 & 11123 & 65 & 18795 & 538 \\
    $2^{4}$ & $2^{5}$ & 69177 & 8 & 23296 & 98 & 54246 & 1153 \\
    $2^{5}$ & $2^{4}$ & 149171 & 9 & 27400 & 124 & 28649 & 1332 \\
    $2^{5}$ & $2^{5}$ & 315653 & 10 & 51158 & 187 & 80995 & 2922 \\
    \bottomrule
    \end{tabular}
\end{table}

\subsubsection{GMRES method performance}

We proceed to solving the linear system \eqref{eq:lin-syst} with the GMRES method, reporting the results in Table \ref{table:2D_gmres}. When the maximum number of iterations, which is the minimum between the size of the matrix and 1000, is reached before convergence, we simply write $>1000$. The initial guess is set to the zero vector, the tolerance is $\texttt{tol}=10^{-8}$. The regularization parameter is set to $\lambda= 5\cdot10^{-3}$ and the noise level to $\varepsilon=0.01$.

In the columns related to the nonpreconditioned linear system, we observe consistently high iteration counts for all values of $\mi{n}$ and $S$, in accordance with the severe ill-conditioning. In fact, the number of iterations quickly becomes unmanageable, even using a coarse discretization.

In contrast, the columns associated with the preconditioned system show iteration counts that are consistent with the spectral analysis. $P_{\mi{n},S}$ delivers an excellent performance, particularly when $\beta$ approaches 1 and $\omega$ approaches 2, where the number of iterations remains very moderate and nearly independent from $\mi{n}$ and $S$. It is less effective when $\beta\to 0$ and $\omega\to 1$, in accordance with Table \ref{table:2D_condnum}, which shows higher condition numbers for such values of the fractional orders.

\begin{table}[ht]
    \centering
    \caption{2D case - Iterations required for solving the nonpreconditioned and preconditioned linear system with the GMRES method.}
    \label{table:2D_gmres}
    \begin{tabular}{cccccccc}
    \toprule
    & & \multicolumn{2}{c}{$\beta=0.1$, $\omega=1.9$} & \multicolumn{2}{c}{$\beta=0.5$, $\omega=1.5$} & \multicolumn{2}{c}{$\beta=0.9$, $\omega=1.1$} \\
    \cmidrule(lr){3-4}\cmidrule(lr){5-6}\cmidrule(lr){7-8}
    $\mi{n}$ & $S$ & - & $P_{\mi{n},S}$ & - & $P_{\mi{n},S}$ & - & $P_{\mi{n},S}$ \\
    \midrule
    $2^{4}$ & $2^{4}$ & 867 & 10 & 460 & 20 & 486 & 41 \\
    $2^{4}$ & $2^{5}$ & $>$1000 & 10 & 597 & 22 & 716 & 48 \\
    $2^{4}$ & $2^{6}$ & $>$1000 & 10 & 795 & 25 & $>$1000 & 58 \\
    $2^{4}$ & $2^{7}$ & $>$1000 & 11 & $>$1000 & 27 & $>$1000 & 68 \\
    $2^{4}$ & $2^{8}$ & $>$1000 & 11 & $>$1000 & 30 & $>$1000 & 81 \\
    \midrule
    $2^{5}$ & $2^{4}$ & $>$1000 & 11 & 861 & 25 & 760 & 60 \\
    $2^{5}$ & $2^{5}$ & $>$1000 & 11 & $>$1000 & 29 & $>$1000 & 72\\
    $2^{5}$ & $2^{6}$ & $>$1000 & 12 & $>$1000 & 32 & $>$1000 & 87 \\
    $2^{5}$ & $2^{7}$ & $>$1000 & 12 & $>$1000 & 35 & $>$1000 & 104 \\
    $2^{5}$ & $2^{8}$ & $>$1000 & 13 & $>$1000 & 39 & $>$1000 & 124 \\
    \midrule
    $2^{6}$ & $2^{4}$ & $>$1000 & 12 & $>$1000 & 32 & $>$1000 & 87 \\
    $2^{6}$ & $2^{5}$ & $>$1000 & 13 & $>$1000 & 36 & $>$1000 & 105 \\
    $2^{6}$ & $2^{6}$ & $>$1000 & 13 & $>$1000 & 40 & $>$1000 & 127 \\
    \bottomrule
    \end{tabular}
\end{table}

\section{Conclusions} \label{sec:concl}

We investigated the effectiveness of block triangular preconditioners for accelerating and stabilizing the numerical solution of inverse source problems governed by time-space fractional diffusion equations (TSFDEs). Specifically, we addressed the recovery of an unknown spatial source function in a multi-dimensional TSFDE incorporating Caputo time-fractional derivatives and the fractional Laplacian. To tackle the inherent ill-posedness of the problem, we employed a quasi-boundary value regularization approach, followed by a finite difference discretization. The resulting linear systems are characterized by various levels of structure, which we exploited to design and analyze a block triangular preconditioning strategy. Numerical experiments using the GMRES solver shown that the proposed preconditioner significantly improves convergence rates, robustness, and accuracy, making it well-suited for large-scale, real world inverse problems involving fractional models.

As future research directions, we identify two promising open questions:
\begin{enumerate}[a) ]
    \item an analysis of the eigenvalue and singular value distributions in the Weyl sense of the resulting matrix-sequences, both original and preconditioned, using the framework of Generalized Locally Toeplitz (GLT) sequences \cite{book1,book2} (see \cite{tutorial} for a gentle guide on the GLT theory in several practical discretizations of differentiaal operators);
    \item the study of item a) within the context of block-structured matrix analysis, as developed in \cite{block1,block2}, so allowing a great degree of generality in the approximation schemes.
\end{enumerate}
Pursuing these directions may lead to the development of even more effective preconditioning techniques and provide a deeper understanding of the spectral properties of the matrix sequences arising from such fractional inverse problems.

\section*{Acknowledgments}
The research of Stefano Serra-Capizzano is supported by the PRIN-PNRR project \lq\lq MATH-ematical tools for predictive maintenance and PROtection of CULTtural heritage (MATHPROCULT)'' (code P20228HZWR, CUP J53D23003780006), by INdAM-GNCS Project \lq\lq Analisi e applicazioni di matrici strutturate (a blocchi)'' (CUP E53C23001670001), and by the European High-Performance Computing Joint Undertaking (JU) under Grant Agreement 955701. The JU receives support from the European Union’s Horizon 2020 research and innovation programme and Belgium, France, Germany, Switzerland. Furthermore Stefano Serra-Capizzano is grateful for the support of the Laboratory of Theory, Economics and Systems – Department of Computer Science at Athens University of Economics and Business. The research of Rosita L. Sormani is funded by the PRIN-PNRR project \lq\lq A mathematical approach to
inverse problems arising in cultural heritage preservation and dissemination'' (code P2022PMEN2, CUP F53D23010100001). Finally Stefano Serra-Capizzano and Rosita L. Sormani are partly supported by Italian National Agency INdAM-GNCS.


\begin{thebibliography}{99}

\bibitem{block1}
Adriani A., Schiavoni-Piazza A.J.A., Serra-Capizzano S.,
\emph{Blocking structures, g.a.c.s. approximation, and distributions},
Bol. Soc. Mat. Mex. (3)  31(2), 41 (2025).

\bibitem{block2}
Barakitis N., Donatelli M., Ferri S., Loi V., Serra-Capizzano S., Sormani R.L., \emph{
Block structures, approximation, and preconditioning}, Numer. Alg. (2025), online 8-7-25,
\\
https://link.springer.com/article/10.1007/s11075-025-02157-y.



    \bibitem{Bini-tau}
    Bini D., Di Benedetto F., \emph{A new preconditioner for the parallel solution of positive definite Toeplitz systems}, Proceedings of the second annual ACM symposium on Parallel algorithms and architectures (SPAA '90), New York, 220--223 (1990).


    \bibitem{chen2018}
    Chen H., Zhang T., Lv W., \emph{Block preconditioning strategies for time-space fractional diffusion equations}, Appl. Math. Comput. 337, 41--53 (2018).

% \bibitem{tau-better1}
% Di Benedetto, F. \emph{Analysis of preconditioning techniques for ill-conditioned Toeplitz matrices}, SIAM J. Sci. Comput. 16(3), 682-–697 (1995).

\bibitem{book1}
Garoni C., Serra-Capizzano S.,
\emph{The theory of Generalized Locally Toeplitz sequences:  theory and applications - Vol I.}
SPRINGER - Springer Monographs in Mathematics,  Berlin, (2017).
%\\  http://www.springer.com/gp/book/9783319536781


\bibitem{book2}
Garoni C., Serra-Capizzano S.,
\emph{The theory of Generalized Locally Toeplitz sequences: theory and applications - Vol II.}
SPRINGER - Springer Monographs in Mathematics,  Berlin, (2018).

\bibitem{tutorial}
Garoni C., Serra-Capizzano S.,
\emph{Generalized locally Toeplitz sequences: a spectral analysis tool for discretized differential equations.}
Splines and PDEs: from approximation theory to numerical linear algebra, 161–236, Lecture Notes in Math., 2219, Fond. CIME/CIME Found. Subser., Springer, Cham, 2018.

SPRINGER - Springer Monographs in Mathematics,  Berlin, (2018).

    \bibitem{diff-lapl}
    Hao Z., Zhang Z., Du R., \emph{Fractional centered difference scheme for high-dimensional integral fractional Laplacian}, J. Comput. Phys. 424, 109851 (2021).

% \bibitem{tau-appl2}
% Hon, S., Fung, P.Y., Dong, J., Serra-Capizzano, S. \emph{A sine transform based preconditioned MINRES method for all-at-once systems from constant and variable-coefficient evolutionary PDEs}, Numer. Algorithms 95(4), 1769–-1799 (2024).

% \bibitem{tau-appl3}
% Hon, S., Dong, J., Serra-Capizzano, S. \emph{A preconditioned MINRES method for optimal control of wave equations and its asymptotic spectral distribution theory}, SIAM J. Matrix Anal. Appl. 44(4) 1477-–1509 (2023).

% \bibitem{tau-appl1}
% Lin, X.-L., Ng, M.K, \emph{A $\tau$-preconditioner for space fractional diffusion equation with non-separable variable coefficients}, J. Sci. Comput. 100(1), 28 (2024).

    \bibitem{luo2022}
    Luo J.-M., Li H.-B., Wei W.-B., \emph{Block splitting preconditioner for time-space fractional diffusion equations}, Electron. Res. Arch. 30(3), 780--797 (2022).


    \bibitem{mariarosa24}
    Mazza M., Serra-Capizzano S., Sormani R.L., \emph{Algebra preconditionings for 2D Riesz distributed-order space-fractional diffusion equations on convex domains}, Numer. Linear Algebra Appl. 31(3), e2536 (2024).


    % \bibitem{ng2004iterative}
    % Ng M.K., \emph{Iterative methods for Toeplitz systems}, Oxford University Press, New York (2004).


    \bibitem{pang2024}
    Pang, H.-K., Qin, H.-H. Ni, S., \emph{Sine transform based preconditioning for an inverse source problem of time-space fractional diffusion equations}, J. Sci. Comput. 100, 74 (2024).


    \bibitem{qiao25}
    Qiao, Y., Xiong, X., Han, J., \emph{A variational approach to recover the unknown source and initial condition for a time-space fractional diffusion equation}, J. Appl. Math. Comput. 71, 3445--3476 (2025).


    \bibitem{ruan2025}
    Ruan Z., Wan G., Zhang W., \emph{Reconstruction of a space-dependent source term for a time fractional diffusion equation by a modified quasi-boundary value regularization method}, Taiwanese J. Math. 29(3), 467--489 (2025).


    \bibitem{ruan2018}
    Ruan, Z., Zhang, S., Xiong, S., \emph{Solving an inverse source problem for a time fractional diffusion equation by a modified quasi-boundary value method}, Evol. Equ. Control Theory 7(4), 669--682 (2018).

    % \bibitem{tau-better3}
    % Serra-Capizzano S., \emph{Toeplitz preconditioners constructed from linear approximation processes},
    % SIAM J. Matrix Anal. Appl. 20(2), 446--465 (1999).

    % \bibitem{tau-better2}
    % Serra-Capizzano S., \emph{Superlinear PCG methods for symmetric Toeplitz systems}, Math. Comp. 68(226), 793--803 (1999).


    \bibitem{fde-diff-book}
    Sun Z.-Z., Gao G.-H., \emph{Fractional differential equations: finite difference methods}, De Gruyter, Berlin, Science Press Beijing (2020).


    % \bibitem{tian2015}
    % Tian W., Zhou H., Deng W., \emph{A class of second order difference approximations for solving space fractional diffusion equations}, Math. Comp. 84(294), 1703--1727 (2015).

    \bibitem{vanloan-book}
    Van Loan C., \emph{Computational frameworks for the fast Fourier transform}, SIAM, Philadelphia (PA) (1992).

\end{thebibliography}
\end{document}